\documentclass{article}
\usepackage{latexsym,amsfonts,amsmath,amsthm,makeidx}

\makeindex
\newtheorem{definition}{Definition}[section]
\newtheorem{theorem}{Theorem}[section]
\newtheorem{lemma}{Lemma}[section]
\newtheorem{corollary}{Corollary}[section]

\newtheorem{remark}{Remark}[section]
\newcommand{\RN}{\mathbb R^N}
\newcommand{\Om}{\Omega}
\newcommand{\om}{\omega}

\newcommand{\iy}{\infty}
\newcommand{\q}{\theta}
\newcommand{\qu}{\quad}
\newcommand{\s}{\section}
\newcommand{\dd}{\delta}
\newcommand{\G}{\Gamma}

\newcommand{\la}{\lambda}
\newcommand{\La}{\Lambda}

\newcommand{\R}{\mathbb R}
\newcommand{\al}{\alpha}

\newcommand{\wu}{\widetilde u}
\newcommand{\wv}{\widetilde v}
\newcommand{\wphi}{\widetilde \phi}
\newcommand{\ou}{\overline u}
\newcommand{\ov}{\overline v}
\newcommand{\oJ}{\overline J}

\newcommand{\bb}{\beta}

\newcommand{\e}{\varepsilon}
\newcommand{\vp}{\varphi}

\newcommand{\bt}{\begin{theorem}}
\newcommand{\et}{\end{theorem}}
\newcommand{\bl}{\begin{lemma}}
\newcommand{\el}{\end{lemma}}
\newcommand{\bd}{\begin{definition}}
\newcommand{\ed}{\end{definition}}
\newcommand{\bc}{\begin{corollary}}
\newcommand{\ec}{\end{corollary}}
\newcommand{\bp}{\begin{proof}}
\newcommand{\ep}{\end{proof}}
\newcommand{\bx}{\begin{example}}
\newcommand{\ex}{\end{example}}
\newcommand{\bi}{\begin{exercise}}
\newcommand{\ei}{\end{exercise}}
\newcommand{\bo}{\begin{prop}}
\newcommand{\eo}{\end{prop}}
\newcommand{\br}{\begin{remark}}
\newcommand{\er}{\end{remark}}
\newcommand{\be}{\begin{equation}}
\newcommand{\ee}{\end{equation}}
\newcommand{\ba}{\begin{align}}
\newcommand{\ea}{\end{align}}
\newcommand{\bn}{\begin{enumerate}}
\newcommand{\en}{\end{enumerate}}
\newcommand{\bg}{\begin{align*}}
\newcommand{\bcs}{\begin{cases}}
\newcommand{\ecs}{\end{cases}}

\newcommand{\M}{{\mathcal M}}
\newcommand{\N}{{\mathcal N}}

\newcommand{\BD}{{\mathbb D}}

\newcommand{\BS}{{\mathbb S}}

\newcommand{\intR}[1]{\int_{\R^N}#1\,  dx}
\newcommand{\Sg}{\Sigma}

\newcommand{\ga}{\gamma}
\newcommand{\sg}{\sigma}
\newcommand{\bean}{\begin{eqnarray*}}
\newcommand{\eean}{\end{eqnarray*}}
\newcommand{\loc}{\operatorname{\rm loc}}

\numberwithin{equation}{section}

\begin{document}

\title{\bf Existence and symmetry of positive ground states for a doubly critical Schr\"{o}dinger system\thanks{Supported by NSFC (11025106).  E-mail address:
chenzhijie1987@sina.com(Chen);\quad \quad wzou@math.tsinghua.edu.cn
(Zou)}}
\date{}
\author{{\bf Zhijie Chen, Wenming Zou}\\
\footnotesize {\it  Department of Mathematical Sciences, Tsinghua
University,}\\
\footnotesize {\it Beijing 100084, China} }

\maketitle \vskip0.16in
\begin{center}
\begin{minipage}{120mm}
\begin{center}{\bf Abstract}\end{center}

We study the following doubly critical Schr\"{o}dinger system
$$\begin{cases}-\Delta u -\frac{\la_1}{|x|^2}u=u^{2^\ast-1}+
\nu \al u^{\al-1}v^\bb, \quad x\in \RN,\\
  -\Delta v -\frac{\la_2}{|x|^2}v=v^{2^\ast-1} +
\nu \bb u^{\al}v^{\bb-1},     \quad x\in \RN,\\
u,\, v\in D^{1, 2}(\RN),\quad u,\, v>0\,\,\hbox{in $\RN\setminus\{0\}$},\end{cases}$$
where $N\ge 3$, $\la_1, \la_2\in (0, \frac{(N-2)^2}{4})$, $2^\ast=\frac{2N}{N-2}$ and $\al>1, \bb>1$ satisfying $\al+\bb=2^\ast$. This problem
is related to coupled nonlinear Schr\"{o}dinger equations with critical exponent
for Bose-Einstein
condensate.
For different ranges of $N$, $\al$, $\bb$ and $\nu>0$,
we obtain positive ground state solutions via some quite different variational methods,
which are all radially symmetric. It turns out that the least energy level depends heavily on the relations among $\al, \,\bb$ and $2$.
Besides, for sufficiently small $\nu>0$, positive solutions
are also obtained via a variational perturbation approach.
Note that the Palais-Smale condition can not hold for any positive energy level,
which makes the study via variational methods rather complicated.\\
\end{minipage}
\end{center}

\vskip0.16in \s{Introduction}

In this paper we consider solitary wave solutions of coupled nonlinear
Schr\"{o}dinger equations, known in the literature as Gross-Pitaevskii equations (\cite{HMEWC, TV}):
\be\label{eq0001}
\begin{cases}
-i\frac{\partial}{\partial t}\Phi_1=\Delta \Phi_1- a(x)\Phi_1+
\mu_1 |\Phi_1|^2 \Phi_1+\nu |\Phi_2|^2\Phi_1,\,\, x\in \RN, \,\,t>0,\\
-i\frac{\partial}{\partial t}\Phi_2=\Delta \Phi_2- b(x)\Phi_2+
\mu_2|\Phi_2|^2 \Phi_2+\nu |\Phi_1|^2\Phi_2,\,\, x\in \RN, \,\,t>0,\\
\Phi_j=\Phi_j(x,t)\in\mathbb{C},\quad j=1,2,\\
\Phi_j(x,t)\to 0,  \quad\hbox{as}\,\, |x|\to+\iy, \,\,t>0, \,\,j=1,2,
\end{cases}\ee
where $i$ is the imaginary unit, $a(x), b(x)$ are potential functions,
$\mu_1,\mu_2 >0$ and $\nu\neq 0$ is a coupling constant.
System (\ref{eq0001}) appears in many physical problems, especially in nonlinear optics.
Physically, the solution $\Phi_j$ denotes the $j^{th}$ component of the beam
in Kerr-like photorefractive media (see \cite{AA}).
The positive constant $\mu_j$ is for self-focusing in the $j^{th}$ component of
the beam. The coupling constant $\nu$ is
the interaction between the two components of the beam. Problem (\ref{eq0001}) also arises
in  the Hartree-Fock theory for a double condensate, i.e., a binary mixture of Bose-Einstein
condensates in two different hyperfine states $|1\rangle$ and $|2\rangle$ (see \cite{EGBB}).
Physically, $\Phi_j$ are the corresponding condensate amplitudes, $\mu_j$ and $\nu$ are the intraspecies and
interspecies scattering lengths. The sign of $\nu$ determines whether the interactions of
states $|1\rangle$ and $|2\rangle$
are repulsive or attractive, i.e., the interaction is attractive if $\nu>0$, and the interaction is repulsive
if $\nu<0$, where the two states are in strong competition.

To obtain solitary wave solutions of system (\ref{eq0001}), we set $\Phi_1(x, t)=e^{i \la_1 t}u(x)$
and $\Phi_2(x, t)=e^{i \la_2 t}v(x)$. Write $V_1(x)=a(x)+\la_1$ and $V_2(x)=b(x)+\la_2$ for convenience, and we are only interested in nonnegative solutions, then system (\ref{eq0001})
is reduced to the following elliptic system
\be\label{eq0002}
\begin{cases}-\Delta u +V_1(x) u =
\mu_1 u^3+\nu uv^2, \quad x\in \R^N,\\
-\Delta v +V_2(x) v =\mu_2 v^3+\nu vu^2,  \quad   x\in
\R^N,\\
u\ge 0, v\ge 0 \,\,\hbox{in $\RN$},\quad u(x), v(x)\to0\,\,\,\hbox{as}\,\,\,|x|\to\iy.\end{cases}\ee
This Bose-Einstein condensate type system (\ref{eq0002}) is a special case of the following problem
\be\label{eq0-0-2}
\begin{cases}-\Delta u +V_1(x) u =
 \mu_1 u^{2p-1}+\nu u^{p-1}v^{p}, \quad x\in \RN,\\
-\Delta v +V_2(x) v =\mu_2 v^{2p-1}+\nu v^{p-1} u^{p},  \quad   x\in
\RN,\\
u\ge 0, v\ge 0 \,\,\hbox{in $\RN$},\quad u(x), v(x)\to0\,\,\,\hbox{as}\,\,\,|x|\to\iy,\end{cases}\ee
where $p>1$ and $p\le 2^\ast/2$ if $N\ge 3$, $2^\ast=\frac{2N}{N-2}$ is the
critical Sobolev exponent. If $p=2$, then (\ref{eq0-0-2}) turns to be the cubic system (\ref{eq0002}).
For further introduction about this problem, readers can also see the survey articles \cite{F, KLu}, which also contain information
about the physical relevance of non-cubic nonlinearities (e.g. quintic).
For the subcritical case $p <2^\ast/2$, the existence and multiplicity of solutions to (\ref{eq0-0-2}) have been widely studied under different
assumptions on $V_i$ and $\nu$, see \cite{AC2, LW1, LW, MMP, Pomponio, S, TV} and references therein.

\vskip0.10in

Remark that, all the papers mentioned above deal with the subcritical case. To the
best of our knowledge, there is no existence results for (\ref{eq0-0-2}) in the critical case $2p=2^\ast$ in
the literature.

\vskip0.10in

In this paper, we study (\ref{eq0-0-2})
in critical case where $N\ge 3$ and $2p=2^\ast$. In this case, if $V_i(x)=\la_i$ are nonzero constants with the same sign,
then by Pohozaev identity,
we easily conclude that any solution $(u, v)$ of (\ref{eq0-0-2}) satisfy $\int_{\R^N}\la_1 u^2+\la_2 v^2=0$, and so $(u, v)=(0, 0)$.
Hence we do not consider the case $V_i(x)=\la_i$ here, and in the sequel we assume that $V_i(x)=-\frac{\la_i}{|x|^2}$ are Hardy type potentials.
The Hardy type potentials,
which arise in several physical contexts (e.g., in nonrelativistic quantum mechanics, molecular physics, quantum cosmology,
linearization of combustion models), do not belong to the Kato's class, so they cannot be regarded as a lower order perturbation term. In particular,
 any nontrivial solutions of (\ref{eq0-0-2}) with $V_i(x)=-\frac{\la_i}{|x|^2}$ are singular
at $x=0$. For the sake of simplicity, in the sequel we assume that $\mu_1=\mu_2=1$.
Then, to study (\ref{eq0-0-2}) with $2p=2^\ast$, $V_i(x)=-\frac{\la_i}{|x|^2}$ and $\mu_1=\mu_2=1$, we turn to study the following general problem
\be\label{eq0-0}\begin{cases}-\Delta u -\frac{\la_1}{|x|^2}u=u^{2^\ast-1}+
\nu \al u^{\al-1}v^\bb, \quad x\in \RN,\\
  -\Delta v -\frac{\la_2}{|x|^2}v=v^{2^\ast-1} +
\nu \bb u^{\al}v^{\bb-1},     \quad x\in \RN,\\
u,\, v\in D^{1, 2}(\RN),\quad u,\, v>0\,\,\hbox{in $\RN\setminus\{0\}$},\end{cases}\ee
where $N\ge 3$, $\la_1, \la_2\in (0, \La_N)$, $\La_N:=\frac{(N-2)^2}{4}$,
\be\label{eq1-8}\al>1,\quad \bb>1,\quad\al+\bb=2^\ast,\ee
and $D^{1, 2}(\RN)$ is
the completion of $C_0^\iy(\RN)$ with respect to the norm
$$\|u\|:=\left(\intR{|\nabla u|^2}\right)^{1/2}.$$
Note that if $\al=\bb=p=2^\ast/2$, then (\ref{eq0-0}) turns to be (\ref{eq0-0-2}) with $V_i(x)=-\frac{\la_i}{|x|^2}$ and $\mu_1=\mu_2=1$.
The mathematical interest in system (\ref{eq0-0}) relies on their double criticality, due to the fact that both the exponent
of the nonlinearities (which is critical in the sense of the Sobolev embedding)
and the singularities share the same order of homogeneity as the Laplacian.
{\it The main goal of this paper is to study the existence and radial symmetry of ground state solutions to system (\ref{eq0-0}) }, where the ground state solution is defined in Definition \ref{definition} below.

Recall that $\la_1, \la_2\in (0, \La_N)$, from Hardy inequality
\be\label{Hardy}\La_N\intR{\frac{u^2}{|x|^2}} \le\intR{|\nabla u|^2},\quad \forall\, u\in D^{1, 2}(\RN),\ee
we see that $\|\cdot\|_{\la_i}, i=1, 2,$ are equivalent norms
to $\|\cdot\|$, where
\be\label{eq1-4}\|u\|_{\la_i}^2:=\intR{|\nabla u|^2-\frac{\la_i}{|x|^2}u^2}.\ee
Denote the norm of $L^p(\RN)$ by $|u|_p =
(\int_{\RN}|u|^p\,dx)^{\frac{1}{p}}$. The case of a single equation has been deeply investigated in the literature. In particular,
by \cite{Terracini}, the problem
\be\label{eq1-2}\begin{cases}-\Delta u -\frac{\la_i}{|x|^2}u=u^{2^\ast-1}, \quad x\in \RN,\\
u(x)\in D^{1, 2}(\RN),\quad u>0\,\,\hbox{in $\RN\setminus\{0\}$},\end{cases}\ee
has exactly an one-dimensional $C^2$ manifold of positive solutions given by
\be\label{eq1-3}Z_{i}=\left\{z_\mu^i(x)=\mu^{-\frac{N-2}{2}}z_1^i\left(\frac{x}{\mu}\right),\quad\mu>0\right\},\ee
where
\be\label{eq1--3}z_1^i(x)=\frac{A(N, \la_i)}{|x|^{a_{\la_i}}\left(1+|x|^{2-\frac{4 a_{\la_i}}{N-2}}\right)^{\frac{N-2}{2}}},\ee
$a_{\la_i}=\frac{N-2}{2}-\sqrt{\frac{(N-2)^2}{4}-\la_i}$ and $A(N, \la_i)=\frac{N(N-2-2a_{\la_i})^2}{N-2}$.
Moreover, all positive solutions of (\ref{eq1-2}) satisfy
\be\label{eq1-5}S(\la_i):=\inf_{u\in D^{1, 2}(\RN)\setminus\{0\}}\frac{\|u\|_{\la_i}^2}{|u|_{2^\ast}^2}
=\frac{\|z_\mu^i\|_{\la_i}^2}{|z_\mu^i|_{2^\ast}^2}=\left(1-\frac{4\la_i}{(N-2)^2}\right)^{\frac{N-1}{N}}S,\ee
and
\be\label{eq1-6}I_{\la_i}(z_\mu^i)=\frac{1}{N}\|z_\mu^i\|_{\la_i}^2=\frac{1}{N}|z_\mu^i|_{2^\ast}^{2^\ast}=\frac{1}{N}S(\la_i)^{N/2},\ee
where
$S$ is the
sharp constant of $D^{1,2}(\RN)\hookrightarrow L^{2^\ast}(\RN)$
\be\label{sobolev}\intR{|\nabla u|^2}\ge
S\left(\intR{|u|^{2^\ast}}\right)^{\frac{2}{2^\ast}},\ee
and
\be\label{eq1-7}I_{\la_i}(u):=\frac{1}{2}\|u\|^2_{\la_i}-\frac{1}{2^\ast}\intR{|u|^{2^\ast}},\quad i=1, 2,\ee
see \cite{Terracini}.
There are also many papers working on related equations with a Hardy type potential and a critical nonlinearity,
we refer readers to \cite{AFP1, FP, Smets} and references therein.

We call a solution $(u, v)$ of (\ref{eq0-0})  {\it nontrivial} if both $u\not\equiv 0$ and $v\not\equiv 0$; we call a solution $(u, v)$ {\it positive} if
both $u>0$ and $v>0$ in $\RN\setminus\{0\}$; we call a solution $(u, v)$ {\it semi-trivial} if $(u, v)$ is a type of $(u, 0)$ or $(0, v)$.

One of the difficulties in the study of (\ref{eq0-0})
is that it has semi-trivial solutions $(z_\mu^1, 0)$ and $(0, z_\mu^2)$. Here, we are only
interested in nontrivial solutions of (\ref{eq0-0}). Define $\mathbb{D}:=D^{1, 2}(\RN)\times D^{1, 2}(\RN)$ with norm
$$\|(u, v)\|_{\BD}^2:=\|u\|_{\la_1}^2+\|v\|_{\la_2}^2.$$
Then nontrivial solutions of (\ref{eq0-0}) can be found as nontrivial critical points of the $C^1$ functional $J_\nu: \BD\to \R$, where
$$J_\nu(u, v):=I_{\la_1}(u)+I_{\la_2}(v)-\nu\intR{|u|^\al |v|^\bb}.$$

Another difficulty is the failure of the Palais-Smale condition, which makes the study of (\ref{eq0-0}) very tough. Since (\ref{eq0-0})
is invariant under the
transformation $(u(x), v(x))\mapsto (\mu^{\frac{N-2}{2}}u(\mu x), \mu^{\frac{N-2}{2}}v(\mu x))$, where $\mu>0$,
it is easy to see that the Palais-Smale condition ($(PS)$ condition for short) can not hold for any energy level $c>0$. In fact,
assume by contradiction that $(PS)$ condition holds for some $c>0$, and let
$(u_n, v_n)$ be a $(PS)_c$ sequence, that is,
$J_\nu (u_n, v_n)\to c$ and $ J_\nu'(u_n, v_n)\to 0$ as $n\to\iy.$
Then up to a subsequence, we may assume that $(u_n, v_n)\to (u, v)$ strongly in $\BD$. Define
$(\wu_n (x), \wv_n (x)):=(n^{\frac{N-2}{2}}u_n(n x), n^{\frac{N-2}{2}}v_n(n x))$,
then it is easy to check that $(\wu_n, \wv_n)$ is also a $(PS)_c$ sequence and
$(\wu_n, \wv_n)\rightharpoonup (0, 0)$ weakly in $\BD$. Since $(PS)_c$ condition holds,
we have $(\wu_n, \wv_n)\to (0, 0)$ strongly in $\BD$, which contradicts with $c>0$.

\bd\label{definition} We say a solution $(u_0, v_0)$ of (\ref{eq0-0}) is a ground state solution if $(u_0, v_0)$ is nontrivial
and $J_\nu (u_0, v_0)\le J_{\nu}(u, v)$ for any other nontrivial solution $(u, v)$ of (\ref{eq0-0}).\ed

To obtain ground state solutions of (\ref{eq0-0}), as in \cite{LW1}, we define
{\allowdisplaybreaks
\begin{align}\label{nehari}\N_\nu:=\Big\{ &(u, v)\in \BD \,:\, u\not\equiv 0, v\not\equiv 0, \|u\|_{\la_1}^2
=\int_{\RN}( |u|^{2^\ast}+\nu\al |u|^{\al}|v|^\bb),\nonumber\\
& \| v\|_{\la_2}^2=\int_{\RN}( |v|^{2^\ast}+\nu\beta |u|^\al|v|^\bb)\Big\}.\end{align}
}%
Then any nontrivial solutions of (\ref{eq0-0}) has to belong to $\N_\nu$. Take $\varphi, \psi\in C^{\iy}_0 (\RN)$ with $\varphi, \psi\not\equiv 0$ and $supp(\varphi)\cap supp(\psi)=\emptyset$,
then there exist $t_1,t_2>0$ such that $(t_1\varphi, t_2\psi)\in \mathcal{N}_\nu$ for any $\nu\neq 0$.
So $\N_\nu\neq \emptyset$. We set
\be\label{least-energy} c_\nu:=\inf_{(u, v)\in \N_\nu} J_\nu(u, v)=\inf_{(u, v)\in \N_\nu}\frac{1}{N}\left(\|u\|_{\la_1}^2+\|v\|_{\la_2}^2\right).\ee
By (\ref{eq1-5}) we have
\be\label{eq1-5-1}\|u\|_{\la_i}^2\ge S(\la_i) |u|_{2^\ast}^2\quad\forall\,u\in D^{1, 2}(\RN),\,\,i=1,2.\ee
Then it is easy to see that $c_\nu>0$ for all $\nu$.
Moreover, if $(u_0, v_0)$ is a nontrivial solution satisfying $J_\nu (u_0, v_0)=c_\nu$, then $(u_0, v_0)$ is a ground state solution.
Our first result is concerned with ground state solutions with energy below $\frac{1}{N}\min\{S(\la_1)^{N/2},\, S(\la_2)^{N/2}\}$.

\bt\label{th1} Assume that $N\ge 3$, $\la_1, \la_2\in (0, \La_N)$ and (\ref{eq1-8}) hold.

\begin{itemize}

\item [(1)] If $\nu<0$, then $c_\nu\equiv\frac{1}{N}S(\la_1)^{N/2}+\frac{1}{N}S(\la_2)^{N/2}$, and $c_\nu$ can not be attained.

\item [(2)] Let
\be\label{eq2-31}\nu_0:=\frac{1}{2^\ast}\left[\left(1+\max\left\{\frac{\La_N-\la_1}{\La_N-\la_2}, \frac{\La_N-\la_2}{\La_N-\la_1}\right\}\right)^{\frac{2^\ast}{2}}-2\right]>0,\ee
then for all $\nu> \nu_0$, (\ref{eq0-0}) has a positive ground state solution $(u_\nu, v_\nu)\in \BD$, which is radially symmetric and satisfies
\be\label{eq1-10}J_\nu(u_\nu, v_\nu)=c_\nu<\frac{1}{N}\min\left\{S(\la_1)^{N/2},\,\, S(\la_2)^{N/2}\right\}.\ee

\item [(3)] If one of the following conditions
\begin{itemize}

\item[$({\bf C_1})$] $N\ge 5$ and $\max\{\al, \bb\}<2$,

\item[$({\bf C_2})$] $\la_1\le\la_2$ and $\al<2$,

\item[$({\bf C_3})$] $\la_2\le\la_1$ and $\bb<2$,
\end{itemize}

holds, then for all $\nu>0$, (\ref{eq0-0}) has a positive ground state solution $(u_\nu, v_\nu)\in\BD$, which is radially symmetric and satisfies (\ref{eq1-10}).
\end{itemize}

\et

Now we want to obtain ground state solutions with energy above the value $\frac{1}{N}\max\{S(\la_1)^{N/2},\, S(\la_2)^{N/2}\}$, which seems
much more interesting to us.
To this goal, by Theorem \ref{th1} (2)-(3) we have to assume that $\min\{\al, \bb\}\ge 2$ and $\nu>0$ is small. In this case,
since $4\le\al+\bb=2^\ast$, so $N=3$ or $N=4$. Moreover, if $N=4$, then we must have $\al=\bb=2$.
Note that if $N=4$ and $\al=\bb=2$, then (\ref{eq0-0}) turns to be the following cubic system
\be\label{eq4-1}\begin{cases}-\Delta u -\frac{\la_1}{|x|^2}u=u^{3}+
2\nu  uv^2, \quad x\in \R^4,\\
  -\Delta v -\frac{\la_2}{|x|^2}v=v^{3} +
2\nu u^{2}v,     \quad x\in \R^4,\\
u,\, v\in D^{1, 2}(\R^4),\quad u,\, v>0\,\,\hbox{in $\R^4\setminus\{0\}$},\end{cases}\ee
which is just the Bose-Einstein condensate type system (\ref{eq0002}) with $V_i(x)=-\frac{\la_i}{|x|^2}$ in critical case $N=4$.
Note that $\La_4=1$. Then we have the following results.

\bt\label{th2}Assume that $N=4$, $\al=\bb=2$ and $\la_1, \la_2\in (0, 1)$. Define
\be\label{th2constant}\nu_1:=\min\frac{1}{2}\left\{\frac{1-\la_1}{1-\la_2},\,\,\frac{1-\la_2}{1-\la_1},\,
\,\frac{(1-\la_1)^{\frac{3}{4}}(1-\la_2)^{\frac{3}{4}}}{(1-\la_1)^{\frac{3}{2}}+(1-\la_2)^{\frac{3}{2}}}\right\}.\ee
Then for any $\nu\in (0, \nu_1)$, (\ref{eq4-1}) has a positive ground state solution $(u_\nu, v_\nu)\in \BD$, which satisfies
\be\label{eq1-10-1}c_\nu=J_\nu(u_\nu, v_\nu)\to \frac{1}{4}\left(S(\la_1)^2+S(\la_2)^2\right),\quad\hbox{as $\nu\to 0$}.\ee
\et

\bt\label{theorem2}Assume that $N=3$, $\al+\bb=2^\ast$, $\al\ge 2$, $\bb\ge 2$ and $\la_1, \la_2\in (0, \La_3)$.
Then there exists $\widetilde{\nu}_1>0$ such that
 for any $\nu\in (0, \widetilde{\nu}_1)$, (\ref{eq0-0}) has a positive ground state solution $(u_\nu, v_\nu)\in \BD$, which satisfies
\be\label{eq1-10-2}J_\nu(u_\nu, v_\nu)\to \frac{1}{3}\left(S(\la_1)^{3/2}+S(\la_2)^{3/2}\right),\quad\hbox{as $\nu\to 0$}.\ee
\et

\br

\begin{itemize}

\item[(1)]  (\ref{eq1-10-1})-(\ref{eq1-10-2}) yield $c_\nu>\frac{1}{N}\max\{S(\la_1)^{N/2},\, S(\la_2)^{N/2}\}$ for $\nu>0$ small appropriately,
that is, we obtain positive ground state solutions with energy above $\frac{1}{N}\max\{S(\la_1)^{N/2},\, S(\la_2)^{N/2}\}$,
so the case $\min\{\al, \bb\}\ge 2$ is completely different
from the cases studied in Theorem \ref{th1}-(3).
As we can see in the following sections, the case $\min\{\al, \bb\}\ge 2$ is much more complicated (see Theorem \ref{th2-1} for example). Besides,
if $\la_1=\la_2$, then Theorem \ref{th2} will be improved by Theorem \ref{th5} in Section 7.

\item[(2)] The case $N= 3$ is much more tough than the case $N=4$, and we can not give an accurate definition of $\widetilde{\nu}_1$ in Theorem \ref{theorem2} as (\ref{th2constant}) unfortunately. As we will see in the following sections, the idea of proving Theorem \ref{th2} takes full use of the fact $\al=\bb=2$, and can not be used in the case $N=3$. Meanwhile, the idea of proving Theorem \ref{theorem2} is quite different and more general.

\end{itemize}

\er

As we will see in Section 2, the radial symmetry of ground state solutions obtained in Theorem \ref{th1} is an easy
corollary of the Schwartz symmetrization. However,
the Schwartz symmetrization can not be used to
prove the radial symmetry of ground state solutions obtained in Theorems \ref{th2} and \ref{theorem2}.
Here, to get the radial symmetry of solutions obtained in Theorems \ref{th2} and \ref{theorem2}, we will use the moving planes method. Precisely, we have
the following result.

\bt\label{th4} Assume that $N=3$ or $N=4$, $\al+\bb=2^\ast$, $\al\ge 2$, $\bb\ge 2$, $\la_1, \la_2\in (0, \La_N)$ and $\nu>0$. Then any positive solutions of (\ref{eq0-0}) is radially symmetric with respect to the origin.
Therefore, ground state solutions $(u_\nu, v_\nu)$ obtained in Theorems \ref{th2} and Theorems \ref{theorem2} are radially symmetric.\et

There are some other special cases,
such as the case in which $N=3, 4, 5$, $1<\al< 2\le\bb$, $\al+\bb=2^\ast$, $\la_2<\la_1$ and $\nu>0$ sufficiently small,
where we have no idea whether the ground state solutions exist or not. This remains to be an interesting open question.
Here we can obtain positive solutions for these cases if $\nu>0$ is sufficiently small. Precisely, we have the
following result, which plays a crucial role in the proof of Theorem \ref{theorem2}.

\bt\label{th1-4}Assume that $N\ge 3$, $\la_1, \la_2\in (0, \La_N)$ and (\ref{eq1-8}) hold.
Then there exists $\nu_2>0$ such that for any $\nu\in (0, \nu_2]$, (\ref{eq0-0}) has a positive solution $(u_\nu, v_\nu)\in \BD$, which is radially symmetric with respect to the origin and satisfies
\be\label{eq1-11} J_\nu (u_\nu, v_\nu)< \frac{1}{N}\left(S(\la_1)^{\frac{N}{2}}+S(\la_2)^{\frac{N}{2}}\right).\ee
\et

We should mention that, Abdellaoui, Felli and Peral \cite{AFP} studied the following class of weakly coupled nonlinear elliptic equations
\be\label{eq1-1}\begin{cases}-\Delta u -\frac{\la_1}{|x|^2}u=u^{2^\ast-1}+
\nu h(x) \al u^{\al-1}v^\bb , \quad x\in \RN,\\
  -\Delta v -\frac{\la_2}{|x|^2}v=v^{2^\ast-1} +
\nu h(x) \bb u^{\al}v^{\bb-1},     \quad x\in \RN,\\
u,\, v\in D^{1, 2}(\RN),\quad u,\, v>0\,\,\hbox{in $\RN\setminus\{0\}$}.\end{cases}\ee
Note that if $h(x)\equiv 1$, then (\ref{eq1-1}) turns to be (\ref{eq0-0}).
For the case (\ref{eq1-8}),
they assumed the following condition on $h(x)$
\begin{itemize}
\item [$({\bf H_1})$]  $h\in L^\iy(\RN),\,h\ge 0,\, h\not\equiv 0$, $h$ is continuous in a neighborhood of $0$ and $\iy$,
and $h(0)=\lim_{|x|\to \iy}h(x)=0.$
\end{itemize}
and then they proved some existence results of ground state solutions
for (\ref{eq1-1}) in case $\nu>0$ (see Section 4 in \cite{AFP}). We should point out that,
under condition $({\bf H_1})$, the Palais-Smale condition holds for energy level $c$ with
\be\label{eq1-9} c<\frac{1}{N}\min\left\{S(\la_1)^{N/2},\,\, S(\la_2)^{N/2}\right\},\ee
(see \cite[Lemmas 4.1 and 4.3]{AFP}), which plays a crucial role in
obtaining ground state solutions in \cite{AFP}. Therefore, {\it problem (\ref{eq0-0}) is completely different from (\ref{eq1-1})}. Moreover, there are no results about the existence of ground state solutions to (\ref{eq1-1}) with energy above $\frac{1}{N}\max\{S(\la_1)^{N/2},\, S(\la_2)^{N/2}\}$ in \cite{AFP}.
It was only pointed out in \cite[Remark 4.6]{AFP} that for the case where $h(x)\equiv 1$, $\nu>0$, $\al+\bb=2^\ast$ and $\la_1=\la_2=\la$,  it is easy to
construct, by a direct computation, positive solutions to (\ref{eq0-0}) of the form $(\phi, c\phi), c>0$.
Remark that whether these solutions $(\phi, c\phi)$ are ground state solutions  is not known in \cite{AFP},
and there are no any conclusions about (\ref{eq0-0}) for the general case $\la_1\neq \la_2$ in \cite{AFP}.

The rest of this paper proves these theorems, and we give some notations here. In the sequel, we denote positive
constants (possibly different in different places) by $C, C_1, C_2,\cdots$, and $B(x,r):=\{y\in \RN : |x-y|<r\}$.
Denote $B_r:= B(0, r)$ for convenience.
The paper is organized as follows.

We give the proof of Theorem \ref{th1} in Section 2, where we will use the concentration-compactness
principle from Lions (\cite{Lions, Lions1}) and some ideas from \cite{AFP}. The proof of Theorem \ref{th2}
is given in Section 3, where we will borrow some ideas from  \cite{Terracini} and the authors' paper \cite{CZ1}.

In Section 4, we will prove Theorem \ref{th1-4} via a perturbation method, where we will use some ideas from Byeon and Jeanjean \cite{BJ}. In order to
construct a spike solution of the following nonlinear elliptic problem
$$-\e^2\Delta u+V(x)u=f(u),\ u\in H^1(\RN),$$
for a general subcritical nonlinearity $f(u)$ and sufficiently small $\e>0$ ,
Byeon and Jeanjean \cite{BJ} developed a new variational approach.
We will mainly follow this variational approach to prove Theorem \ref{th1-4}.
Note that this approach can not be used directly, and we need some
crucial modifications for our proof. For example, we will define a special mountain-pass
value $a_{\nu}$, where all paths are required to be bounded in
$\BD$ by a same constant which is independent of
$\nu$. This special $a_\nu$ is essential to our proof.
Moreover, we should point out
that the lack of compactness is the main difficulty because of the failure
of the $(PS)$ condition of (\ref{eq0-0}),
and especially because $Z_i, i=1,2$ are not compact in $\BD$.

In Section 5, we will prove Theorem \ref{theorem2} with the help of Theorem \ref{th1-4}. Here, quite different ideas are needed comparing to those of proving
Theorem \ref{th2} in Section 3.

In Section 6, we will prove Theorem \ref{th4} via the moving planes method. The moving planes method has been used by many authors to prove
symmetry and monotonicity of positive solutions to various nonlinear elliptic problems, we refer readers to \cite{CL, CLO, GNN1, GNN2} and references therein.

Finally, by following some arguments from the authors' papers \cite{CZ1, CZ2}, we will
give some remarks for the special case $\la_1=\la_2$, $\al=\bb=2^\ast/2$ and $\nu>0$ in Section 7, where
some uniqueness results about the ground state solutions will be obtained, see Theorems \ref{th5} and \ref{th6}.
In authors' papers \cite{CZ1, CZ2}, we studied the following Bose-Einstein condensation system for critical case
\be\label{CZ}
\begin{cases}-\Delta u +\la_1 u =
\mu_1 u^{2^\ast-1}+\nu u^{\frac{2^\ast}{2}-1}v^{\frac{2^\ast}{2}}, \quad x\in \Omega,\\
-\Delta v +\la_2 v =\mu_2 v^{2^\ast-1}+\nu v^{\frac{2^\ast}{2}-1} u^{\frac{2^\ast}{2}},  \quad   x\in
\Om,\\
u\ge 0, v\ge 0 \,\,\hbox{in $\Om$},\quad
u=v=0  \,\,\hbox{on $\partial\Om$}.\end{cases}\ee
Here, $\Om\subset \R^N (N\ge 4)$ is a smooth bounded domain and $\mu_1, \mu_2>0$ (the special case $N=4$ was studied in \cite{CZ1},
and the general case $N\ge 5$ in \cite{CZ2}). When $\nu=0$, (\ref{CZ}) turns to be
the well-known Brezis-Nirenberg problem (\cite{BN}).
Thanks to the celebrated idea from Brezis and Nirenberg \cite{BN},
we can show that the $(PS)$ condition of (\ref{CZ}) holds for some ranges of energy level (see \cite{CZ1, CZ2} for details).
Therefore, problem (\ref{eq0-0}) is also completely different from (\ref{CZ}).
Fortunately, some ideas of studying (\ref{CZ}) in \cite{CZ1, CZ2} can be used in this paper.

\vskip0.1in

\s{Proof of Theorem \ref{th1}}
\renewcommand{\theequation}{2.\arabic{equation}}

\subsection{The case $\nu<0$}

\bl\label{lemma1} If $c_\nu$ is attained by a couple $(u, v)\in\N_\nu$, then
this couple is a critical point of $J_\nu$, provided $\nu<0$.\el

\noindent {\bf Proof. } This proof is standard. Let $\nu<0$. Assume that $(u, v)\in\mathcal{N}_\nu$ such that $c_\nu=J_\nu(u, v)$. Define
{\allowdisplaybreaks
\begin{gather*}
    G_1(u, v)=\| u\|_{\la_1}^2
-\int_{\RN}( |u|^{2^\ast}+\nu\al |u|^{\al}|v|^\bb),\\
G_2(u, v)=\|v\|_{\la_2}^2-\int_{\RN}( |v|^{2^\ast}+\nu\bb |u|^{\al}|v|^\bb).
\end{gather*}
}%
Then there exist two Lagrange multipliers $K_1, K_2\in\R$ such that
\be\label{eq2-6}J_\nu'(u, v)+K_1 G_1'(u, v)+K_2G_2'(u, v)=0.\ee
Testing (\ref{eq2-6}) with $(u, 0)$ and $(0, v)$ respectively, we conclude from $(u, v)\in\N_\nu$ that
{\allowdisplaybreaks
\begin{align*}
&\left((2^\ast-2)|u|_{2^\ast}^{2^\ast}+\al(2-\al)|\nu|\int_{\RN}|u|^\al |v|^\bb\right)K_1-\al\bb |\nu|\left(\int_{\RN}|u|^\al |v|^\bb\right) K_2=0,\\
&\left((2^\ast-2)|v|_{2^\ast}^{2^\ast}+\bb(2-\bb)|\nu|\int_{\RN}|u|^\al |v|^\bb\right)K_2-\al\bb |\nu|\left(\int_{\RN}|u|^\al |v|^\bb\right) K_1=0.
\end{align*}
}%
Recall that
$|u|_{2^\ast}^{2^\ast}>\al|\nu|\int_{\RN}|u|^\al |v|^\bb$ and $ |v|_{2^\ast}^{2^\ast}>\bb|\nu|\int_{\RN}|u|^\al |v|^\bb$, we see from $\al+\bb=2^\ast$ that
{\allowdisplaybreaks
\begin{align*}
&\left((2^\ast-2)|u|_{2^\ast}^{2^\ast}+\al(2-\al)|\nu|\int_{\RN}|u|^\al |v|^\bb\right)\\
&\quad\times\left((2^\ast-2)|v|_{2^\ast}^{2^\ast}+\bb(2-\bb)|\nu|\int_{\RN}|u|^\al |v|^\bb\right)>\left(\al\bb |\nu|\int_{\RN} |u|^\al|v|^\bb\right)^2.
\end{align*}
}%
From this we deduce that $K_1=K_2=0$ and so $J_\nu'(u, v)=0$.\hfill$\square$

\bl\label{lemma2} Let $\nu<0$. For any $(u, v)\in \BD$ with $u\not\equiv 0$ and $v\not\equiv 0$, if
\be\label{eq2-7}\left(\int_{\RN}|u|^{2^\ast}\right)^\al \left(\int_{\RN}|v|^{2^\ast}\right)^\bb>\al^\al \bb^\bb \left(|\nu|\int_{\RN}|u|^{\al}|v|^\bb\right)^{2^\ast},\ee
then there exist $t_1>0, s_1>0$, such that $(t_1 u, s_1 v)\in\N_\nu$.\el

\noindent {\bf Proof. } For simplicity, we denote
$$A_1=\|u\|_{\la_1}^2,\,\, B_1=|u|_{2^\ast}^{2^\ast},\,\, C=|\nu|\int_{\RN}|u|^\al |v|^\bb,\,\, A_2=\|v\|_{\la_2}^2,\,\, B_2=|v|_{2^\ast}^{2^\ast}.$$
Recall the definition (\ref{nehari}) of $\N_\nu$, we see that $(t u, s v)\in\N_\nu$ for $t, s>0$ is equivalent to $t, s>0$ satisfy
\begin{align}\label{eq2-8}
A_1 t^{2-\al}= B_1 t^\bb-\al Cs^\bb,\quad A_2 s^{2-\bb}= B_2 s^\al-\bb C t^\al.
\end{align}
If $C=0$, then it is trivial to see that (\ref{eq2-8}) has a solution $(t_1, s_1)$ with $t_1, s_1>0$.
So we may assume that $C>0$. Then the equation $A_1 t^{2-\al}= B_1 t^\bb-\al Cs^\bb$ is equivalent to
$$s=g(t):=\left(\frac{B_1 t^\bb- A_1 t^{2-\al}}{\al C}\right)^{1/\bb}>0, \quad t>t_0:=\left(\frac{A_1}{B_1}\right)^{\frac{1}{2^\ast-2}}.$$
Therefore, it suffices to prove that
\be\label{eq2-9}A_2 \left(\frac{B_1 t^\bb- A_1 t^{2-\al}}{\al C}\right)^{\frac{2-\bb}{\bb}}
-B_2 \left(\frac{B_1 t^\bb- A_1 t^{2-\al}}{\al C}\right)^{\frac{\al}{\bb}}+\bb C t^\al=0\ee
has a solution $t>t_0$. Note that (\ref{eq2-9}) is equivalent to
$$f(t):=A_2 \left(\frac{B_1 - A_1 t^{2-2^\ast}}{\al C}\right)^{\frac{2-\bb}{\bb}}+t^{2^\ast-2}\left[\bb C-
B_2\left(\frac{B_1 - A_1 t^{2-2^\ast}}{\al C}\right)^{\frac{\al}{\bb}}\right]=0.$$
Note that (\ref{eq2-7}) implies that
$$\bb C-B_2\left(\frac{B_1}{\al C}\right)^{\frac{\al}{\bb}}<0,$$
then it is easy to check that
$\lim_{t\searrow t_0}f(t)>0$ and $ \lim_{t\to+\iy} f(t)=-\iy.$
So there exists $t_1>t_0>0$ such that $f(t_1)=0$. Let $s_1=g(t_1)$, then $s_1>0$ and $(t_1 u, s_1 v)\in \N_\nu$. This completes the proof.
\hfill$\square$\\

\noindent {\bf Proof of Theorem \ref{th1}-(1). } Fix any $\nu<0$.
Recall (\ref{eq1-3})-(\ref{eq1-6}), it is easy to see that
 $z_\mu^2\rightharpoonup 0$ weakly in $D^{1,2}(\RN)$ and so $(z_\mu^2)^\bb\rightharpoonup 0$ weakly in $L^{2^\ast/\bb}(\RN)$ as $\mu\to +\iy$.
That is,
$$\lim_{\mu\to+\iy}|\nu|\intR{(z_1^1)^{\al} (z_\mu^2)^\bb}=0.$$
Then (\ref{eq2-7}) holds for $(z_1^1, z_\mu^2)$ when $\mu>0$ sufficiently large, and so there exist $t_\mu, s_\mu>0$
such that $(t_\mu z_1^1, s_\mu z_\mu^2)\in\N_\nu$. Denote $F_\mu:=|\nu|\intR{(z_1^1)^{\al} (z_\mu^2)^\bb}$,
Then
\be\label{eq3-1}t_\mu^2 S(\la_1)^{\frac{N}{2}}=t_\mu^{2^\ast}S(\la_1)^{\frac{N}{2}}-\al t_\mu^\al s_\mu^\bb F_\mu,
\,\, s_\mu^2 S(\la_2)^{\frac{N}{2}} =s_\mu^{2^\ast} S(\la_2)^{\frac{N}{2}}-\bb t_\mu^\al s_\mu^\bb F_\mu.\ee
Assume that, up to a subsequence, $t_\mu\to +\iy$ as $\mu\to \iy$, then by $\bb(t_\mu^{2^\ast}-t_\mu^2)S(\la_1)^{\frac{N}{2}}=\al (s_\mu^{2^\ast}-s_\mu^2)S(\la_2)^{\frac{N}{2}}$ we also have $s_\mu\to+\iy$.
Then
$$t_\mu^{2^\ast}-t_\mu^2\ge\frac{1}{2}t_\mu^{2^\ast},\quad s_\mu^{2^\ast}-s_\mu^2\ge\frac{1}{2}s_\mu^{2^\ast}\,\,\,\hbox{for $\mu$ large enough.} $$
Combining this with (\ref{eq3-1}) we see that
$$\left(\frac{t_\mu}{s_\mu}\right)^\bb\le 2\al S(\la_1)^{-\frac{N}{2}} F_\mu\to 0,
\quad\left(\frac{s_\mu}{t_\mu}\right)^\al\le 2\bb S(\la_2)^{-\frac{N}{2}} F_\mu\to 0,\,\,\hbox{as $\mu\to +\iy$},$$
a contradiction. Therefore, $t_\mu$ and $s_\mu$ are uniformly bounded. Then by (\ref{eq3-1}) and $F_\mu\to 0$ as $\mu\to\iy$, we get that
$\lim_{\mu\to+\iy}(t_\mu, s_\mu)=(1, 1).$
Note that $(t_{\mu}z_1^1, s_{\mu} z_\mu^2)\in \mathcal{N}_\nu$, we see from (\ref{eq1-6}) and (\ref{least-energy}) that
{\allowdisplaybreaks
\begin{align*}c_\nu &\le J_\nu(t_{\mu}z_1^1, s_\mu z_\mu^2)=\frac{1}{N}\left(t_\mu^2\|z_1^1\|_{\la_1}^2+s_\mu^2\|z_\mu^2\|_{\la_2}^2\right)\\
&=\frac{1}{N}\left(t_\mu^2 S(\la_1)^{N/2}+s_\mu^2 S(\la_2)^{N/2}\right).
\end{align*}
}%
Letting $\mu\to+\iy$, we get that $c_\nu\le\frac{1}{N}(S(\la_1)^{N/2}+S(\la_2)^{N/2})$.
On the other hand, for any $(u, v)\in\mathcal{N}_\nu$, we see from $\nu<0$ and (\ref{eq1-5-1}) that
$$ \|u\|_{\la_1}^2
\le \int_{\RN} |u|^{2^\ast}\,dx\le  S(\la_1)^{-2^\ast/2}\|u\|_{\la_1}^{2^\ast},$$
and so $\|u\|_{\la_1}^2\ge S(\la_1)^{N/2}$. Similarly, $\|v\|_{\la_2}^2\ge S(\la_2)^{N/2}$.
Combining these with (\ref{least-energy}), we get that $c_\nu\ge\frac{1}{N}(S(\la_1)^{N/2}+S(\la_2)^{N/2})$. Hence,
\be\label{eq3-2}c_\nu=\frac{1}{N}\left(S(\la_1)^{N/2}+S(\la_2)^{N/2}\right).\ee

Now, assume that $c_\nu$ is attained by some $(u, v)\in \mathcal{N}_\nu$, then $(|u|, |v|)\in \mathcal{N}_\nu$ and $J_\nu(|u|, |v|)=c_\nu$.
By Lemma \ref{lemma1}, we know that $(|u|, |v|)$ is a nontrivial solution of (\ref{eq0-0}). By the
maximum principle, we may assume that $u>0, v>0$ in $\RN\setminus\{0\}$ and so $\intR{u^\al v^\bb}>0$. Then
$$ \|u\|_{\la_1}^2
< \int_{\RN} |u|^{2^\ast}\,dx\le  S(\la_1)^{-2^\ast/2}\|u\|_{\la_1}^{2^\ast},$$
Therefore, it is easy to see that
$c_\nu=J_\nu(u, v)>\frac{1}{N}(S(\la_1)^{N/2}+S(\la_2)^{N/2}),$
which is a contradiction. This completes the proof.\hfill$\square$

\subsection{The case $\nu>0$}

In this subsection, we let $\nu>0$. Define
\be\label{eq3-3} c'_\nu:=\inf_{(u, v)\in\mathcal{N}'_\nu}J_\nu(u, v),\ee
where
\begin{align}\label{eq3-4}\mathcal{N}'_\nu:=\Big\{(u, v)\in D\setminus\{ (0,0)\} : J_\nu'(u, v)(u, v)=0\Big\}. \end{align}
Note that $\mathcal{N}_\nu\subset\mathcal{N}'_\nu$, so
$c'_\nu\le c_\nu$. By (\ref{eq1-5-1}) it is easy to prove that $c'_\nu>0$. Moreover, it is standard to prove that
{\allowdisplaybreaks
\begin{align}\label{eq3-5}c'_\nu &=\inf_{(u, v)\in\BD\setminus\{(0, 0)\}}\max_{t>0}J_\nu(t u, t v)\nonumber\\
&=\inf_{(u, v)\in\BD\setminus\{(0, 0)\}}\frac{1}{N}\left[\frac{\|u\|_{\la_1}^2+\|v\|_{\la_2}^2}{\left(\int_{\RN}|u|^{2^\ast}+2^\ast \nu |u|^\al |v|^\bb +|v|^{2^\ast}\right)^{\frac{2}{2^\ast}}}\right]^{\frac{N}{2}}.\end{align}
}%
Define
$E(u, v):= |\nabla u|^2+|\nabla v|^2-\frac{\la_1}{|x|^2}|u|^2-\frac{\la_2}{|x|^2}|v|^2$
and $F(u, v):=|u|^{2^\ast}+2^\ast \nu |u|^\al |v|^\bb +|v|^{2^\ast}$ for simplicity, then
\be\label{eq3-6}\intR{E(u, v)}\ge \left(N c_\nu'\right)^{\frac{2}{N}}\left(\int_{\RN}F(u, v)\,dx\right)^{\frac{2}{2^\ast}},\quad\forall\,(u, v)\in \BD.\ee

The following lemma is the counterpart of Brezis-Lieb Lemma (\cite{BL}) for $(u, v)$, and the idea of its proof comes from \cite{BL}
(see also \cite[Lemma 1.32]{W}).

\bl\label{lemma3}Let $\Om\subset\RN$ be an open set and $(u_n, v_n)$ be a bounded sequence in $L^{2^\ast}(\Om)\times L^{2^\ast}(\Om)$.
If $(u_n, v_n)\to (u, v)$ almost everywhere in $\Om$, then
$$\lim_{n\to \iy}\int_{\Om}\left(|u_n|^\al |v_n|^\bb-|u_n -u|^\al |v_n- v|^\bb\right)\,dx=\int_{\Om}|u|^\al |v|^\bb\,dx.$$\el

\noindent {\bf Proof. } Fatou Lemma yields
$$\int_{\Om}|u|^{2^\ast}\le\liminf_{n\to\iy}\int_{\Om}|u_n|^{2^\ast}<\iy,\,\,\,\int_{\Om}|v|^{2^\ast}\le\liminf_{n\to\iy}\int_{\Om}|v_n|^{2^\ast}<\iy.$$
Recall that $\al, \bb$ satisfy (\ref{eq1-8}). For any $a_1, a_2, b_1, b_2\in\R$ and $\e>0$, we deduce from the mean value theorem and Young inequality that
{\allowdisplaybreaks
\begin{align*}
&\Big||a_1+a_2|^\al |b_1+b_2|^\bb- |a_1|^\al |b_1|^\bb\Big|\\
\le &\Big||a_1+a_2|^\al-|a_1|^\al\Big||b_1+b_2|^\bb+|a_1|^\al\Big||b_1+b_2|^\bb-|b_1|^\bb\Big|\\
\le & C\left[(|a_1|+|a_2|)^{\al-1}(|b_1|+|b_2|)^{\bb}|a_2|+|a_1|^\al (|b_1| +|b_2|)^{\bb-1}|b_2|\right]\\
\le & C\e\left[(|a_1|+|a_2|)^{2^\ast}+(|b_1|+|b_2|)^{2^\ast}\right]+C \e^{1-2^\ast}\left(|a_2|^{2^\ast}+|b_2|^{2^\ast}\right)\\
\le & C\e\left(|a_1|^{2^\ast}+|a_2|^{2^\ast}+|b_1|^\ast+|b_2|^{2^\ast}\right)+C \e^{1-2^\ast}\left(|a_2|^{2^\ast}+|b_2|^{2^\ast}\right).
\end{align*}
}%
Denote $\om_n=u_n-u$ and $\sg_n=v_n-v$, then
{\allowdisplaybreaks
\begin{align*}
f_n^\e:=&\bigg[\Big||u_n|^\al |v_n|^\bb-|\om_n|^\al |\sg_n|^\bb-|u|^\al |v|^\bb\Big|\\
&-C\e\left(|\om_n|^{2^\ast}+|u|^{2^\ast}+|\sg_n|^\ast+|v|^{2^\ast}\right)\bigg]_+\\
\le& |u|^\al |v|^\bb+C \e^{1-2^\ast}\left(|u|^{2^\ast}+|v|^{2^\ast}\right),
\end{align*}
}%
and so the dominated convergence theorem yields $\int_{\Om} f_n^\e\,dx\to 0$ as $n\to \iy$. Note that
$$\Big||u_n|^\al |v_n|^\bb-|\om_n|^\al |\sg_n|^\bb-|u|^\al |v|^\bb\Big|\le f_n^\e+C\e\left(|\om_n|^{2^\ast}+|u|^{2^\ast}+|\sg_n|^\ast+|v|^{2^\ast}\right),$$
we obtain
$$\limsup_{n\to\iy}\int_{\Om} \Big||u_n|^\al |v_n|^\bb-|\om_n|^\al |\sg_n|^\bb-|u|^\al |v|^\bb\Big|\le C\e.$$
Since $C>0$ is independent of $\e>0$, the proof is complete.\hfill$\square$\\

The following lemma is the counterpart of Lions' concentration-compactness principle (\cite{Lions, Lions1}) for problem (\ref{eq0-0}).

\bl\label{lemma7} Let $(u_n, v_n)\in\BD$ be a sequence such that
{\allowdisplaybreaks
\be\label{s3concentrate}\begin{cases}
(u_n, v_n)\rightharpoonup (u, v)\quad \hbox{weakly in $\BD$},\\
(u_n, v_n)\to (u, v)\quad \hbox{almost everywhere on $\RN$},\\
E(u_n-u, v_n-v)\rightharpoonup\mu \quad\hbox{in the sense of measures},\\
F(u_n- u, v_n- v)\rightharpoonup\rho \quad\hbox{in the sense of measures}.
\end{cases}\ee
}%
Define
{\allowdisplaybreaks
\begin{align}\label{s3concentrate1}&\mu_\iy:=\lim_{R\to\iy}\limsup_{n\to\iy}\int_{|x|\ge R}E(u_n, v_n)\,dx,\nonumber\\
&\rho_\iy:=\lim_{R\to\iy}\limsup_{n\to\iy}\int_{|x|\ge R}F(u_n, v_n)\,dx.
\end{align}
}%
Then it follows that
{\allowdisplaybreaks
\begin{align}
\label{s3eq2-16}&\|\mu\|\ge \left(N c_\nu'\right)^{\frac{2}{N}}\|\rho\|^{\frac{2}{2^\ast}},\\
\label{s3eq2-17}&\mu_\iy\ge \left(N c_\nu'\right)^{\frac{2}{N}}\rho_\iy^{\frac{2}{2^\ast}},\\
\label{s3eq2-18}&\limsup_{n\to\iy}\intR{E(u_n, v_n)}= \intR{E(u, v)}+\|\mu\|+\mu_\iy,\\
\label{s3eq2-19}&\limsup_{n\to\iy}\intR{F(u_n, v_n)}= \intR{F(u, v)}+\|\rho\|+\rho_\iy.
\end{align}
}%
Moreover, if $(u, v)=(0, 0)$ and $\|\mu\|= \left(N c_\nu'\right)^{\frac{2}{N}}\|\rho\|^{\frac{2}{2^\ast}}$, then $\mu$ and $\rho$ are concentrated at a single point.
\el

\noindent {\bf Proof. } In this proof we mainly follow the argument of \cite[Lemma 1.40]{W}. First we assume $(u, v)=(0, 0)$. For
any $h\in C_0^\iy (\RN)$, we see from (\ref{eq3-6}) that
\be\label{eq2-30}\intR{E( h u_n, h v_n)}\ge \left(N c_\nu'\right)^{\frac{2}{N}}\left(\int_{\RN}|h|^{2^\ast} F(u_n, v_n)\,dx\right)^{\frac{2}{2^\ast}}.\ee
Since $u_n\to 0, v_n \to 0$ in $L_{loc}^2(\RN)$, we have
$$\intR{E( h u_n, h v_n)}-\intR{|h|^2 E(  u_n,  v_n)} \to 0\quad\hbox{ as} \,\,\,n\to\iy.$$
Then by letting $n\to\iy$ in (\ref{eq2-30}), we obtain
\be\label{s3eq2-20}\int_{\RN}|h|^2\,d\mu\ge\left(N c_\nu'\right)^{\frac{2}{N}}\left(\int_{\RN}|h|^{2^\ast}\,d\rho\right)^{\frac{2}{2^\ast}},\ee
that is, (\ref{s3eq2-16}) holds.

For $R>1$, let $\psi_R\in C^1(\RN)$ be such that $0\le \psi_R\le 1$, $\psi_R (x)=1$ for $|x|\ge R+1$ and $\psi_R (x)=0$ for $|x|\le R$.
Then we see from (\ref{eq3-6}) that
$$\intR{E(\psi_R u_n, \psi_R v_n)}\ge \left(N c_\nu'\right)^{\frac{2}{N}}\left(\int_{\RN}|\psi_R|^{2^\ast} F(u_n, v_n)\,dx\right)^{\frac{2}{2^\ast}}.$$
Since $u_n\to 0, v_n \to 0$ in $L_{loc}^2(\RN)$, so
{\allowdisplaybreaks
\begin{align}\label{eq02-21}\limsup_{n\to\iy}&\intR{|\psi_R|^2 E(u_n, v_n)}\nonumber\\
&\ge
\left(N c_\nu'\right)^{\frac{2}{N}}\limsup_{n\to\iy}\left(\int_{\RN}|\psi_R|^{2^\ast} F(u_n, v_n)\,dx\right)^{\frac{2}{2^\ast}}.\end{align}
}%
Note that
$\int_{|x|\ge R+1} F(u_n, v_n) \le\int_{\RN}|\psi_R|^{2^\ast} F(u_n, v_n)\le\int_{|x|\ge R} F(u_n, v_n),$
so
\begin{align}\label{eq-lions1}\rho_\iy=\lim_{R\to\iy}\limsup_{n\to\iy}\int_{\RN} |\psi_R|^{2^\ast}F(u_n, v_n)\,dx.
\end{align}
On the other hand,
{\allowdisplaybreaks
\begin{align*}
&\limsup_{n\to\iy}\int_{\RN}|\psi_R|^2 E(u_n, v_n)\,dx\\
=&\limsup_{n\to\iy}\int_{|x|\ge R+1} E(u_n, v_n)\,dx+\limsup_{n\to\iy}\int_{R\le|x|\le R+1}|\psi_R|^2 E(u_n, v_n)\,dx\\
=&\limsup_{n\to\iy}\int_{|x|\ge R+1} E(u_n, v_n)\,dx+\limsup_{n\to\iy}\int_{R\le|x|\le R+1}|\psi_R|^2 (|\nabla u_n|^2+|\nabla v_n|^2)\,dx\\
\ge &\limsup_{n\to\iy}\int_{|x|\ge R+1} E(u_n, v_n)\,dx.
\end{align*}
}%
Letting $R\to\iy$ we see that $\mu_\iy\le\lim_{R\to\iy}\limsup_{n\to\iy}\int_{\RN}|\psi_R|^2 E(u_n, v_n)\,dx.$ Similarly,
{\allowdisplaybreaks
\begin{align*}
&\limsup_{n\to\iy}\int_{\RN}|\psi_R|^2 E(u_n, v_n)\,dx\\
=&\limsup_{n\to\iy}\int_{|x|\ge R} E(u_n, v_n)\,dx-\liminf_{n\to\iy}\int_{R\le|x|\le R+1}(1-|\psi_R|^2) E(u_n, v_n)\,dx\\
\le &\limsup_{n\to\iy}\int_{|x|\ge R} E(u_n, v_n)\,dx.
\end{align*}
}%
Letting $R\to\iy$ we see that $\mu_\iy\ge\lim_{R\to\iy}\limsup_{n\to\iy}\int_{\RN}|\psi_R|^2 E(u_n, v_n)\,dx.$ Hence
\be\label{eq-lions2}\mu_\iy=\lim_{R\to\iy}\limsup_{n\to\iy}\int_{\RN}|\psi_R|^2 E(u_n, v_n)\,dx.\ee
Then (\ref{s3eq2-17}) follows directly from (\ref{eq02-21}), (\ref{eq-lions1}) and (\ref{eq-lions2}).

Assume moreover that $\|\mu\|= \left(N c_\nu'\right)^{\frac{2}{N}}\|\rho\|^{\frac{2}{2^\ast}}$. Then by H\"{o}lder inequality and (\ref{s3eq2-20}), we have
$$\int_{\RN}|h|^{2^\ast}\,d\rho\le (N c_\nu')^{-\frac{2}{N-2}}\|\mu\|^{\frac{2}{N-2}}\int_{\RN}|h|^{2^\ast}\,d\mu,\quad\forall\,h\in C_0^\iy(\RN).$$
From this we deduce that $\rho=(N c_\nu')^{-\frac{2}{N-2}}\|\mu\|^{\frac{2}{N-2}}\mu$. So $\mu=(N c_\nu')^{\frac{2}{N}}\|\rho\|^{-\frac{2}{N}}\rho$,
and we see from (\ref{s3eq2-20}) that
$$\|\rho\|^{\frac{2}{N}}\left(\int_{\RN}|h|^{2^\ast}\,d\rho\right)^{\frac{2}{2^\ast}}\le\int_{\RN}|h|^2\,d\rho,\quad\forall\,h\in C_0^\iy(\RN).$$
That is, for each open set $\Om$, we have
$\rho(\Om)^{\frac{2}{2^\ast}}\rho(\RN)^{\frac{2}{N}}\le\rho(\Om).$ Therefore, $\rho$ is concentrated at a single point.

For the general case, we denote $\om_n=u_n-u$ and $\sg_n=v_n-v$, then $(\om_n, \sg_n)\rightharpoonup (0, 0)$ weakly in $\BD$.
From Brezis-Lieb Lemma (\cite{BL}) and Lemma \ref{lemma3}, we obtain for nonnegative $h\in C_0(\RN)$ that
{\allowdisplaybreaks
\begin{align*}
&\int_{\RN}h E(u, v)\,dx=\lim_{n\to\iy}\left(\int_{\RN}h E(u_n, v_n)\,dx-\int_{\RN}h E(\om_n, \sg_n)\,dx\right),\\
&\int_{\RN}h F(u, v)\,dx=\lim_{n\to\iy}\left(\int_{\RN}h F(u_n, v_n)\,dx-\int_{\RN}h F(\om_n, \sg_n)\,dx\right),
\end{align*}
}%
so
\begin{align}\label{eq-lions3}
&E(u_n, v_n)
\rightharpoonup E(u, v)+\mu, \,\,
F(u_n, v_n)\rightharpoonup F(u, v)+\rho, \,\,\hbox{in the sense of measures}.
\end{align}
Inequality (\ref{s3eq2-16}) follows from the corresponding one for $(\om_n, \sg_n)$. From Brezis-Lieb Lemma (\cite{BL}) and Lemma \ref{lemma3} again, it is easy
to prove that
{\allowdisplaybreaks
\begin{align*}&\mu_\iy:=\lim_{R\to\iy}\limsup_{n\to\iy}\int_{|x|\ge R}E(\om_n, \sg_n)\,dx,\\
&\rho_\iy:=\lim_{R\to\iy}\limsup_{n\to\iy}\int_{|x|\ge R}F(\om_n, \sg_n)\,dx,
\end{align*}
}%
and so inequality (\ref{s3eq2-17}) follows from the corresponding one for $(\om_n, \sg_n)$. For any $R>1$, we deduce from (\ref{eq-lions3}) that
{\allowdisplaybreaks
\begin{align*}
&\limsup_{n\to\iy}\int_{\RN}F(u_n, v_n)\\
=&\limsup_{n\to\iy}\left(\int_{\RN}|\psi_R|^{2^\ast} F(u_n, v_n)+\int_{\RN}(1-|\psi_R|^{2^\ast}) F(u_n, v_n)\right)\\
=&\limsup_{n\to\iy}\int_{\RN}|\psi_R|^{2^\ast} F(u_n, v_n)+\int_{\RN}(1-|\psi_R|^{2^\ast}) F(u, v)+\int_{\RN}(1-|\psi_R|^{2^\ast})\,d\rho.
\end{align*}
}%
Letting $R\to\iy$, we see from (\ref{eq-lions1}) that (\ref{s3eq2-19}) holds.
The proof of (\ref{s3eq2-18}) is similar. This completes the proof.\hfill$\square$

\bl\label{lemma8} Let $\nu>0$. Then (\ref{eq0-0}) has a solution $(u, v)\in\BD\setminus\{(0, 0)\}$ (maybe {\it  semi-trivial}), such that $J_\nu(u, v)=c_\nu'$ and
$u, v\ge 0$ are radially symmetric with respect to the origin. Moreover, if $c'_\nu<\frac{1}{N}\min\{S(\la_1)^{N/2}, S(\la_2)^{N/2}\}$, then
$(u, v)\in\BD$ is a positive ground state solution of (\ref{eq0-0}), and $c_\nu=c_\nu'=J_\nu(u, v)$.\el

\noindent{\bf Proof. }  For $(u, v)\in \mathcal{N}'_\nu$ with $u\ge 0, v\ge 0$,  we denote by $(u^*, v^*)$ its Schwartz
symmetrization. Then by the properties of Schwartz
symmetrization (see \cite{LL} for example), we see from $\la_1, \la_2, \nu>0$ that
$$\int_{\RN}
(|\nabla u^*|^2+|\nabla v^*|^2-\frac{\la_1}{|x|^2}|u^*|^2-\frac{\la_2}{|x|^2}|v^*|^2)
\le \int_{\RN}( |u^*|^{2^\ast}+2^\ast\nu |u^*|^{\al}|v^*|^{\bb}+ |v^*|^{2^\ast}).$$
Therefore, there exists $0<t^*\le 1$ such that $(t^*u^*, t^* v^*)\in\mathcal{N}'_\nu$, and then
{\allowdisplaybreaks
\begin{align}\label{s3eq2-22}
J_\nu(t^*u^*, t^* v^*)&=\frac{1}{N}(t^*)^2
(\| u^*\|_{\la_1}^2+\| v^*\|_{\la_2}^2)\nonumber\\
&\le\frac{1}{N}(\| u\|_{\la_1}^2+\| v\|_{\la_2}^2)=J_\nu(u, v).
\end{align}
}%
Therefore, we may take a minimizing sequence $(\wu_n, \wv_n)\in\mathcal{N}'_\nu$ of $c'_\nu$ such
that $(\wu_n, \wv_n)=(\wu_n^*, \wv_n^*)$ and $J_\nu(\wu_n, \wv_n) \to c'_\nu$ as $n\to\iy$. Define the L\'{e}vy concentration functions
$$Q_n(R):=\sup_{y\in\RN}\int_{B(y, R)}F(\wu_n, \wv_n)\,dx.$$
Since $\wu_n, \wv_n\ge 0$ are radially symmetric non-increasing,
one has that $Q_n(R)=\int_{B(0, R)}F(\wu_n, \wv_n)\,dx$. Then there exists $R_n>0$ such that
$$Q_n(R_n)=\int_{B(0, R_n)}F(\wu_n, \wv_n)\,dx=\frac{1}{2}\intR{F(\wu_n, \wv_n)}.$$
Define
$$(u_n(x), v_n(x)):=\left(R_n^{\frac{N-2}{2}}\wu_n(R_n x), \, R_n^{\frac{N-2}{2}}\wv_n(R_n x)\right),$$
Then by a direct computation, we see that $(u_n, v_n)\in\N_\nu'$, $J_\nu(u_n, v_n) \to c'_\nu$, $u_n, v_n\ge 0$ are radially symmetric non-increasing, and
\be\label{s3eq2-23}\int_{B(0, 1)}F(u_n, v_n)\,dx=\frac{1}{2}\intR{F(u_n, v_n)}=\sup_{y\in\RN}\int_{B(y, 1)}F(u_n, v_n)\,dx.\ee
From (\ref{s3eq2-22}) we know that $(u_n, v_n)$ are uniformly bounded in $\BD$. Then passing to a subsequence, there exist $(u, v)\in\BD$ and finite
measures $\mu, \rho$ such that (\ref{s3concentrate}) holds. Define $\mu_\iy, \rho_\iy$ as in (\ref{s3concentrate1}),
then by Lemma \ref{lemma7} we see that (\ref{s3eq2-16})-(\ref{s3eq2-19}) hold. Note that
\be\label{s3eq2-24}\|u_n\|_{\la_1}^2+\|v_n\|_{\la_2}^2=\intR{F(u_n, v_n)}\to N c_\nu',\quad\hbox{as $n\to\iy$},\ee
we conclude from (\ref{s3eq2-16})-(\ref{s3eq2-19}) and (\ref{eq3-6}) that
{\allowdisplaybreaks
\begin{align*}
&N c_\nu'=\intR{F(u, v)}+\|\rho\|+\rho_\iy,\\
&N c_\nu'\ge(N c_\nu')^{\frac{2}{N}}\left[\left(\intR{F(u, v)}\right)^{\frac{2}{2^\ast}}+\|\rho\|^{\frac{2}{2^\ast}}+\rho_\iy^{\frac{2}{2^\ast}}\right].
\end{align*}
}%
Therefore, $\intR{F(u, v)}$, $\|\rho\|$ and $\rho_\iy$ are equal either to $0$ or to $N c_\nu'$. By (\ref{s3eq2-23})-(\ref{s3eq2-24}),
we have $\rho_\iy\le\frac{1}{2} N c_\nu'$, so $\rho_\iy=0$. If $\|\rho\|=N c_\nu'$, then one has that $\intR{F(u, v)}=0$ and so $(u, v)=(0, 0)$. Moreover,
since $\|\mu\|\le N c_\nu'$, we deduce from (\ref{s3eq2-16}) that $\|\mu\|= \left(N c_\nu'\right)^{\frac{2}{N}}\|\rho\|^{\frac{2}{2^\ast}}$. Then Lemma \ref{lemma7} implies
that $\rho$ is concentrated at a single point $z$, and we see from (\ref{s3eq2-23})-(\ref{s3eq2-24}) that
\begin{align*}
\frac{1}{2}N c_\nu'=\lim_{n\to\iy}\sup_{y\in\RN}\int_{B(y, 1)}{F(u_n, v_n)}\ge\lim_{n\to\iy}\int_{B(z, 1)}F(u_n, v_n)=\|\rho\|,
\end{align*}
a contradiction. Therefore, $\intR{F(u, v)}=N c_\nu'$. Since $\|u\|_{\la_1}^2+\|v\|_{\la_2}^2\le N c_\nu'$, we deduce from (\ref{eq3-6}) and (\ref{s3eq2-24}) that
$$N c_\nu'=\|u\|_{\la_1}^2+\|v\|_{\la_2}^2=\intR{F(u, v)}=\lim_{n\to\iy}(\|u_n\|_{\la_2}^2+\|v_n\|_{\la_2}^2),$$
that is, $(u_n, v_n)\to (u, v)$ strongly in $\BD$, $(u, v)\in\N_\nu'$ and $J_\nu(u, v)=c_\nu'$. Recall that $c_\nu'>0$, so $(u, v)\neq (0, 0)$.
By the definition (\ref{eq3-1}) of $\N_\nu'$
and using the Lagrange multiplier method,
it is standard to prove that $J_\nu'(u, v)=0$, so $(u, v)$ is a solution of (\ref{eq0-0}). Moreover, $u, v\ge 0$ are radially symmetric.

Now, assume that $c_\nu'<\frac{1}{N}\min\{S(\la_1)^{N/2}, \,S(\la_2)^{N/2}\}$,
then it is easy to prove that both $u\not\equiv 0$ and $v\not\equiv 0$, that is, $(u, v)\in\N_\nu$, and so
$J_\nu(u, v)=c_\nu'=c_\nu$. Hence, $(u, v)$ is a ground state solution of (\ref{eq0-0}).
By the maximum principle, $u, v> 0$ in $\RN\setminus\{0\}$ and are radially symmetric.
This completes the proof.\hfill$\square$\\

Since $(z_\mu^1, 0)$ and $(0, z_\mu^2)$ belong to $\N_\nu'$, so $c_\nu'\le \frac{1}{N}\min\{S(\la_1)^{\frac{N}{2}}, \,S(\la_2)^{\frac{N}{2}}\}$ always holds.
However,
the following result says that the conclusion $c_\nu'< \frac{1}{N}\min\{S(\la_1)^{\frac{N}{2}}, \,S(\la_2)^{\frac{N}{2}}\}$
can not always hold unfortunately.

\bt\label{th2-1} Assume that $\al, \bb \ge 2$, then there exists $\widetilde{\nu}>0$ such that for all $\nu\in (0, \widetilde{\nu})$ there hold
$$c_\nu'=\frac{1}{N}\min\left\{S(\la_1)^{\frac{N}{2}}, \,S(\la_2)^{\frac{N}{2}}\right\}.$$
Moreover $c_\nu'$ is achieved by and only by
$$\begin{cases}
(0, \pm z_\mu^2),\quad\mu>0,\quad\hbox{if $\la_1<\la_2$},\\
( \pm z_\mu^1, 0),\quad\mu>0,\quad\hbox{if $\la_1>\la_2$},\\
(0, \pm z_\mu^2), \,\,( \pm z_\mu^1, 0),\quad\mu>0,\quad\hbox{if $\la_1=\la_2$}.
\end{cases}$$\et

\noindent {\bf Proof. } Thanks to Lemma \ref{lemma8}, this proof is completely the same as that of \cite[Theorem 3.4]{AFP}, and we omit the details here. \hfill$\square$\\

\noindent {\bf Proof of Theorem \ref{th1} (2)-(3). } Let $\nu>0$. By Lemma \ref{lemma8} and Theorem \ref{th2-1} we know that,
we have to require further assumptions
on $\al, \bb$ and $\nu$ to obtain positive ground state solutions with energy below $\frac{1}{N}\min\{S(\la_1)^{\frac{N}{2}}, \,S(\la_2)^{\frac{N}{2}}\}$.

Denote
$$d_1:=\frac{\La_N-\la_2}{\La_N-\la_1},\quad d_2:=\frac{\La_N-\la_1}{\La_N-\la_2}.$$
Recall (\ref{eq2-31}), we let $\nu>\nu_0$. Then
\be\label{eq2-32}1+\max\{d_1, d_2\}<(2+2^\ast \nu)^{\frac{2}{2^\ast}}.\ee
Without loss of generality, we may assume that $\la_1\le \la_2$. Then (\ref{eq1-5}) yields $S(\la_2)\le S(\la_1)$.
By Hardy inequality (\ref{Hardy}) we have
$\|u\|_{\la_1}^2\le d_2\|u\|_{\la_2}^2$ for all $u\in D^{1, 2}(\RN)$.
Then we deduce from (\ref{eq1-5}), (\ref{eq3-5})  and (\ref{eq2-32}) that
{\allowdisplaybreaks
\begin{align*}c'_\nu
&\le\frac{1}{N}\left[\frac{\|z_\mu^2\|_{\la_1}^2+\|z_\mu^2\|_{\la_2}^2}{\left(\int_{\RN}|z_\mu^2|^{2^\ast}+2^\ast \nu |z_\mu^2|^{2^\ast} +|z_\mu^2|^{2^\ast}\right)^{\frac{2}{2^\ast}}}\right]^{\frac{N}{2}}\\
&\le\frac{1}{N}\left[\frac{1+d_2}{(2+2^\ast \nu)^{\frac{2}{2^\ast}}}\cdot\frac{\|z_\mu^2\|_{\la_2}^2}{\left(\int_{\RN}|z_\mu^2|^{2^\ast}\right)^{\frac{2}{2^\ast}}}\right]^{\frac{N}{2}}\\
&<\frac{1}{N}\left[\frac{\|z_\mu^2\|_{\la_2}^2}{\left(\int_{\RN}|z_\mu^2|^{2^\ast}\right)^{\frac{2}{2^\ast}}}\right]^{\frac{N}{2}}
=\frac{1}{N}S(\la_2)^{\frac{N}{2}}\\
&=\frac{1}{N}\min\left\{S(\la_1)^{\frac{N}{2}}, \,S(\la_2)^{\frac{N}{2}}\right\}.
\end{align*}
}%
Hence, conclusion (2) follows from Lemma \ref{lemma8}.

Repeating the proof of \cite[Theorem 2.2 (iii)-(iv)]{AFP} with minor modifications, we can show that,
if $\al<2$, then for all $\mu>0$, $(0, z_\mu^2)$ is a saddle point for $J_\nu$ in $\N_\nu'$, and so $c'_\nu<\frac{1}{N}S(\la_2)^{N/2}$;
if $\bb<2$, then for all $\mu>0$, $(z_\mu^1, 0)$ is a saddle point for $J_\nu$ in $\N_\nu'$, and so $c'_\nu<\frac{1}{N}S(\la_1)^{N/2}$.
Meanwhile, by (\ref{eq1-5}) we see that, $\la_1<\la_2$ implies $S(\la_1)>S(\la_2)$ and $\la_1>\la_2$ implies $S(\la_1)<S(\la_2)$.
Then under any one conditions of $({\bf C_1})$, $({\bf C_2})$ and $({\bf C_3})$, we have that $c'_\nu<\frac{1}{N}\min\{S(\la_1)^{N/2}, S(\la_2)^{N/2}\}$, and so
conclusion (3) follows from
Lemma \ref{lemma8}. This completes the proof.\hfill$\square$

\vskip0.1in

\s{Proof of Theorem \ref{th2}: The case $N=4$}
\renewcommand{\theequation}{3.\arabic{equation}}

In this section, we assume that $N=4$, $\al=\bb=2$ and $\la_1, \la_2\in (0, 1)$. By Theorem \ref{th2-1} we know that, the ideas of proving Theorem \ref{th1} can not be used here, and we need to
use a different approach, which is much more complicated. This approach will take full use of the fact $\al=\bb=2$,
and so can not be used in the case $N=3$ unfortunately.

\subsection{The special case $\la_1=\la_2=0$}

Consider the following problem
\be\label{dim4}\begin{cases}-\Delta u = u^{3}+2\nu uv^{2}, \quad x\in \R^4,\\
  -\Delta v = v^{3}+2\nu v u^{2},     \quad x\in \R^4,\\
u,\, v\in D^{1, 2}(\R^4),\quad u,\, v>0\,\hbox{in $\R^4$}.\end{cases}\ee

For $\varepsilon > 0$ and $y \in \R^4$, we consider the Aubin-Talenti instanton
\cite{A, T} $U_{\varepsilon,y} \in D^{1,2}(\mathbb{R}^4)$ defined by
\be\label{A-T} U_{\varepsilon,y}(x):= \frac{2\sqrt{2}\varepsilon}{\varepsilon^2 +
|x-y|^2}.\ee  Then $U_{\varepsilon,y}$ satisfies $-\Delta u = u^3$ in $\R^4$ and
\be\label{A-T1}\int_{\mathbb{R}^4}|\nabla U_{\varepsilon,y}|^2\,dx =
\int_{\mathbb{R}^4} |U_{\varepsilon,y}|^{4}\,dx = S^2.\ee Furthermore,
$\{U_{\varepsilon,y} : \varepsilon>0, y \in \R^4\}$ contains all positive solutions
of the equation $-\Delta u = u^3$ in $\R^4$. Note
that (\ref{dim4}) has semi-trivial solutions $(U_{\e,y}, 0)$ and $(0, U_{\e,y})$. Here we are only
interested in nontrivial solutions of (\ref{dim4}), which can be found as nontrivial critical points of the $C^2$ functional $L_\nu: \BD\to \R$, where
\begin{align}\label{s4eq2-2}
L_\nu(u, v)=\frac{1}{2}\left(\| u\|^2+\| v\|^2\right)-\frac{1}{4}\int_{\R^4}\left( u^{4}+4\nu u^2 v^2 + v^{4}\right).
\end{align}

\bd We say a solution $(u_0, v_0)$ of (\ref{dim4}) is a ground state solution if $(u_0, v_0)$ is nontrivial
and $L_\nu (u_0, v_0)\le L_{\nu}(u, v)$ for any other nontrivial solution $(u, v)$ of (\ref{dim4}).\ed

Define the general Nehari manifold of (\ref{s4eq2-2}) as
{\allowdisplaybreaks
\begin{align*}\M_\nu:=\Big\{ &(u, v)\in \BD\,:\, u\not\equiv 0, v\not\equiv 0, \int_{\R^4}|\nabla u|^2
=\int_{\R^4}( u^4+2 \nu u^2 v^2),\\
& \int_{\R^4}|\nabla v|^2=\int_{\R^4}(v^4+2\nu u^2 v^2)\Big\}.\end{align*}
}%
Then any nontrivial solutions of (\ref{dim4}) has to belong to $\M_\nu$. Similarly as $\N_\nu$, we see that $\M_\nu\neq \emptyset$. We set
\be\label{s4eq2-3} m_\nu:=\inf_{(u, v)\in \M_\nu} L_\nu(u, v)=\inf_{(u, v)\in \M_\nu}\frac{1}{4}\int_{\R^4}{\left(|\nabla u|^2+|\nabla v|^2\right)}\,dx.\ee
By Sobolev inequality (\ref{sobolev}), it is easily seen that $m_\nu>0$ for all $\nu$.
Moreover, if $(u_0, v_0)$ is a nontrivial solution of (\ref{dim4}) satisfying $L_\nu (u_0, v_0)=m_\nu$, then $(u_0, v_0)$ is a ground state solution.
Then we have the following result, which will play a crucial role in the proof of Theorem \ref{th2}.
Part of this result comes from the authors' paper \cite{CZ1}.

\bt\label{th3} Let $\nu>0$.
\begin{itemize}

\item[(1)] If $\nu\neq1/2$, then for any $\e>0$, $y\in\R^4$, $((1+2\nu)^{-1/2} U_{\e, y},\,  (1+2\nu)^{-1/2} U_{\e, y})$ is
a positive ground state solution  of (\ref{dim4}), with
\be\label{s4eq2-28}m_\nu=L_\nu \left((1+2\nu)^{-1/2} U_{\e, y},  (1+2\nu)^{-1/2} U_{\e, y}\right)=\frac{1}{2(1+2\nu)}S^{2}.\ee
Moreover, the set $\{(1+2\nu)^{-1/2} U_{\e, y},\,  (1+2\nu)^{-1/2} U_{\e, y}): \e>0, y\in\R^4\}$ contains all positive ground state solutions of (\ref{dim4}).

\item[(2)] If $\nu=1/2$, then for any $\e>0$, $y\in\R^4$, $\theta\in (0, \pi/2)$, $(\sin \theta\, U_{\e, y}, \cos\theta\, U_{\e, y})$ is a ground state solution  of (\ref{dim4}), and $m_{1/2}=\frac{1}{4}S^2$. Moreover, the set $\{(\sin\theta\, U_{\e, y},\, \cos\theta \,U_{\e, y}): \e>0, y\in\R^4, \theta\in (0, \pi/2)\}$ contains all positive ground state solutions of (\ref{dim4}).

\end{itemize}\et

\noindent {\bf Proof. } (1) Let $\nu>0$ and $\nu\neq 1/2$.
Then this result follows directly from \cite[Theorem 1.5 and Theorem 4.1]{CZ1}.

(2) Let $\nu=1/2$. Firstly, note that for any $\theta\in (0, \pi/2)$, $(\sin\theta\,U_{\e, y},\, \cos\theta\, U_{\e, y})$ is a positive solution of (\ref{dim4}) when $\nu=1/2$, so
\be\label{s4eq2-26} m_{1/2}\le L_{1/2}(\sin\theta\,U_{\e, y},\, \cos\theta\, U_{\e, y})=\frac{1}{4}S^{2}.\ee
Secondly, take any $(u, v)\in \M_{1/2}$. If
\be\label{s4eq2-25}\left(\int_{\R^4}{u^2v^2}\,dx\right)^2=\int_{\R^4}{u^4}\,dx\int_{\R^4}{v^4}\,dx,\ee
then by H\"{o}lder inequality, we may assume that $v=C u$ for some constant $C\neq 0$. Recall the definition of $\M_\nu$,
we see from (\ref{sobolev}) that
$$\int_{\R^4}|\nabla u|^2=(1+C^2)\int_{\R^4}u^4\le (1+C^2)S^{-2}\left(\int_{\R^4}|\nabla u|^2\right)^{2},$$
so $\int_{\R^4}|\nabla u|^2\ge (1+C^2)^{-1}S^2$ and then
$$L_{1/2}(u, v)=L_{1/2}(u, C u)\ge\frac{1}{4}S^2.$$
If (\ref{s4eq2-25}) does not hold, then
$$\left(\int_{\R^4}{u^2v^2}\,dx\right)^2<\int_{\R^4}{u^4}\,dx\int_{\R^4}{v^4}\,dx,$$
and it is easy to prove that for any $\nu\in (0, 1/2)$, there exist $t_\nu, s_\nu>0$ such that $(\sqrt{t_\nu}u, \sqrt{s_\nu} v)\in \M_{\nu}$
and $(t_\nu, s_\nu)\to (1, 1)$ as $\nu\to 1/2$. This implies that
$$L_{1/2}(u, v)=\lim_{\nu\nearrow 1/2}L_\nu (\sqrt{t_\nu}u, \sqrt{s_\nu} v)\ge \lim_{\nu\nearrow 1/2} m_\nu=\frac{1}{4}S^2.$$
Therefore, for any $(u, v)\in\M_{1/2}$, we have $L_{1/2}(u, v)\ge\frac{1}{4}S^2$, and so $m_{1/2}\ge \frac{1}{4}S^2$.
Combining this with (\ref{s4eq2-26}), we see that $(\sin\theta\,U_{\e, y},\, \cos\theta\, U_{\e, y})$ is
a ground state solution of (\ref{dim4}) when $\nu=1/2$.

Now assume that $(u, v)$ is any a positive ground state solution of (\ref{dim4}) when $\nu=1/2$. Then
we deduce from (\ref{sobolev}) that
\begin{align*}
S|u|_4^2\le\|u\|^2=|u|_4^4+\int_{\R^4}u^2v^2\,dx\le|u|_4^4+|u|_4^2|v|_4^2,
\end{align*}
that is, $|u|_4^2+|v|_4^2\ge S$. Meanwhile, since $L_{1/2}(u, v)=m_{1/2}=\frac{1}{4}S^2$, we have
$$S^2=\|u\|^2+\|v\|^2\ge S|u|_4^2+ S|v|_4^2,$$
that is, $|u|_4^2+|v|_4^2\le S$. So $|u|_4^2+|v|_4^2= S$, that is,
\begin{align*}
S|u|_4^2=\|u\|^2=|u|_4^4+\int_{\R^4}u^2v^2=|u|_4^4+|u|_4^2|v|_4^2.
\end{align*}
First this means that $v=Cu$ for some $C>0$. Secondly, combining (\ref{sobolev}) with $\|u\|^2=S |u|_4^2$, it is well known that $u=C_1U_{\e, y}$
for some $C_1>0$, $\e>0$ and $y\in\R^4$ (see \cite{A, T} for example). Hence, $(u, v)=(C_1 U_{\e, y}, C_2 U_{\e, y})$ for some $C_1, C_2>0$. Then $L_{1/2}(u, v)=\frac{1}{4}S^2$ yields that
$C_1^2+C_2^2=1$. Therefore, there exists $\theta\in (0, \pi/2)$ such that $C_1=\sin\theta$ and $C_2=\cos\theta$.
This completes the proof.\hfill$\square$

\subsection{The general case $\la_1, \la_2\in (0, 1)$}

Recall the definition (\ref{th2constant}) of $\nu_1$, we have the following important energy estimate,
and the idea of the proof comes from the authors' paper \cite{CZ1}.

\bl\label{lemma3-1}For any $\nu\in (0, \nu_1)$, there holds
$$c_\nu<\min\left\{\frac{1}{4}S(\la_1)^2+\frac{1}{4}S(\la_2)^2,\quad m_\nu\right\}.$$\el

\noindent{\bf Proof. } Define
\be\label{equa1}G(u, v):=\left(
                           \begin{array}{cc}
                             \int_{\R^4}u^4\,dx & 2\nu\int_{\R^4}u^2v^2\,dx \\
                             2\nu\int_{\R^4}u^2v^2\,dx & \int_{\R^4}v^4\,dx \\
                           \end{array}
                         \right).
\ee
When $\det G(u, v)>0$, the inverse matrix of $G(u, v)$ is
\be\label{equa2}G^{-1}(u, v):=\frac{1}{\det G(u, v)}\left(
                           \begin{array}{cc}
                             \int_{\R^4}v^4\,dx & -2\nu\int_{\R^4}u^2v^2\,dx \\
                             -2\nu\int_{\R^4}u^2v^2\,dx & \int_{\R^4}u^4\,dx \\
                           \end{array}
                         \right).
\ee
Assume $\nu\in\big(0,\,\nu_1)$. Obviously, one has that $2\nu<1$ and so $ \det G(z_1^1, z_1^2)>0$.
Recall that $\|z_1^i\|_{\la_i}^2=|z_1^i|_4^4=4 S(\la_i)^2, i=1, 2$,
we see that $(\sqrt{t_0}z_1^1, \sqrt{s_0}z_1^2)\in\mathcal{N}_\nu$ for some $t_0>0, s_0>0$  is equivalent to
{\allowdisplaybreaks
\begin{align}\label{equa3}
&\left(
  \begin{array}{c}
    t_0 \\
    s_0\\
  \end{array}
\right)
:=G^{-1}(z_1^1, z_1^2)\left(
                               \begin{array}{c}
                                 |z_1^1|_4^4 \\
                                 |z_1^2|_4^4 \\
                               \end{array}
                             \right)\nonumber\\
&=\frac{1}{\det G(z_1^1, z_1^2)}\left(
                                          \begin{array}{c}
                                            |z_1^2|_4^4\big(|z_1^1|_4^4-2\nu\int_{\R^4} (z_1^1)^2 (z_1^2)^2\big)\\
                                            |z_1^1|_4^4\big(|z_1^2|_4^4-2\nu\int_{\R^4} (z_1^1)^2 (z_1^2)^2\big)\\
                                          \end{array}
                                        \right)>\left(
                                                  \begin{array}{c}
                                                    0 \\
                                                    0 \\
                                                  \end{array}
                                                \right).
\end{align}
}%
Here and in the following, $\left(
                              \begin{array}{c}
                                a \\
                                b \\
                              \end{array}
                            \right)>\left(
                                                  \begin{array}{c}
                                                    0 \\
                                                    0 \\
                                                  \end{array}
                                                \right)
$ means both $a>0$ and $b>0$.
Meanwhile, (\ref{eq1-5}) yields $S_{\la_i}=(1-\la_i)^{3/4} S$, so we deduce from (\ref{th2constant}) that
$$2\nu<\min\left\{\frac{1-\la_1}{1-\la_2},\,\,\frac{1-\la_2}{1-\la_1}\right\}\le\min\left\{\frac{S(\la_1)}{S(\la_2)},\,\,\frac{S(\la_2)}{S(\la_1)}\right\}.$$
Then
\begin{align*}
2\nu\int_{\R^4} (z_1^1)^2 (z_1^2)^2< \min\left\{\frac{S(\la_1)}{S(\la_2)},\,\frac{S(\la_2)}{S(\la_1)}\right\}\left|z_1^1\right|_4^2 \left|z_1^2\right|_4^2
=\min\left\{\left|z_1^1\right|_4^4, \left|z_1^2\right|_4^4\right\}.
\end{align*}
So (\ref{equa3}) holds
and $(\sqrt{t_0}z_1^1,\sqrt{s_0}z_1^2)\in \mathcal{N}_\nu$ for $(t_0, s_0)$ defined in (\ref{equa3}). Then
{\allowdisplaybreaks
\begin{align*}
c_\nu &\le J_\nu(\sqrt{t_0}z_{1}^1,\sqrt{s_0} z_1^2)
=\frac{t_0}{4}\|z_1^1\|_{\la_1}^2+\frac{s_0}{4}\|z_1^2\|_{\la_2}^2\\
&=\frac{t_0}{4}\int_{\R^4} (z_1^1)^4
        +\frac{s_0}{4}\int_{\R^4}(z_1^2)^4\\
&<\frac{t_0}{4}\int_{\R^4}\Big((z_1^1)^4+2\nu (z_1^1)^2 (z_1^2)^2\Big)+
        \frac{s_0}{4}\int_{\R^4}\Big((z_1^2)^4+2\nu (z_1^1)^2 (z_1^2)^2\Big)\\
&=\frac{1}{4}\|z_1^1\|_{\la_1}^2+\frac{1}{4}\|z_1^2\|_{\la_2}^2
= \frac{1}{4}S(\la_1)^2+\frac{1}{4}S(\la_2)^2.
\end{align*}
}%
Hence $c_\nu< \frac{1}{4}S(\la_1)^2+\frac{1}{4}S(\la_2)^2$. It remains to prove $c_\nu < m_\nu$. Take $y_0\in\R^4$ such that $|y_0|=2$.
Let $\psi\in C_0^\iy(B(y_0, 1), \R)$ be a function with $0\le \psi\le 1$, $\psi\equiv 1$ for $x\in B(y_0, 1/2)$.
Recall $U_{\e, y_0}$ in (\ref{A-T}) and (\ref{A-T1}), we define $U_\e:=\psi U_{\e, y_0}$. Then by \cite{BN} or \cite[Lemma 1.46]{W},
we have the following inequalities
{\allowdisplaybreaks
\begin{gather*}
    \int_{\R^4}|\nabla U_\e|^2=S^2+O(\e^2),\quad\int_{\R^4}|U_\e|^4=S^2+O(\e^4),\\
    \int_{\R^4}\frac{|U_\e|^2}{|x|^2}\,dx\ge\frac{1}{9}\int_{B(y_0, 1)}|U_\e|^2\ge C \e^2|\ln \e|+O(\e^2),
\end{gather*}
}%
where $C$ is a positive constant. Recalling that $\la_1, \la_2>0$, we have
{\allowdisplaybreaks
\begin{align}\label{s4eq13}
J_\nu(\sqrt{t}U_\e, \sqrt{s}U_\e)&=\frac{1}{2}t\int_{\R^4}\left(|\nabla U_\e|^2-\frac{\la_1}{|x|^2} U_\e^2\right)
   +\frac{1}{2}s\int_{\R^4}\left(|\nabla U_\e|^2-\frac{\la_2}{|x|^2} U_\e^2\right)\nonumber\\
&\quad-\frac{1}{4}(t^2+4\nu t s+s^2)\int_{\R^4} U_\e^4\nonumber\\
&\le \frac{1}{2}\left(t+s\right)\left(S^2- C\e^2|\ln \e|+ O(\e^2)\right)\nonumber\\
&\quad-\frac{1}{4}\left(t^2+4\nu t s+s^2\right)\left(S^2 +O(\e^4)\right).
\end{align}
}%
Denote
$$A_\e=S^2- C\e^2|\ln \e |+ O(\e^2),\quad B_\e=S^2 +O(\e^4),$$
then $0<A_\e<B_\e$ and $A_\e< S^2$ for $\e>0$ small enough. Consider
$$f_\e(t, s):=\frac{1}{2}A_\e( t+s)-\frac{1}{4}B_\e(t^2+4\nu t s+s^2),$$
then it is easy to see that there exists $t_\e, s_\e>0$ such that
$$f_\e(t_\e, s_\e)=\max_{t, s>0}f_\e(t,s).$$
By
$\frac{\partial}{\partial t}f_\e(t, s)|_{(t_\e, s_\e)}
=\frac{\partial}{\partial s}f_\e(t, s)|_{(t_\e, s_\e)}=0$, we see that
$$t_\e=s_\e=\frac{A_\e}{(1+2\nu)B_\e}.$$
Then it follows from (\ref{s4eq2-28}) and (\ref{s4eq13}) that
{\allowdisplaybreaks
\begin{align}\label{equa7}
\max_{t,s>0}J_\nu(\sqrt{t}U_\e, \sqrt{s}U_\e)&\le\max_{t, s>0}f_\e(t,s)=f_\e(t_\e, s_\e)\nonumber\\
&=\frac{1}{2(1+2\nu)}\frac{A_\e^2}{B_\e}<\frac{1}{2(1+2\nu)}A_\e\nonumber\\
&<\frac{S^2}{2(1+2\nu)}=m_\nu\quad\hbox{holds for $\e$ small enough.}
\end{align}
}%
Similarly as above, we have $\det G(U_{\e}, U_{\e})>0$. Moreover,
$(\sqrt{\tilde{t}_\e}U_{\e}, \sqrt{\tilde{s}_\e}U_{\e})\in\mathcal{N}_\nu$ for some $\tilde{t}_\e>0, \tilde{s}_\e>0$  is equivalent to
{\allowdisplaybreaks
\begin{align}\label{equa4}
\left(
  \begin{array}{c}
    \tilde{t}_\e\\
    \tilde{s}_\e\\
  \end{array}
\right)
&=\frac{|U_\e|_4^4}{\det G(U_\e, U_\e)}\left(
                                          \begin{array}{c}
                                            \|U_\e\|_{\la_1}^2-2\nu\|U_\e\|_{\la_2}^2\\
                                             \|U_\e\|_{\la_2}^2-2\nu\|U_\e\|_{\la_1}^2\\
                                          \end{array}
                                        \right)
                                        >\left(
                                                  \begin{array}{c}
                                                    0 \\
                                                    0 \\
                                                  \end{array}
                                                \right).
\end{align}
}%
On the other hand, by (\ref{Hardy}) we have
{\allowdisplaybreaks
\begin{align*}
\|U_\e\|_{\la_1}^2-2\nu\|U_\e\|_{\la_2}^2
=&(1-2\nu)\int_{\R^4}|\nabla U_\e|^2-(\la_1-2\nu\la_2)\int_{\R^4}\frac{U_\e^2}{|x|^2}\\
\ge&(1-2\nu)\int_{\R^4}\frac{U_\e^2}{|x|^2}-(\la_1-2\nu\la_2)\int_{\R^4}\frac{U_\e^2}{|x|^2}\\
=&\Big[(1-\la_1)-2\nu(1-\la_2)\Big]\int_{\R^4}\frac{U_\e^2}{|x|^2}>0.
\end{align*}
}%
Similarly, $\|U_\e\|_{\la_2}^2-2\nu\|U_\e\|_{\la_1}^2>0$.
Hence, (\ref{equa4}) holds and $(\sqrt{\tilde{t}_{\e}}U_\e, \sqrt{\tilde{s}_{\e}}U_\e)\in \mathcal{N}_\nu$ for $(\tilde{t}_\e, \tilde{s}_\e)$ defined in (\ref{equa4}). Then we see from (\ref{equa7}) that
$$c_\nu\le J_\nu\left(\sqrt{\tilde{t}_{\e}}U_\e, \sqrt{\tilde{s}_{\e}}U_\e\right)\le \max_{t,s>0}J_\nu(\sqrt{t}U_\e, \sqrt{s}U_\e)<m_\nu.$$
This completes the proof.\hfill$\square$

\bl\label{lemma3-2}Assume that $\nu\in (0,\,\nu_1)$.
Then there exist $C_2>C_1>0$
such that for any $(u, v)\in \mathcal{N}_\nu$ with $J_\nu(u, v)\le \frac{1}{4}S(\la_1)^{2}+\frac{1}{4}S(\la_2)^{2}$, there holds
\be\label{eq4-0} C_1\le\int_{\R^4} u^4\,dx,\,\, \int_{\R^4}v^4\,dx\le C_2.\ee\el

\noindent{\bf Proof. } Take any $(u, v)\in \mathcal{N}_\nu$ with $J_\nu(u, v)\le \frac{1}{4}S(\la_1)^{2}+\frac{1}{4}S(\la_2)^{2}$. By (\ref{eq1-5-1})
and H\"{o}lder inequality, one has
{\allowdisplaybreaks
\begin{align*}
    S(\la_1)|u|_4^2\le\|u\|_{\la_1}^2
    &=\int_{\R^4}( u^4+2\nu u^2 v^2)\le |u|_4^4+2\nu |u|_4^2|v|_4^2,\\
    S(\la_2)|v|_4^2\le\|v\|_{\la_2}^2
    &=\int_{\R^4}( v^4+2\nu u^2 v^2)\le |v|_4^4+2\nu |u|_4^2|v|_4^2,
\end{align*}
}%
Therefore, there exists $C_2>0$ such that $\int_{\R^4} u^4, \int_{\R^4}v^4\le C_2$.
Moreover,
{\allowdisplaybreaks
\begin{align}
\label{eq15} &|u|_4^2+2\nu |v|_4^2\ge S(\la_1),\\
\label{eq17} &2\nu |u|_4^2+|v|_4^2\ge S(\la_2),\\
\label{eq18} &S(\la_1)|u|_4^2+S(\la_2)|v|_4^2\le S(\la_1)^{2}+S(\la_2)^{2}.
\end{align}
}%
Recall that $S(\la_i)=(1-\la_i)^{3/4}S$. Since
 $\nu\in (0, \nu_1)$ and $\nu_1$ is defined in (\ref{th2constant}), by (\ref{eq15}) and (\ref{eq18}) we have
\be\label{eq4-14}|u|_4^2\ge\frac{S(\la_1)S(\la_2)-2\nu\left(S(\la_1)^2+S(\la_2)^2\right)}{S(\la_2)-2\nu S(\la_1)}>0,\ee
and by (\ref{eq17}) and (\ref{eq18}) we have
$$|v|_4^2\ge\frac{S(\la_1)S(\la_2)-2\nu\left(S(\la_1)^2+S(\la_2)^2\right)}{S(\la_1)-2\nu S(\la_2)}>0.$$
This completes the proof.\hfill$\square$\\

The following lemma is motivated by \cite{Terracini}, and some ideas of the proof come from \cite{S1}.

\bl\label{lemma3-3} Assume that $\nu\in (0, \nu_1)$. Let $(u_n, v_n)\in\N_\nu$ be a minimizing
sequence of $c_\nu$, and $(u_n, v_n)\rightharpoonup (0, 0)$ weakly in $\BD$.
Then for any $r>0$ and for every $\e\in (-r, 0)\cup(0, r)$, there exists $\rho\in (\e, 0)\cup (0, \e)$ such that,
up to a subsequence,
\be\label{eq400}\hbox{either}\,\,\int_{B_{r+\rho}}(|\nabla u_n|^2+|\nabla v_n|^2)\to 0
\,\,\hbox{or}\,\,\int_{\RN\setminus B_{r+\rho}}(|\nabla u_n|^2+|\nabla v_n|^2)\to 0.\ee\el

\noindent{\bf Proof. } Without loss of generality, we only consider the case $\e\in (0, r)$
(the proof for the case $\e\in (-r, 0)$ is similar). Since $(u_n, v_n)\in\N_\nu$ is a minimizing
sequence of $c_\nu$, then $(u_n, v_n)$ are uniformly bounded in $\BD$.
Moreover, by Lemmas \ref{lemma3-1} and \ref{lemma3-2} we may assume that $(u_n, v_n)$ satisfies (\ref{eq4-0}) for all $n\in\mathbb{N}$.\\

\noindent{\bf Step 1.} We prove (\ref{eq400}) by further assuming that $J_\nu'(u_n, v_n)\to 0$ as $n\to\iy$.

In the following, some arguments are borrowed from \cite{S1} (see also \cite[Lemma III.3.3]{S2} or \cite[Proposition 5.2]{Terracini}).
Denote $\BS$ as the unit sphere of $\R^4$. Since
$$\int_r^{r+\e}\,d\rho\int_{\rho\BS}\left(|\nabla u_n|^2+|\nabla v_n|^2\right)=\int_{r\le |x|\le r+\e}\left(|\nabla u_n|^2+|\nabla v_n|^2\right)$$
is bounded, we can find $\rho\in (0, \e)$ such that
$$\int_{(r+\rho)\BS}\left(|\nabla u_n|^2+|\nabla v_n|^2\right)\le\frac{3}{\e}\int_{r\le |x|\le r+\e}\left(|\nabla u_n|^2+|\nabla v_n|^2\right)$$
holds for infinitely many $n$'s.
Therefore, as $H^1((r+\rho)\BS)$ is compactly embedded into $H^{1/2}((r+\rho)\BS)$, up to a subsequence we can assume that $u_n\to u, v_n\to v$
strongly in $H^{1/2}((r+\rho)\BS)$. On the other hand, by the continuity of the embedding $H^1(B_{r+\rho})\hookrightarrow H^{1/2}((r+\rho)\BS)$ and by
the weak convergence to $(0, 0)$ of $(u_n, v_n)$, we deduce that $(u, v)=(0, 0)$,
that is, $u_n\to 0$ and $v_n\to 0$
strongly in $H^{1/2}((r+\rho)\BS)$. Let $w_{i, n}, i=1, 2$ be the solutions to the Dirichlet problems
\be\label{5eq10}\begin{cases}\Delta w_{1, n}=0 &\hbox{in}\,\, B_{r+\e}\setminus B_{r+\rho}\\
w_{1, n}=0 &\hbox{on}\,\,(r+\e)\BS\\
w_{1, n}=u_n &\hbox{on}\,\,(r+\rho)\BS\end{cases},\quad \begin{cases}\Delta w_{2, n}=0 &\hbox{in}\,\, B_{r+\rho}\setminus B_{r-\e}\\
w_{2, n}=0 &\hbox{on}\,\,(r-\e)\BS\\
w_{2, n}=u_n &\hbox{on}\,\,(r+\rho)\BS\end{cases},\ee
and let $\sg_{i, n}, i=1, 2$ be the solutions to the Dirichlet problems
\be\label{5eq11}\begin{cases}\Delta \sg_{1, n}=0 &\hbox{in}\,\, B_{r+\e}\setminus B_{r+\rho}\\
\sg_{1, n}=0 &\hbox{on}\,\,(r+\e)\BS\\
\sg_{1, n}=v_n &\hbox{on}\,\,(r+\rho)\BS\end{cases},\quad \begin{cases}\Delta \sg_{2, n}=0 &\hbox{in}\,\, B_{r+\rho}\setminus B_{r-\e}\\
\sg_{2, n}=0 &\hbox{on}\,\,(r-\e)\BS\\
\sg_{2, n}=v_n &\hbox{on}\,\,(r+\rho)\BS\end{cases}.\ee
By the continuity of the inverse Laplace operator from $H^{1/2}(\partial\Om)$ to $H^1(\Om)$, it follows from the above discussion
that $w_{1, n}\to 0,\, \sg_{1, n}\to 0$ strongly in $H^1(B_{r+\e}\setminus B_{r+\rho})$
and $w_{2, n}\to 0,\, \sg_{2, n}\to 0$ strongly in $H^1(B_{r+\rho}\setminus B_{r-\e})$. Define
{\allowdisplaybreaks
\begin{align}
\label{5eq12}u_{1, n}(x)=\begin{cases}u_n(x) &\hbox{if}\,\, x\in B_{r+\rho}\\
w_{1, n} &\hbox{if}\,\,x\in B_{r+\e}\setminus B_{r+\rho}\\
0 &\hbox{elsewhere},\end{cases}\\
v_{1, n}(x)=\begin{cases}v_n(x) &\hbox{if}\,\, x\in B_{r+\rho}\\
\sg_{1, n} &\hbox{if}\,\,x\in B_{r+\e}\setminus B_{r+\rho}\\
0 &\hbox{elsewhere},\end{cases}\\
u_{2, n}(x)=\begin{cases}0 &\hbox{if}\,\, x\in B_{r-\e}\\
w_{2, n} &\hbox{if}\,\,x\in B_{r+\rho}\setminus B_{r-\e}\\
u_n(x) &\hbox{elsewhere},\end{cases}\\
\label{5eq13}v_{2, n}(x)=\begin{cases}0 &\hbox{if}\,\, x\in B_{r-\e}\\
\sg_{2, n} &\hbox{if}\,\,x\in B_{r+\rho}\setminus B_{r-\e}\\
v_n(x) &\hbox{elsewhere}.\end{cases}
\end{align}
}%
Then it is easy to see that
{\allowdisplaybreaks
\begin{align}\label{eq4-2}&\|u_n\|_{\la_1}^2=\|u_{1, n}\|_{\la_1}^2+\|u_{2,n}\|_{\la_1}^2+o(1),\\
\label{eq4-3}&\|v_n\|_{\la_2}^2=\|v_{1, n}\|_{\la_2}^2+\|v_{2,n}\|_{\la_2}^2+o(1).\end{align}
}%
Moreover, we can easily obtain
{\allowdisplaybreaks
\begin{align}
\label{eq4-4} &J_\nu'(u_{1, n}, v_{1, n})(u_{1, n}, 0)=J_\nu'(u_{n}, v_{n})(u_{1, n}, 0)+o(1)=o(1),\\
\label{eq4-5} &J_\nu'(u_{1, n}, v_{1, n})(0, v_{1, n})=J_\nu'(u_{n}, v_{n})( 0, v_{1, n})+o(1)=o(1),\\
\label{eq4-6} &J_\nu'(u_{2, n}, v_{2, n})(u_{2, n}, 0)=J_\nu'(u_{n}, v_{n})(u_{2, n}, 0)+o(1)=o(1),\\
\label{eq4-7} &J_\nu'(u_{2, n}, v_{2, n})(0, v_{2, n})=J_\nu'(u_{n}, v_{n})( 0, v_{2, n})+o(1)=o(1).
\end{align}
}%
Then we claim that
\be\label{eq4-8}\hbox{either}\,\,\lim_{n\to\iy}(\|u_{1, n}\|^2+\|v_{1, n}\|^2)=0\,\,\hbox{or}\,\,\lim_{n\to\iy}(\|u_{2, n}\|^2+\|v_{2, n}\|^2)=0.\ee

In fact, if (\ref{eq4-8}) does not hold, then up to a subsequence,
\be\label{eq44-8}\hbox{both}\,\,\lim_{n\to\iy}(\|u_{1, n}\|^2+\|v_{1, n}\|^2)>0\,\,\hbox{and}\,\,\lim_{n\to\iy}(\|u_{2, n}\|^2+\|v_{2, n}\|^2)>0.\ee
We have the following several cases.

{\bf Case 1.} Up to a subsequence, both $\lim\limits_{n\to\iy}\|u_{1, n}\|^2>0$ and $\lim\limits_{n\to\iy}\|v_{1, n}\|^2>0$.

Since norms $\|\cdot\|_{\la_i}, i=1, 2$ are equivalent to $\|\cdot\|$, and (\ref{eq4-4})-(\ref{eq4-5}) yield
{\allowdisplaybreaks
\begin{align}
\label{eq4-9}\|u_{1, n}\|_{\la_1}^2=|u_{1, n}|_4^4 + 2\nu \int_{\R^4}u_{1, n}^2 v_{1, n}^2+o(1),\\
\label{eq4-10}\|v_{1, n}\|_{\la_2}^2=|v_{1, n}|_4^4 + 2\nu \int_{\R^4}u_{1, n}^2 v_{1, n}^2+o(1).
\end{align}
}%
Hence, both $\liminf\limits_{n\to\iy}|u_{1, n}|_4^4>0$ and $\liminf\limits_{n\to\iy}|v_{1, n}|_4^4>0$. Since $2\nu<2\nu_1\le 1$,
so by H\"{o}lder inequality we have
$$\liminf_{n\to\iy} \left[|u_{1, n}|_4^4|v_{1, n}|_4^4-\left(2\nu \int_{\R^4}u_{1, n}^2 v_{1, n}^2\right)^2\right]>0.$$
Combining this with (\ref{eq4-9})-(\ref{eq4-10}), it is easy to prove that there exist $t_n, s_n>0$ such that $(\sqrt{t_n} u_{1, n},\,\sqrt{s_n} v_{1, n})\in\N_\nu$ and $(t_n, s_n)\to (1, 1)$, and so we conclude from (\ref{eq4-2})-(\ref{eq4-3}) and (\ref{eq44-8}) that
{\allowdisplaybreaks
\begin{align*}
c_\nu &=\lim_{n\to\iy}J_\nu(u_n, v_n)= \lim_{n\to\iy}\frac{1}{4}(\|u_n\|_{\la_1}^2+\|v_n\|_{\la_2}^2)\\
&=\lim_{n\to\iy}\frac{1}{4}(t_n\|u_{1, n}\|_{\la_1}^2+s_n\|v_{1,n}\|_{\la_2}^2)+\lim_{n\to\iy}\frac{1}{4}(\|u_{2, n}\|_{\la_1}^2+\|v_{2, n}\|_{\la_2}^2)\\
&>\lim_{n\to\iy}\frac{1}{4}(t_n\|u_{1, n}\|_{\la_1}^2+s_n\|v_{1,n}\|_{\la_2}^2)\\
&=\lim_{n\to\iy}J_\nu(\sqrt{t_n} u_{1, n},\,\sqrt{s_n} v_{1, n})\ge c_\nu,
\end{align*}
}%
a contradiction. So Case 1 is impossible.

{\bf Case 2.} Up to a subsequence, $\lim\limits_{n\to\iy}\|u_{1, n}\|^2=0$ and $\lim\limits_{n\to\iy}\|v_{1, n}\|^2>0$.

Then (\ref{eq4-10}) yields that
$$\|v_{1, n}\|_{\la_2}^2=|v_{1, n}|_4^4+o(1)\le S(\la_2)^{-2}\|v_{1, n}\|_{\la_2}^4+o(1), $$
and so $\lim\limits_{n\to\iy}\|v_{1, n}\|_{\la_2}^2\ge S(\la_2)^2$. By (\ref{eq4-0}) we have
$$\liminf_{n\to\iy} \|u_{2, n}\|^2= \liminf_{n\to\iy} \|u_{n}\|^2-\lim_{n\to\iy}\|u_{1, n}\|^2>0.$$
If up to a subsequence, $\lim_{n\to\iy} \|v_{2, n}\|^2>0$, then we can get a contradiction just as Case 1. Therefore, $\lim_{n\to\iy} \|v_{2, n}\|^2=0$.
Then similarly as above, we can deduce from (\ref{eq4-6}) that $\lim\limits_{n\to\iy}\|u_{2, n}\|_{\la_1}^2\ge S(\la_1)^2$. Then
{\allowdisplaybreaks
\begin{align*}
c_\nu &= \lim_{n\to\iy}\frac{1}{4}(\|u_n\|_{\la_1}^2+\|v_n\|_{\la_2}^2)\\
&=\lim_{n\to\iy}\frac{1}{4}(\|u_{2, n}\|_{\la_1}^2+\|v_{1,n}\|_{\la_2}^2)
\ge\frac{1}{4}\left(S(\la_1)^2+ S(\la_2)^2\right),
\end{align*}
}%
a contradiction with Lemma \ref{lemma3-1}. So Case 2 is impossible.

{\bf Case 3.} Up to a subsequence, $\lim\limits_{n\to\iy}\|u_{1, n}\|^2>0$ and $\lim\limits_{n\to\iy}\|v_{1, n}\|^2=0$.

By a similar argument as Case 2, we get a contradiction. So Case 3 is impossible.

Since none of Cases 1, 2 and 3 is true, we see that (\ref{eq44-8}) is impossible, that is, (\ref{eq4-8}) holds.
Recall the definition of $(u_{i, n}, v_{i, n})$, (\ref{eq400}) follows directly from (\ref{eq4-8}). This completes the proof of Step 1.\\

\noindent{\bf Step 2.}  We prove (\ref{eq400}) without assuming that $J_\nu'(u_n, v_n)\to 0$ as $n\to\iy$.

By the Ekeland variational principle (see \cite[Theorem 5.1]{S2} for example), there exists a sequence $\{(\wu_n, \wv_n)\}\in\N_\nu$ such that
{\allowdisplaybreaks
\begin{gather}
\label{s4eq20}J_\nu(\wu_n, \wv_n)\le J_\nu (u_n, v_n),\quad \|(u_n, v_n)-(\wu_n, \wv_n)\|\le \frac{1}{n},\\
\label{s4eq21}J_\nu(u, v)\ge J_\nu (\wu_n, \wv_n)-\frac{1}{n}\|(\wu_n, \wv_n)-(u, v)\|,\quad \forall (u, v)\in\mathcal{N}_\nu.
\end{gather}
}%
Here, $\|(u, v)\|:=(\int_{\R^4}(|\nabla u|^2+|\nabla v|^2)\,dx)^{1/2}$ is also a norm of $\BD$, which is equivalent to $\|(u, v)\|_{\BD}$.
Recall that $(u_n, v_n)\rightharpoonup (0, 0)$ weakly in $\BD$, by (\ref{s4eq20}) we also have $J_\nu(\wu_n, \wv_n)\to c_\nu$ and
$(\wu_n, \wv_n)\rightharpoonup (0, 0)$ weakly in $\BD$. Moreover, by Lemma \ref{lemma3-1} we may assume that $(\wu_n, \wv_n)$ satisfies (\ref{eq4-0}) for all $n\in\mathbb{N}$. Then by repeating the proof of \cite[Theorem 1.3 (1)-(2)]{CZ1}, we can prove that $J'_\nu(\wu_n, \wv_n)\to 0$ as $n\to\iy$.
Hence Step 1 yields that (\ref{eq400}) holds for $(\wu_n, \wv_n)$. Combining this with (\ref{s4eq20}), we see that (\ref{eq400}) holds for $(u_n, v_n)$. This completes the proof.\hfill$\square$\\

\noindent {\bf Proof of Theorem \ref{th2}. } Fix any $\nu\in (0, \nu_1)$.
Take a sequence $(\ou_n, \ov_n)\in\mathcal{N}_\nu$ such
that $J_\nu(\ou_n, \ov_n) \to c_\nu$ as $n\to\iy$. Recall that
$E(u, v)= |\nabla u|^2+|\nabla v|^2-\frac{\la_1}{|x|^2}|u|^2-\frac{\la_2}{|x|^2}|v|^2$, there exists $R_n>0$ such that
$$\int_{B_{R_n}} E(\ou_n, \ov_n)=\int_{\RN\setminus B_{R_n}} E(\ou_n, \ov_n)=\frac{1}{2}(\|\ou_n\|_{\la_1}^2+\|\ov_n\|_{\la_2}^2).$$
Define
$$(\wu_n(x), \wv_n(x)):=\left(R_n^{\frac{N-2}{2}}\ou_n(R_n x), \, R_n^{\frac{N-2}{2}}\ov_n(R_n x)\right),$$
Then by a direct computation, we see that $(\wu_n, \wv_n)\in\N_\nu$ and $J_\nu(\wu_n, \wv_n) \to c_\nu$. Moreover,
\be\label{eq4-11}\int_{B_{1}} E(\wu_n, \wv_n)=\int_{\RN\setminus B_{1}} E(\wu_n, \wv_n)=\frac{1}{2}(\|\wu_n\|_{\la_1}^2+\|\wv_n\|_{\la_2}^2)\to 2 c_\nu>0.\ee
By the Ekeland variational principle (see \cite[Theorem 5.1]{S2} for example), there exists a sequence $\{(u_n, v_n)\}\in\N_\nu$ such that
{\allowdisplaybreaks
\begin{gather}
\label{s4eq22}J_\nu(u_n, v_n)\le J_\nu (\wu_n, \wv_n),\quad \|(u_n, v_n)-(\wu_n, \wv_n)\|\le \frac{1}{n},\\
\label{s4eq23}J_\nu(u, v)\ge J_\nu (u_n, v_n)-\frac{1}{n}\|(u_n, v_n)-(u, v)\|,\quad \forall (u, v)\in\mathcal{N}_\nu.
\end{gather}
}%
Similarly as Step 2 in the proof of Lemma \ref{lemma3-3}, we have that $J_\nu(u_n, v_n)\to c_\nu$ and $J_\nu'(u_n, v_n)\to 0$ as $n\to\iy$.
Moreover, (\ref{eq4-11}) and (\ref{s4eq22}) yield that
{\allowdisplaybreaks
\begin{align}\label{eq4-12}&\lim_{n\to\iy}\int_{B_{1}} E(u_n, v_n)=\lim_{n\to\iy}\int_{B_{1}} E(\wu_n, \wv_n)=2c_\nu,\\
\label{eq4-12-1}&\lim_{n\to\iy}\int_{\RN\setminus B_{1}} E(u_n, v_n)=\lim_{n\to\iy}\int_{\RN\setminus B_{1}} E(\wu_n, \wv_n)=2 c_\nu.
\end{align}
}%

Note that $(u_n, v_n)$ are uniformly bounded in $\BD$. Then up to a subsequence, we assume that $(u_n, v_n)\rightharpoonup (u, v)$ weakly in $\BD$.
Then $J_\nu'(u, v)=0$.

{\bf Case 1.} $(u, v)\equiv (0, 0)$.

Then we can apply twice Lemma \ref{lemma3-3} with $r=1$ and $\e=\pm 1/4$ respectively, and there exist $\rho^+\in (0, 1/4)$ and $\rho^-\in (-1/4, 0)$
such that the alternative (\ref{eq400}) holds. By (\ref{eq4-12})-(\ref{eq4-12-1}) we can rule out all possibilities other than
\be\label{eq4-13}\int_{B_{1+\rho^-}}(|\nabla u_n|^2+|\nabla v_n|^2)\to 0
\,\,\hbox{and}\,\,\int_{\RN\setminus B_{1+\rho^+}}(|\nabla u_n|^2+|\nabla v_n|^2)\to 0.\ee
Now let $\eta\in C^\iy_0(\R^4)$ such that $0\le\eta\le 1$, $\eta(x)=1$ for $|x|\in[3/4, 5/4]$ and $\eta(x)=0$ for $|x|\not\in [1/2, 3/2]$.
Recall that $(u_n, v_n)\rightharpoonup (0, 0)$ weakly in $\BD$, so $u_n, v_n\to 0$ strongly in $L_{\loc}^2(\R^4)$. Combining this with (\ref{eq4-13}), we
obtain that
$$\|(\eta u_n)-u_n\|\to 0,\quad\|(\eta v_n)-v_n\|\to 0,\quad\hbox{as $n\to\iy$}.$$
By Hardy inequality (\ref{Hardy}), we have
{\allowdisplaybreaks
\begin{align*}\int_{\R^4}\frac{u_n^2}{|x|^2}&=\int_{\R^4}\frac{(1-\eta)^2 u_n^2}{|x|^2}+\int_{\R^4}\frac{(2-\eta)\eta u_n^2}{|x|^2}\\
&\le \|(\eta u_n)-u_n\|^2+8 \int_{1/2\le|x|\le 3/2}u_n^2=o(1),
\end{align*}
}%
Similarly, $\int_{\R^4}\frac{v_n^2}{|x|^2}=o(1)$. Therefore, we see from $(u_n, v_n)\in \N_\nu$ that
{\allowdisplaybreaks
\begin{align*}
&\int_{\R^4}|\nabla u_n|^2=\int_{\R^4}u_n^4+2\nu\int_{\R^4}u_n^2 v_n^2+ o(1),\\
&\int_{\R^4}|\nabla v_n|^2=\int_{\R^4}v_n^4+2\nu\int_{\R^4}u_n^2 v_n^2+ o(1).
\end{align*}
}%
From Lemma \ref{lemma3-2} we may assume that $|u_n|_4^4, |v_n|_4^4\ge C>0$, where $C$ is independent of $n$. Since $2\nu<2\nu'\le1$, then it is easy to prove that there exist $t_n, s_n>0$ such that $(\sqrt{t_n}u_n, \sqrt{s_n} v_n)\in\M_\nu$ and $(t_n, s_n)\to (1, 1)$ as $n\to\iy$. Hence,
{\allowdisplaybreaks
\begin{align*}
c_\nu &=\lim_{n\to\iy}J_\nu (u_n, v_n)=\lim_{n\to\iy}\frac{1}{4}(\|u_n\|_{\la_1}^2+\|v_n\|_{\la_2}^2)\\
&=\lim_{n\to\iy}\frac{1}{4}\left(t_n\int_{\R^4}|\nabla u_n|^2 +s_n\int_{\R^4}|\nabla v_n|^2\right)\\
&=\lim_{n\to\iy}L_\nu (\sqrt{t_n}u_n, \sqrt{s_n} v_n)\ge m_\nu,
\end{align*}
}%
a contradiction with Lemma \ref{lemma3-1}. So Case 1 is impossible.

{\bf  Case 2.} Either $u\equiv 0,\,v\not\equiv 0$ or $u\not\equiv 0, v\equiv 0$.

Without loss of generality, we assume that $u\not\equiv 0, v\equiv 0$. Note that $J_\nu'(u, v)(u, 0)=0$ yields
$$\|u\|_{\la_1}^2=|u|_4^4\le S(\la_1)^{-2}\|u\|_{\la_1}^4,$$
which implies $\|u\|_{\la_1}^2\ge S(\la_1)^2$.

{\bf Case 2.1} Up to a subsequence, $\lim_{n\to\iy}\|u_n-u\|>0$.

Denote $w_n=u_n-u$. Note that $(u_n, v_n)\in \N_\nu$. Then by Brezis-Lieb Lemma (\cite{BL}) and Lemma \ref{lemma3} we conclude that
{\allowdisplaybreaks
\begin{align*}
&\|w_n\|_{\la_1}^2=\int_{\R^4}w_n^4+2\nu\int_{\R^4}w_n^2 v_n^2+ o(1),\\
&\|v_n\|_{\la_2}^2=\int_{\R^4}v_n^4+2\nu\int_{\R^4}w_n^2 v_n^2+ o(1),
\end{align*}
}%
Similarly as above, it is easy to prove that there exist $t_n, s_n>0$ such that $(\sqrt{t_n}w_n, \sqrt{s_n} v_n)\in\N_\nu$ and $(t_n, s_n)\to (1, 1)$ as $n\to\iy$. Hence,
{\allowdisplaybreaks
\begin{align*}
c_\nu &=\lim_{n\to\iy}J_\nu (u_n, v_n)=\lim_{n\to\iy}\frac{1}{4}(\|u_n\|_{\la_1}^2+\|v_n\|_{\la_2}^2)\\
&=\frac{1}{4}\|u\|_{\la_1}^2+\lim_{n\to\iy}\frac{1}{4}\left(t_n\|w_n\|_{\la_1}^2+s_n\|v_n\|_{\la_2}^2\right)\\
&>\lim_{n\to\iy}J_\nu (\sqrt{t_n}w_n, \sqrt{s_n} v_n)\ge c_\nu,
\end{align*}
}%
a contradiction. So Case 2.1 is impossible.

{\bf Case 2.2} $u_n\to u$ strongly in $D^{1, 2}(\RN)$.

Then $u_n^2\to u^2$ strongly in $L^2(\R^4)$. Recall that $v_n\rightharpoonup 0$ in $D^{1, 2}(\R^4)$, so
$v_n^2\rightharpoonup 0$ weakly in $L^2(\R^4)$, which easily implies
$$\int_{\R^4}u_n^2 v_n^2\le\int_{\R^4}u^2v_n^2+\int_{\R^4}|u_n^2-u^2|v_n^2=o(1).$$
Then we have
$$\|v_n\|_{\la_2}^2=|v_n|_4^4+o(1)\le S(\la_2)^{-2}\|v_n\|_{\la_2}^4+o(1),$$
Since Lemma \ref{lemma3-2} yields $\lim_{n\to \iy}\|v_n\|_{\la_2}^2>0$, so $\lim_{n\to \iy}\|v_n\|_{\la_2}^2\ge S(\la_2)^2$, and
{\allowdisplaybreaks
\begin{align*}
c_\nu &=\lim_{n\to\iy}J_\nu (u_n, v_n)=\lim_{n\to\iy}\frac{1}{4}(\|u_n\|_{\la_1}^2+\|v_n\|_{\la_2}^2)\\
&=\frac{1}{4}\|u\|_{\la_1}^2+\lim_{n\to\iy}\frac{1}{4}\|v_n\|_{\la_2}^2
\ge \frac{1}{4}\left(S(\la_1)^2+S(\la_2)^2\right),
\end{align*}
}%
a contradiction with Lemma \ref{lemma3-1}. So Case 2.2 is impossible, and so Case 2 is impossible.

Since neither Case 1 nor Case 2 is true, we obtain that $u\not\equiv 0$ and $v\not\equiv 0$. Since $J'_\nu(u, v)=0$, so $(u, v)\in\N_\nu$.
Meanwhile, Fatou Lemma implies that
$$c_\nu\le J_\nu(u, v)\le \liminf_{n\to\iy} J_\nu(u_n, v_n)=c_\nu,$$
so $J_\nu(u, v)=c_\nu$. Then $(|u|, |v|)\in \N_\nu$ and $J_\nu(|u|, |v|)=c_\nu$. Since $2\nu<1$ and $\al=\bb=2$, then
by repeating the proof of Lemma \ref{lemma1}, we can prove that $J_\nu'(|u|, |v|)=0$.
By the maximum principle, $|u|, |v|>0$ in $\R^4\setminus \{0\}$.
Therefore, $(|u|, |v|)$ is a positive ground state solution of (\ref{eq4-1}).

To finish the proof, it suffices to prove $c_\nu\to \frac{1}{4}\left(S(\la_1)^2+S(\la_2)^2\right)$ as $\nu\to 0$.
From the above argument, we may assume that $(u_\nu, v_\nu)$ is a positive ground state solution of (\ref{eq4-1}) with $c_\nu=J_\nu(u_\nu, v_\nu)$ for any $\nu\in (0, \nu_1)$. Since $c_\nu< \frac{1}{4}\left(S(\la_1)^2+S(\la_2)^2\right)$, we see that $(u_\nu, v_\nu)$ are uniformly bounded in $\BD$.
Then
$$\|u_\nu\|_{\la_1}^2=|u_\nu|_4^4+2\nu\int_{\R^4}u_\nu^2 v_\nu^2\le S(\la_1)^{-2}\|u\|_{\la_1}^4 + O(\nu).$$
From (\ref{eq4-14}) we see that $\liminf_{\nu\to 0}\|u_\nu\|_{\la_1}^2>0$, so $\liminf_{\nu\to 0}\|u_\nu\|_{\la_1}^2\ge S(\la_1)^2$. Similarly,
we can prove that $\liminf_{\nu\to 0}\|v_\nu\|_{\la_2}^2\ge S(\la_2)^2$, and so
$\liminf_{\nu\to 0}c_\nu\ge \frac{1}{4}\left(S(\la_1)^2+S(\la_2)^2\right).$ That is,
$\lim_{\nu\to 0}c_\nu=\frac{1}{4}\left(S(\la_1)^2+S(\la_2)^2\right).$ This completes the proof.\hfill$\square$

\vskip0.1in

\s{Proof of Theorem \ref{th1-4}: A variational perturbation approach}
\renewcommand{\theequation}{4.\arabic{equation}}

In this section, we give the proof of Theorem \ref{th1-4}, and this result will be used in the proof of Theorem \ref{theorem2}.
Assume that $N\ge 3$, $\la_1, \la_2\in (0, \La_N)$ and  (\ref{eq1-8}) hold. Let $\nu>0$. To obtain positive solutions of
(\ref{eq0-0}), we consider the following modified problem
\be\label{s7eq2-0}\begin{cases}-\Delta u -\frac{\la_1}{|x|^2}u-u_+^{2^\ast-1} =
\nu \al u_+^{\al-1}v_+^\bb , \quad x\in \RN,\\
  -\Delta v -\frac{\la_2}{|x|^2}v-v_+^{2^\ast-1} =
\nu \bb u_+^{\al}v_+^{\bb-1},     \quad x\in \RN,\\
u(x),\,\, v(x)\in D^{1, 2}(\RN),\end{cases}\ee
where $u_\pm(x):=\max\{\pm u(x), 0\}$ and so does $v_\pm$.
The associated energy functional of (\ref{s7eq2-0}) is
\be\label{s7eq2-1}\oJ_\nu (u, v):=\frac{1}{2}\|u\|_{\la_1}^2+\frac{1}{2}\|v\|_{\la_2}^2-\frac{1}{2^\ast}\int_{\RN}\left(u_+^{2^\ast}
+v_+^{2^\ast}\right)-\nu\int_{\RN}{u_+^\al v_+^\bb},\ee
Then it is standard to prove that $\oJ\in C^1(\BD, \R)$.
Define
\begin{align*}& C^\iy_{0, r}(\RN):=\{u\in C_0^\iy(\RN) : \quad u \,\,\hbox{is radially symmetric}\},\\
& D^{1, 2}_{r}(\RN):=\{u\in D^{1,2}(\RN) : \quad u \,\,\hbox{is radially symmetric}\},\end{align*}
and $\BD_r:=D^{1, 2}_r (\RN)\times D^{1, 2}_r (\RN)$. Then $\BD_r$ is a subspace of $\BD$ with norm $\|\cdot\|_\BD$.
In this section, we consider the functional $\oJ_\nu$ restricted to $\BD_r$. By Palais's Symmetric Criticality Principle,
any critical points of $\oJ_\nu: \BD_r\to \R$ are radially symmetric
solutions of (\ref{s7eq2-0}).

Without loss of generality, we assume that
$S(\la_1)\le S(\la_2)$. Then we see from (\ref{eq1-6}) that
\be\label{s7eq2-2}S(\la_1)^{N/4}=\|z_\mu^1\|_{\la_1}\le \|z_\mu^2\|_{\la_2}=S(\la_2)^{N/4},\quad \forall\,\mu>0.\ee
Define
$$P_{\la_i}(u):=\begin{cases}\frac{|u|^{2^\ast}_{2^\ast}}{\|u\|_{\la_i}^2}, & \hbox{if}\quad u\in D^{1, 2}(\RN)\setminus\{0\},\\
0, & \hbox{if}\quad u=0, \end{cases}\quad i=1, 2.$$
By (\ref{eq1-5-1}) it is easy to prove that $P_{\la_i}\in C(D^{1, 2}(\RN),\, \R)$. Note that $P_{\la_i}(u)=1$ is equivalent to $I_{\la_i}'(u)u=0$.
Then by (\ref{eq1-5})-(\ref{eq1-7}) it is easy to check that
\be\label{s7eq2-3}M_i:=\frac{1}{N}S(\la_i)^{N/2}=\inf\limits_{\stackrel{u\in D^{1, 2}(\R^N)}{ P_{\la_i}(u)=1}}I_{\la_i}(u),\quad
i=1, 2.\ee
By (\ref{eq1-6})-(\ref{eq1-7}) we have
{\allowdisplaybreaks
\begin{align}\label{s7eq2-4}I_{\la_i}(t z_1^i)&=\frac{t^2}{2}\|z_1^i\|^2_{\la_i}-\frac{t^{2^\ast}}{2^\ast}|z_1^i|^{2^\ast}_{2^\ast}
=\left(\frac{t^2}{2}-\frac{t^{2^\ast}}{2^\ast}\right)\|z_1^i\|^2_{\la_i},\,\,i=1, 2.\end{align}
}%
Note that
\begin{align}\label{s7eq2-5}I_{\la_i}(z_1^i)=\max_{t>0}I_{\la_i}(t z_1^i)= M_i,\,\,i=1, 2,\quad M_1\le M_2,\end{align}
it is easily seen that
there exist $0<t_0<1<t_1$ such that
\be\label{s7eq2-6} I_{\la_i}(t z_1^i)\le M_1/4 \quad\hbox{for}\,\,\, t\in (0, t_0]\cup[t_1, \iy),\,\, i=1, 2.\ee
Define
{\allowdisplaybreaks
\begin{gather*}
    \widetilde{\gamma}_i(t):=t z_1^i\,\,\,\hbox{for $0 \le t\le t_1$},\,\, i=1, 2;\quad
    \widetilde{\gamma}(t, s):=(\widetilde{\gamma}_1(t), \widetilde{\ga}_2(s)).
\end{gather*}
}%
Then $\widetilde{\gamma}(t, s)\in \BD_r$ for all $(t, s)$ and there exists a constant $\mathcal{C}>0$ such that
\be\label{s7eq2-7}\max_{(t, s)\in [0, t_1]\times [0, t_1]}\|\widetilde{\gamma}(t, s)\|_\BD\le \mathcal{C}.\ee
Denote $Q:=[0, t_1]\times[0, t_1]$ for convenience. For $\nu\ge 0$, we define
$$a_{\nu}:=\inf\limits_{\ga\in \G}\max\limits_{(t,s)\in Q} \oJ_\nu(\ga(t, s)), \quad d_\nu:=\max\limits_{(t,s)\in Q} \oJ_\nu(\widetilde{\ga}(t, s)),$$
where
{\allowdisplaybreaks
\begin{align}\label{s7eq2-8}\G:=\Big\{& \ga\in C(Q, \BD_r)\,\,\,:\,\,\, \max_{(t, s)\in Q}\|\ga(t, s)\|_\BD\le 2S(\la_2)^{N/4}+\mathcal{C},\nonumber\\
&\ga(t, s)=\widetilde{\ga}(t, s)\,\,\,\hbox{for}\,\,\,(t, s)\in Q\backslash (t_0, t_1)\times (t_0, t_1)\Big\}.\end{align}
}%

The definition of $a_{\nu}$ is different
from the definitions of usual mountain-pass values (cf.
\cite{AR-pass}). All paths in $\G$ are required to be
uniformly bounded in $\BD$ by $2S(\la_2)^{N/4}+\mathcal{C} $, which will play a crucial role in the proof of Lemma \ref{s7lemma2} below.

\bl\label{s7lemma2} There hold $d_\nu<d_0$ for $\nu>0$ and $\lim\limits_{\nu\to 0+}a_\nu=\lim\limits_{\nu\to 0+}d_\nu= a_0=d_0=M_1+M_2$.\el

\noindent{\bf Proof. } Note that $z_1^i>0$ in $\RN\setminus\{0\}$, we have
\begin{align*}
d_0=\max\limits_{(t,s)\in Q} \oJ_0(\widetilde{\ga}(t, s))=\max_{t\in (0, t_1)}I_{\la_1}(t z_1^1)+\max_{s\in (0, t_1)}I_{\la_2}(s z_1^2).
\end{align*}
Then by (\ref{s7eq2-4})-(\ref{s7eq2-5}) we see that
\be\label{s7eq2-9}d_0=\oJ_0(\widetilde{\ga}(1, 1))=M_1+M_2, \quad
\oJ_0(\widetilde{\ga}(t, s))<d_0\,\,\,\hbox{for $(t, s)\in Q\setminus\{ (1, 1)\}$}.\ee

Fix any $\nu>0$.  Note that there exists $(t_\nu, s_\nu)\in Q\setminus\{(0, 0)\}$
such that $d_\nu= \oJ_\nu(\widetilde{\ga}(t_\nu, s_\nu))$. If $(t_\nu, s_\nu)=(1, 1)$, then we see from (\ref{s7eq2-9}) that
\begin{align*}
d_\nu=\oJ_\nu(\widetilde{\ga}(1, 1))=\oJ_0(\widetilde{\ga}(1, 1))-\nu\intR{(z_1^1)^\al (z_1^2)^\bb}<d_0;
\end{align*}
if $(t_\nu, s_\nu)\neq(1, 1)$, then we deduce from (\ref{s7eq2-9}) again that
\begin{align*}
d_\nu=\oJ_\nu(\widetilde{\ga}(t_\nu, s_\nu))=\oJ_0(\widetilde{\ga}(t_\nu, s_\nu))-\nu t_\nu^\al s_\nu^\bb\intR{(z_1^1)^\al (z_1^2)^\bb}<d_0,
\end{align*}
that is, $d_\nu<d_0$ for any $\nu>0$. Note that $\widetilde{\ga}\in\G$, we have $a_\nu\le d_\nu$, that is
\be\label{s7eq2-10}\liminf_{\nu\to0+}a_\nu\le\liminf_{\nu\to 0+}d_\nu,\,\,\limsup_{\nu\to0+}a_\nu\le\limsup_{\nu\to 0+}d_\nu\le d_0\quad\text{and}\,\,\, a_0\le d_0.\ee
On the other hand, for any $\ga (t, s)=(\ga_1(t, s), \ga_2(t, s))\in \G$, we define $\Upsilon(\ga): [t_0, t_1]^2\to \R^2$ by
$$\Upsilon(\ga) (t, s):=\big(P_{\la_1}(\ga_1(t, s))-1,\,\,\,P_{\la_2}(\ga_2(t, s))-1\big).$$
By the definitions of $P_{\la_i}$ and $\widetilde{\ga}$, it is easily seen that
$$\Upsilon(\widetilde{\ga}) (t, s)=\left(t^{2^\ast-2}-1, \,s^{2^\ast-2}-1\right).$$
Then $\deg(\Upsilon(\widetilde{\ga}), [t_0, t_1]^2, (0, 0))=1$. By (\ref{s7eq2-8}) we see that
for any $(t, s)\in \partial ([t_0, t_1]^2)$,  $\Upsilon(\ga) (t, s)=\Upsilon(\widetilde{\ga}) (t, s)\neq(0, 0)$. Therefore,
$\deg(\Upsilon(\ga), [t_0, t_1]^2, (0, 0))$ is well defined and
$$\deg(\Upsilon(\ga), [t_0, t_1]^2, (0, 0))=\deg(\Upsilon(\widetilde{\ga}), [t_0, t_1]^2, (0, 0))=1.$$
Then there exists $(t_2, s_2)\in [t_0, t_1]^2$ such that $\Upsilon(\ga) (t_2, s_2)=(0, 0)$, that
is $P_{\la_1}(\ga_1(t_2, s_2))=1$ and $P_{\la_2}(\ga_2(t_2, s_2))=1$.
Combining these with (\ref{s7eq2-3}), we have
{\allowdisplaybreaks
\begin{align}\label{s7eq2-11}
\max_{(t, s)\in Q}\oJ_0(\ga(t, s))&\ge \oJ_0(\ga(t_2, s_2))\ge J_0(\ga(t_2, s_2))\nonumber\\
&=I_{\la_1}(\ga_1(t_2, s_2))+I_{\la_2}(\ga_2(t_2, s_2))\nonumber\\
&\ge M_1+M_2=d_0.
\end{align}
}%
Therefore, $a_0\ge d_0$. By (\ref{s7eq2-10}) one get that
$a_0=d_0$.

Assume by contradiction that $\liminf_{\nu\to 0+} a_\nu <d_0$. Then there exists $\e>0$, $\nu_n\to 0+$ and $\ga_n=(\ga_{n, 1}, \ga_{n, 2})\in\G$ such that
$$\max_{(t, s)\in Q}J_{\nu_n}(\ga_n(t, s))\le d_0-2\e.$$
Recall $\al+\bb=2^\ast$. By (\ref{s7eq2-8}) and H\"{o}lder inequality, there exists $n_0$ large enough such that
$$\max_{(t, s)\in Q}\nu_n\left|\intR{\left(\ga_{n,1}(t, s)_+\right)^\al\left(\ga_{n, 2}(t, s)_+\right)^\bb}\right|\le C\nu_n\le \e,\quad \forall\,\, n\ge n_0,$$
and so
$$a_0\le\max_{(t, s)\in Q}\oJ_{0}(\ga_n(t, s))\le \max_{(t, s)\in Q}\oJ_{\nu_n}(\ga_n(t, s))+\e\le d_0-\e,\quad \forall\,\, n\ge n_0,$$
a contradiction with $a_0=d_0$. Therefore, $\liminf_{\nu\to 0+} a_\nu \ge d_0$. Combining this with (\ref{s7eq2-10}), we complete the proof.
\hfill$\square$\\

Recall (\ref{eq1-3}), we define $X:=Z_1\times Z_2\subset\BD_r$, and
\begin{gather}
    X^\dd:=\{(u, v)\in \BD_r : \hbox{dist}((u, v), X)\le \dd\},\,\,
    \oJ_\nu^{d}:=\{(u, v)\in \BD_r : \oJ_\nu(u, v)\le d\};\nonumber\\
    \label{s7eq2-12}\dd:=\min\left\{\frac{1}{2}, \,\,\frac{1}{16}S(\la_1)^{N/4},\,\,\frac{\mathcal{C}}{4}\right\}.
\end{gather}
Here dist$((u, v), X):=\inf\{\|(u-\vp,\, v-\psi)\|_{\BD} \,:\, (\vp, \psi)\in X\}$.  Define
$$\left\|\oJ_\nu'(u, v)\right\|:=\sup\left\{\oJ_\nu'(u, v)(\psi, \phi) \,:\, (\psi, \phi)\in\BD_r,\,\, \|(\psi, \phi)\|_{\BD}=1\right\}.$$

\bl\label{s7lemma3} Recall $\dd$ in (\ref{s7eq2-12}). Then there
exist $0<\sg<1$ and $\nu_3\in (0, 1)$, such that $\|\oJ_\nu'(u, v)\|\ge \sg$ holds for any
$(u, v)\in \oJ_\nu^{d_\nu}\cap (X^{\dd}\backslash X^{\dd/2})$ and $\nu\in (0, \nu_3]$.  \el

\noindent{\bf Proof. }
Assume by contradiction, there exist $\nu_n\to 0+$ and $(u_n, v_n)\in \oJ_{\nu_n}^{d_{\nu_n}}\cap (X^{\dd}\backslash X^{\dd/2})$
such that $\|\oJ_{\nu_n}'(u_n, v_n)\|\to 0$. Then there exist $\mu_{i, n}>0, i=1, 2, n\in\mathbb{N}$ such that
$$\left\|(u_n, v_n)-\left(z_{\mu_{1, n}}^1, z_{\mu_{2, n}}^2\right)\right\|_{\BD}\le 2\dd,\quad\forall\,n\in\mathbb{N}.$$
Hence $(u_n, v_n)$ are uniformly bounded in $\BD_r$, and up to a subsequence,
we may assume that $(u_n, v_n)\rightharpoonup (\ou, \ov)$ weakly in $\BD_r$. However, since $(z^1_\mu, z^2_\mu)\rightharpoonup (0, 0)$ weakly in $\BD_r$ as $\mu\to\iy$, we know that $X$ is not compact in $\BD_r$. So it seems very difficult for us to show that $\ou\not\equiv 0$ and $\ov\not\equiv 0$.
To overcome this difficulty, let us define
$$\widetilde{u}_n(x):=\mu_{1, n}^{\frac{N-2}{2}}u_n(\mu_{1, n}x),\quad \widetilde{v}_n(x):=\mu_{2, n}^{\frac{N-2}{2}}v_n(\mu_{2, n}x).$$
Note that $\|\cdot\|_{\la_i}, i=1, 2$ are invariant with respect to the
transformation $u(\cdot)\mapsto \mu^{-\frac{N-2}{2}}u(\frac{\cdot}{\mu})$ for all $\mu>0$. Therefore,
\begin{align*}
\left\|\widetilde{u}_n-z_1^1\right\|_{\la_1}=\left\|u_n-z_{\mu_{1, n}}^1\right\|_{\la_1}\le 2\dd,
\quad \left\|\widetilde{v}_n-z_1^2\right\|_{\la_2}=\left\|v_n-z_{\mu_{2, n}}^2\right\|_{\la_2}\le 2\dd.
\end{align*}
This means that $(\widetilde{u}_n, \widetilde{v}_n)$ are uniformly bounded in $\BD_r$. Up to a subsequence,
we may assume that $(\wu_n, \wv_n)\rightharpoonup (\wu, \wv)$ weakly in $\BD_r\cap L^{2^\ast}(\RN)\times L^{2^\ast}(\RN)$.
Then we have $\|\wu-z_1^1\|_{\la_1}\le\liminf_{n\to\iy}\|\widetilde{u}_n-z_1^1\|_{\la_1}\le 2\dd$. Combining this with (\ref{s7eq2-2}) and (\ref{s7eq2-12}),
we get that $\wu\not\equiv 0$. Similarly, $\wv\not\equiv 0$.

Take any $\wphi\in C_{0, r}^\iy(\RN)$ such that $\|\wphi\|_{\la_1}=1$, we define
\be\label{test}\phi_n(x):=\mu_{1, n}^{-\frac{N-2}{2}}\wphi\left(\frac{x}{\mu_{1, n}}\right).\ee
Then $\|\phi_n\|_{\la_1}=\|\wphi\|_{\la_1}=1$. Since $\nu_n\to 0$, by H\"{o}lder inequality and Sobolev inequality, we easily obtain that
$$\lim_{n\to\iy}\left|\nu_n \al\intR{(u_n)_+^{\al-1}(v_n)_+^\bb\phi_n}\right|=0.$$
Therefore,
{\allowdisplaybreaks
\begin{align*}
0&=\lim_{n\to\iy}\oJ_{\nu_n}'(u_n, v_n)(\phi_n, 0)\\
&=\lim_{n\to\iy}\intR{\nabla u_n\nabla \phi_n-\frac{\la_1}{|x|^2}u_n\phi_n-(u_n)_+^{2^\ast-1}\phi_n}\\
&=\lim_{n\to\iy}\intR{\nabla \wu_n\nabla \wphi-\frac{\la_1}{|x|^2}\wu_n\wphi-(\wu_n)_+^{2^\ast-1}\wphi}\\
&=\intR{\nabla \wu\nabla \wphi-\frac{\la_1}{|x|^2}\wu\wphi-\wu_+^{2^\ast-1}\wphi},\,\,\,\hbox{holds for any}\,\,\wphi\in C_{0, r}^\iy(\RN),
\end{align*}
}%
that is, $-\Delta \wu-\frac{\la_1}{|x|^2}\wu=\wu_+^{2^\ast-1}$ and $\wu\in D_r^{1, 2}(\RN)$. By testing this equation with $\wu_-$ we see that $\wu\ge 0$. By the maximum principle, one has that $\wu>0$ in $\RN\setminus\{0\}$, that is $\wu$ is a positive solution of (\ref{eq1-2}) with $i=1$.
Then by (\ref{eq1-2})-(\ref{eq1-3}) we get that $\wu\in Z_1$. Similarly,
we may prove that $\wv\in Z_2$, that is, $(\wu, \wv)\in X$.
Recall that $\oJ_{\nu_n}(u_n, v_n)\le d_{\nu_n}$ and $\al+\bb=2^\ast$, we deduce from Lemma \ref{s7lemma2} that
{\allowdisplaybreaks
\begin{align*}
M_1+M_2 &\ge\lim_{n\to\iy}\left(\oJ_{\nu_n}(u_n, v_n)-\frac{1}{2^\ast}\oJ_{\nu_n}'(u_n, v_n)(u_n, v_n)\right)\\
&=\frac{1}{N}\lim_{n\to\iy}\left(\|u_n\|^2_{\la_1}+\|v_n\|_{\la_2}^2\right)=\frac{1}{N}\lim_{n\to\iy}\left(\|\wu_n\|^2_{\la_1}+\|\wv_n\|_{\la_2}^2\right)\\
&\ge\frac{1}{N}\|\wu\|_{\la_1}^2+\frac{1}{N}\|\wv\|_{\la_2}^2=M_1+M_2,
\end{align*}
}%
which implies that all inequalities above are identities, so $\oJ_{\nu_n}(u_n, v_n)\to M_1+M_2$ and  $(\wu_n, \wv_n)\to (\wu, \wv)\in X$ strongly in $\BD_r$. Then
$$\|(\wu_n, \wv_n)-(\wu, \wv)\|_{\BD}\le \dd/4, \quad\hbox{for $n$ large enough,}$$
and so
$$\|(u_n, v_n)-(\ou_n, \ov_n)\|_{\BD}\le \dd/4, \quad\hbox{for $n$ large enough,}$$
where
$$(\ou_n(x), \ov_n(x)):=\left(\mu_{1, n}^{-\frac{N-2}{2}}\wu\left(\frac{x}{\mu_{1, n}}\right),\, \mu_{2, n}^{-\frac{N-2}{2}}\wv\left(\frac{x}{\mu_{2, n}}\right)\right)\in X.$$
This contradicts with $(u_n, v_n)\not\in X^{\dd/2}$ for any $n$. This completes the
proof.\hfill$\square$

\bl\label{s7lemma4} There exists $\nu_4\in (0, \nu_3]$ and $\e>0$ such that for any $\nu\in (0, \nu_4]$,
$$\oJ_\nu(\widetilde{\ga}(t, s))\ge a_\nu-\e\quad \hbox{implies that}\quad \widetilde{\ga}(t, s)\in X^{\dd/2}.$$\el

\noindent{\bf Proof. } Assume by contradiction that there exist $\nu_n\to 0$, $\e_n\to 0$ and $(t_n, s_n)\in Q$ such that
\be\label{s7eq2-14}\oJ_{\nu_n}(\widetilde{\ga}(t_n, s_n))\ge a_{\nu_n}-\e_n\quad
\hbox{and}\quad \widetilde{\ga}(t_n, s_n)\not\in X^{\dd/2},\,\,\,\forall \,n\in\mathbb{N}.\ee
Passing to a subsequence, we may assume that $(t_n, s_n)\to (\tilde{t}, \tilde{s})\in Q$.
Then by Lemma \ref{s7lemma2} and letting $n\to\iy$ in (\ref{s7eq2-14}), we have
$$J_0(\widetilde{\ga}(\tilde{t}, \tilde{s}))\ge \lim_{n\to\iy}a_{\nu_n}=M_1+M_2.$$
Combining this with (\ref{s7eq2-9}), we obtain that $(\tilde{t}, \tilde{s})=(1, 1)$. Hence,
$$\lim_{n\to\iy}\|\widetilde{\ga}(t_n, s_n)-\widetilde{\ga}(1, 1)\|_{\BD}= 0.$$
However, $\widetilde{\ga}(1, 1)=(z_1^1, z_1^2)\in X$, which is a contradiction with (\ref{s7eq2-14}).\hfill$\square$\\

Let \be\label{s7eq2-15}\e_0 :=\min\left\{\frac{\e}{2}, \frac{M_1}{4}, \frac{1}{8}\dd
\sg^2\right\},\ee where $\dd, \sg$ are seen in Lemma \ref{s7lemma3}. By Lemma \ref{s7lemma2} there exists $\nu_2\in (0, \nu_4]$ such that
\be\label{s7eq2-16}|a_\nu-d_\nu|<\e_0, \quad|a_\nu- (M_1+M_2)|<\e_0,\quad\forall \,\,\nu\in (0, \nu_2].\ee

\bl\label{s7lemma5} For any fixed $\nu\in (0, \nu_2]$, there exists $\{(u_n, v_n)\}_{n=1}^\iy\subset
X^\dd \cap \oJ_{\nu}^{d_\nu}$ such that
$$\left\|\oJ_\nu' (u_n, v_n)\right\|\to 0 \quad\hbox{as}\,\, n\to \iy.$$\el

\noindent{\bf Proof. }Fix any $\nu\in (0, \nu_2]$. Assume by contradiction that there exists
$0<l(\nu)<1$ such that $\|\oJ_\nu' (u, v)\| \ge l(\nu)$
on $X^\dd\cap \oJ_{\nu}^{d_\nu}$. Then
there exists a pseudo-gradient vector field $T_\nu$ in $\BD_r$ which is defined on
a neighborhood $Z_\nu\subset\BD_r$ of $X^\dd\cap
\oJ_{\nu}^{d_{\nu}}$ (cf. \cite{S2}), such that for any $(u,
v)\in Z_\nu$, there holds
{\allowdisplaybreaks
\begin{gather*}
\|T_\nu(u, v)\|_\BD \le 2\min\{1, \|\oJ_\nu'(u, v)\|\},\\
\oJ_{\nu}'(u, v) \big(T_\nu(u, v)\big) \ge \min\{1,
\|\oJ_{\nu}'(u, v)\|\}\|\oJ_{\nu}'(u, v)\|.
\end{gather*}
}%
Let $\eta_\nu$ be a Lipschiz continuous function on $\BD_r$ such
that $0\le \eta_\nu\le 1, \, \eta_\nu\equiv1$ on $X^\dd \cap \oJ_{\nu}^{d_{\nu}}$ and $ \eta_\nu\equiv 0$ on
$\BD_r\backslash Z_\nu$. Let $\xi_\nu$ be a Lipschiz continuous
function on $\R$ such that $0\le \xi_\nu\le 1, \, \xi_\nu(l)\equiv1$
if $|l-a_\nu|\le \frac{\e}{2}$ and $\xi_\nu(l)\equiv 0$ if
$|l-a_{\nu}|\ge \e$. Let
\begin{displaymath}
e_\nu(u, v)=\begin{cases} -\eta_\nu (u, v)\xi_\nu (\oJ_{\nu}(u, v))T_\nu (u, v) & \hbox{if}\quad (u, v)\in Z_\nu, \\
0 & \hbox{if}\quad (u, v)\in \BD_r\backslash Z_\nu, \end {cases}
\end{displaymath}
then there exists a global solution $\psi_\nu: \BD_r\times [0,
+\iy)\to \BD_r$ to the following initial value problem
\begin{displaymath}
\begin{cases}\frac{d}{d\theta}\psi_\nu(u, v, \theta)= e_\nu (\psi_\nu(u, v, \theta)), \\
\psi_\nu(u, v, 0)=(u, v). \end {cases}
\end{displaymath}
It is easy to see from Lemma \ref{s7lemma3} and (\ref{s7eq2-15})-(\ref{s7eq2-16})  that $\psi_\nu$ has following properties:\\
(1) $ \displaystyle \psi_\nu(u, v, \theta)=(u, v)$ if $\theta=0$ or $(u, v)\in \BD_r\backslash Z_\nu$ or $|\oJ_{\nu}(u, v)- a_{\nu}|\ge \e$;\\
(2) $\displaystyle\left\|\frac{d}{d\theta}\psi_\nu(u, v, \theta)\right\|_\BD \le 2$;\\
(3) $ \displaystyle\frac{d}{d\theta}\oJ_{\nu}(\psi_\nu(u, v, \theta))= \oJ_{\nu}'(\psi_\nu(u, v, \theta))\big(e_\nu(\psi_\nu(u, v, \theta))\big)\le 0$;\\
(4) $ \displaystyle\frac{d}{d\theta}J_{\nu}(\psi_\nu(u, v, \theta))\le -l(\nu)^2$ if $ \displaystyle\psi_\nu(u, v, \theta)\in X^{\dd}\cap \left(\oJ_{\nu}^{d_\nu}\backslash \oJ_{\nu}^{a_\nu-\e/2}\right)$;\\
(5) $\displaystyle\frac{d}{d\theta}J_{\nu}(\psi_\nu(u, v, \theta))\le -\sg^2$ if $\displaystyle\psi_\nu(u, v, \theta)\in (X^{\dd}\backslash X^{\dd/2})\cap \left(\oJ_{\nu}^{d_\nu}\backslash \oJ_{\nu}^{a_\nu-\e/2}\right)$.\\

{\noindent }{\bf Step 1.} For any $(t, s)\in Q$, we claim that there exists $\theta_{t, s} \in
[0, +\iy)$ such that $\psi_\nu(\widetilde{\ga}(t, s), \theta_{t, s})\in
\oJ_{\nu}^{a_{\nu}-\e_0}$, where $\e_0$ is seen in
(\ref{s7eq2-15}).

Assume by contradiction that there exists $(t, s)\in Q$ such that
$$\oJ_{\nu}(\psi_\nu(\widetilde{\ga}(t, s), \q))>a_{\nu}-\e_0,\qu\forall\, \q\ge 0.$$
Note that $\e_0<\e$, we see from Lemma
\ref{s7lemma4} that $\widetilde{\ga}(t, s)\in X^{\dd/2}$. Note that
$\oJ_{\nu}(\widetilde{\ga}(t, s))\le d_{\nu}< a_{\nu}+\e_0$, we
see from the property (3) that
$$a_{\nu}-\e_0<\oJ_{\nu}(\psi_\nu(\widetilde{\ga}(t, s), \q))\le d_{\nu}< a_{\nu}+\e_0, \qu\forall\, \q\ge 0.$$
This implies $\xi_\nu(\oJ_{\nu}(\psi_\nu(\widetilde{\ga}(t, s), \q)))\equiv
1$. If $\psi_\nu(\widetilde{\ga}(t, s), \q)\in X^\dd$ for all $\q\ge
0$, then $\eta_\nu(\psi_\nu(\widetilde{\ga}(t, s), \q))\equiv 1$, and
$\|\oJ_{\nu}'(\psi_\nu(\widetilde{\ga}(t, s), \q))\|\ge l(\nu)$ for all
$\q>0$. Then we see from property (4) that
$$\oJ_{\nu}\left(\psi_\nu\left(\widetilde{\ga}(t, s),\frac{\e}{l(\nu)^2}\right)\right)
\le
a_{\nu}+\frac{\e}{2}-\int_{0}^{\frac{\e}{l(\nu)^2}}l(\nu)^2\,dt\le
a_{\nu}-\frac{\e}{2},$$ a contradiction. Thus, there exists
$\q_{t, s}>0$ such that $\psi_\nu(\widetilde{\ga}(t, s), \q_{t,s})\not\in
X^\dd$. Note that $\widetilde{\ga}(t, s)\in X^{\dd/2}$, there
exists $0<\q_{t, s}^1<\q_{t, s}^2\le \q_{t, s}$ such that $\psi_\nu(\widetilde{\ga}(t, s),
\q_{t, s}^1)\in \partial X^{\dd/2}$, $\psi_\nu(\widetilde{\ga}(t, s),
\q_{t, s}^2)\in \partial X^{\dd}$ and $\psi_\nu(\widetilde{\ga}(t, s),
\q)\in X^\dd\backslash X^{\dd/2}$ for all $\q\in (\q_{t, s}^1,
\q_{t, s}^2)$. Then by Lemma \ref{s7lemma3} we have
$\|\oJ_{\nu}'(\psi_\nu(\widetilde{\ga}(t, s), \q))\|\ge \sg$ for all $\q\in
(\q_{t, s}^1, \q_{t, s}^2)$. Then using the property (2) we have
\begin{align*}
\dd/2\le \left\|\psi_\nu(\widetilde{\ga}(t, s), \q_{t, s}^2)-\psi_\nu(\widetilde{\ga}(t, s),
\q_{t, s}^1)\right\|_\BD\le 2 |\q_{t, s}^2 - \q_{t, s}^1|,
\end{align*}
that is, $\q_{t, s}^2-\q_{t, s}^1\ge \dd/4$. This implies from (\ref{s7eq2-15}) and property (5) that
{\allowdisplaybreaks
\begin{align*}
\oJ_{\nu}\left(\psi_\nu(\widetilde{\ga}(t, s), \q_{t, s}^2)\right)&\le
\oJ_{\nu}(\psi_\nu(\widetilde{\ga}(t, s),
\q_{t, s}^1))+\int_{\q_{t,s}^1}^{\q_{t, s}^2}\frac{d}{d\q}\oJ_{\nu}(\psi_\nu(u, v,
\q))\,d\q\\
&< a_{\nu} +\e_0 -\sg^2 (\q_{t, s}^2-\q_{t, s}^1)\le  a_{\nu} +\e_0
-\frac{1}{4}
\dd \sg^2\\
&\le a_{\nu}-\e_0,
\end{align*}
}%
which is a contradiction.

By Step 1 we can define $T(t, s):=\inf\{\q\ge 0 :
\oJ_{\nu}(\psi_\nu(\widetilde{\ga}(t, s), \q))\le a_{\nu}-\e_0\}$ and
let $\ga(t, s):=\psi_\nu(\widetilde{\ga}(t, s), T(t, s))$, then
$\oJ_{\nu}(\ga(t, s))\le a_{\nu}-\e_0$ for all $(t, s)\in Q$.\\

{\noindent }{\bf Step 2.} We shall prove that $\ga(t, s) \in \G$.

For any $(t, s)\in Q\backslash (t_0, t_1)\times (t_0, t_1)$, by (\ref{s7eq2-5})-(\ref{s7eq2-6}) and (\ref{s7eq2-15})-(\ref{s7eq2-16}), we have
{\allowdisplaybreaks
\begin{align*}
\oJ_\nu(\widetilde{\ga}(t, s))&\le \oJ_0(\widetilde{\ga}(t, s))=I_{\la_1}(\widetilde{\ga}_1(t))+I_{\la_2}(\widetilde{\ga}_2(s))\\
&\le M_1/4+M_2\le M_1+M_2-3\e_0<a_\nu-\e_0,
\end{align*}
}%
which implies that $T(t, s)=0$ and so $\ga(t, s)=\widetilde{\ga}(t, s)$.
From the definition of $\G$ in (\ref{s7eq2-8}), it
suffices to prove that $\|\ga(t, s)\|_\BD\le 2 S(\la_2)^{N/4} + \mathcal{C}$ for
all $(t, s)\in Q$ and $T(t, s)$ is continuous with respect to $(t,s)$.

For any $(t, s)\in Q$, if $\oJ_{\nu}(\widetilde{\ga}(t, s))\le
a_{\nu}-\e_0$, we have $T(t, s) =0$ and so $\ga(t, s)=\widetilde{\ga}(t, s)$,
and by (\ref{s7eq2-7}) we see that $\|\ga(t, s)\|_\BD\le
\mathcal{C}<2S(\la_2)^{N/4} +\mathcal{C}$.

If $\oJ_{\nu}(\widetilde{\ga}(t, s))> a_{\nu}-\e_0$, then
$\widetilde{\ga}(t, s)\in X^{\dd/2}$ and
$$a_{\nu}-\e_0<\oJ_{\nu}(\psi_\nu(\widetilde{\ga}(t, s), \q))\le d_{\nu}
< a_{\nu}+\e_0, \quad\forall\, \q\in [0, T(t, s)),$$
This implies
$\xi_\nu(\oJ_{\nu}(\psi_\nu(\widetilde{\ga}(t, s), \q)))\equiv 1$ for $\q\in
[0, T(t, s))$. If $\psi_\nu(\widetilde{\ga}(t, s), T(t, s)) \not\in X^\dd$,
then there exists $0<\q_{t, s}^1<\q_{t, s}^2< T(t, s)$ as above. Then we can prove
that $\oJ_{\nu}(\psi_\nu(\widetilde{\ga}(t, s), \q_{t, s}^2))\le
a_{\nu}-\e_0$ as above, which contradicts the definition of
$T(t, s)$. Therefore, $\ga(t, s)=\psi_\nu(\widetilde{\ga}(t, s), T(t, s))\in
X^\dd$. Then there exists $(u, v)\in X$ such that
$\|\ga(t, s)-(u, v)\|_\BD\le 2\dd\le\mathcal{C}$.
By (\ref{s7eq2-2}) we have
$$\|\ga(t, s)\|_\BD\le\|(u, v)\|_\BD+\mathcal{C}\le 2S(\la_2)^{N/4}+\mathcal{C}.$$

To prove the continuity of $T(t, s)$, we fix any $(\tilde{t}, \tilde{s})\in Q$.
Assume $\oJ_{\nu}(\ga(\tilde{t}, \tilde{s}))< a_{\nu}-\e_0$ first. Then
$T(\tilde{t}, \tilde{s})=0$ from the definition of $T(t, s)$. So
$\oJ_{\nu}(\widetilde{\ga}(\tilde{t}, \tilde{s}))< a_{\nu}-\e_0$. By the continuity
of $\widetilde{\ga}$, there exists $\tau>0$ such that for any $(t, s)\in
(\tilde{t}-\tau, \tilde{t}+\tau)\times(\tilde{s}-\tau, \tilde{s}+\tau)\cap Q$, there holds
$\oJ_{\nu}(\widetilde{\ga}(t, s))<a_{\nu}-\e_0$, that is, $T(t, s)=0$ for any $(t, s)\in
(\tilde{t}-\tau, \tilde{t}+\tau)\times(\tilde{s}-\tau, \tilde{s}+\tau)\cap Q$,
and so $T$ is continuous at $(\tilde{t}, \tilde{s})$. Now we assume that
$\oJ_{\nu}(\ga(\tilde{t}, \tilde{s}))= a_{\nu}-\e_0$, then from the previous
proof we see that $\ga(\tilde{t}, \tilde{s})=\psi_\nu(\widetilde{\ga}(\tilde{t}, \tilde{s}), T(\tilde{t}, \tilde{s}))\in
X^\dd$ and so
$$\left\|\oJ_{\nu}' \left(\psi_\nu\left(\widetilde{\ga}(\tilde{t}, \tilde{s}),
T(\tilde{t}, \tilde{s})\right)\right)\right\| \ge l(\nu)>0.$$
Then for any $\omega>0$, we have
$\oJ_{\nu}(\psi_\nu(\widetilde{\ga}(\tilde{t}, \tilde{s}), T(\tilde{t}, \tilde{s})+\omega))<
a_{\nu}-\e_0$. By the continuity of $\psi_\nu$ and $\widetilde{\ga}$, there exists
$\tau>0$ such that for any $(t, s)\in
(\tilde{t}-\tau, \tilde{t}+\tau)\times(\tilde{s}-\tau, \tilde{s}+\tau)\cap Q$, we have $\oJ_{\nu}(\psi_\nu(\widetilde{\ga}(t, s), T(\tilde{t}, \tilde{s})+\omega))<
a_{\nu}-\e_0$, so $T(t, s)\le T(\tilde{t}, \tilde{s})+\omega$. It follows that
$$0\le
\limsup\limits_{(t,s)\to (\tilde{t}, \tilde{s})}T(t, s)\le T(\tilde{t}, \tilde{s}).$$
If $T(\tilde{t}, \tilde{s})=0$, we have
$\lim_{(t, s)\to (\tilde{t}, \tilde{s})}T(t, s)= T(\tilde{t}, \tilde{s})$ immediately.
If $T(\tilde{t}, \tilde{s})>0$,
then for any $0<\omega< T(\tilde{t}, \tilde{s})$, similarly we have
$\oJ_{\nu}(\psi_\nu(\widetilde{\ga}(\tilde{t}, \tilde{s}), T(\tilde{t}, \tilde{s})-\omega))) >
a_{\nu}-\e_0$. By the continuity of $\psi_\nu$ and $\widetilde{\ga}$ again,  we easily obtain that
 $$\liminf\limits_{(t, s)\to (\tilde{t}, \tilde{s})} T(t, s)\ge T(\tilde{t}, \tilde{s}).$$ So $T$ is
continuous at $(\tilde{t}, \tilde{s})$. This completes the proof of Step 2.

Now, we have proved that $\ga(t, s) \in \G$ and
$\max\limits_{(t,s)\in Q}\oJ_{\nu}(\ga(t, s))\le a_{\nu}-\e_0$,
which contradicts the definition of $a_{\nu}$. This
completes the proof.\hfill$\square$\\

\noindent {\bf Proof of Theorem \ref{th1-4}. } Fix any $\nu\in (0, \nu_2]$. By Lemma \ref{s7lemma5}
there exists $\{(u_n, v_n)\}_{n=1}^\iy\subset
X^\dd \cap \oJ_{\nu}^{d_\nu}$ such that
$$\left\|\oJ_\nu' (u_n, v_n)\right\|\to 0 \quad\hbox{as}\,\, n\to \iy.$$
Note that there exist $\mu_{i, n}>0, i=1, 2, n\in\mathbb{N}$ such that
\be\label{eq2-21}\left\|(u_n, v_n)-\left(z_{\mu_{1, n}}^1, z_{\mu_{2, n}}^2\right)\right\|_{\BD}\le 2\dd,\quad\forall\,n\in\mathbb{N}.\ee

By (\ref{s7eq2-2}), $\{(u_n, v_n), n\ge 1\}$ are uniformly bounded in $\BD_r$. Up to a subsequence, we may assume that
$(u_n, v_n)\rightharpoonup (u, v)$ weakly in $\BD_r$. As pointed out before,
since $X$ is not compact in $\BD_r$, it seems very difficult to prove that $u\not\equiv 0$ and $v\not\equiv 0$.
To overcome this difficulty, let us define
{\allowdisplaybreaks
\begin{align}\label{eq2-20}&\widetilde{u}_n(x):=\mu_{1, n}^{\frac{N-2}{2}}u_n(\mu_{1, n}x),\quad \widetilde{v}_n(x):=\mu_{1, n}^{\frac{N-2}{2}}v_n(\mu_{1, n}x);\nonumber\\
&\ou_n(x):=\mu_{2, n}^{\frac{N-2}{2}}u_n(\mu_{2, n}x),\quad \ov_n(x):=\mu_{2, n}^{\frac{N-2}{2}}v_n(\mu_{2, n}x).
\end{align}
}%
Then by a direct computation, we see that $\|\oJ_\nu' (\wu_n, \wv_n)\|\to 0$ and $\|\oJ_\nu' (\ou_n, \ov_n)\|\to 0$ as $n\to\iy$.
Moreover,
\begin{align*}
\|\widetilde{u}_n-z_1^1\|_{\la_1}\le 2\dd,\quad \|\ov_n-z_1^2\|_{\la_2}\le 2\dd.
\end{align*}
Up to a subsequence,
we may assume that $(\wu_n, \wv_n)\rightharpoonup (\wu, \,\wv)$ and $(\ou_n, \ov_n)\rightharpoonup (\ou, \,\ov)$ weakly in $\BD_r\cap L^{2^\ast}(\RN)\times L^{2^\ast}(\RN)$. Then $\oJ_\nu' (\wu, \wv)= 0$ and $\oJ'_\nu (\ou, \ov)=0$. Moreover, as in the proof of Lemma \ref{s7lemma3},
we get that $\wu\not\equiv 0$ and $\ov\not\equiv 0$.

Now we claim that either $\wv\not\equiv 0$ or $\ou\not\equiv 0$.

Assume by contradiction that both $\wv\equiv 0$ and $\ou\equiv 0$. Then $\wv_n\rightharpoonup 0$ weakly in $D_r^{1, 2}(\RN)\cap L^{2^\ast}(\RN)$.
Hence it is easy to prove that $(\wv_n)_+^{2^\ast-1}\rightharpoonup 0$
in $L^{\frac{2^\ast}{2^\ast-1}}(\RN)$.
Take any $\wphi\in C_{0, r}^\iy(\RN)$ such that $\|\wphi\|_{\la_1}=1$, and $\phi_n$ is defined in (\ref{test}).
Then we see from H\"{o}der inequality and $\al+\bb=2^\ast$ that
{\allowdisplaybreaks
\begin{align}\label{s7eq2-18}
&\lim_{n\to\iy}\left|\nu \al\intR{(u_n)_+^{\al-1}(v_n)_+^{\bb}\phi_n}\right|
=\lim_{n\to\iy}\left|\nu \al\intR{(\wu_n)_+^{\al-1}(\wv_n)_+^{\bb}\wphi}\right|\nonumber\\
\le &\lim_{n\to\iy}C\left(\intR{(\wu_n)_+^{2^\ast-1}|\wphi|}\right)^{\frac{\al-1}{2^\ast-1}}
\left(\intR{(\wv_n)_+^{2^\ast-1}|\wphi|}\right)^{\frac{\bb}{2^\ast-1}}=0.\end{align}
}%
Combining this with $\oJ_\nu'(u_n, v_n)(\phi_n, 0)\to 0$, we may repeat the proof of Lemma \ref{s7lemma3} and then get that $\wu\in Z_1$.
Similarly, we can prove that $\ov\in Z_2$. Then
{\allowdisplaybreaks
\begin{align}\label{s7eq2-19}
d_\nu &\ge\lim_{n\to\iy}\left(\oJ_{\nu}(u_n, v_n)-\frac{1}{2^\ast}\oJ_{\nu}'(u_n, v_n)(u_n, v_n)\right)\nonumber\\
&=\frac{1}{N}\lim_{n\to\iy}\left(\|u_n\|^2_{\la_1}+\|v_n\|_{\la_2}^2\right)=\frac{1}{N}\lim_{n\to\iy}
\left(\|\wu_n\|^2_{\la_1}+\|\ov_n\|_{\la_2}^2\right)\nonumber\\
&\ge\frac{1}{N}\|\wu\|_{\la_1}^2+\frac{1}{N}\|\ov\|_{\la_2}^2=M_1+M_2=d_0>d_\nu,
\end{align}
}%
a contradiction. So either $\wv\not\equiv 0$ or $\ou\not\equiv 0$. Without loss of generality, we may assume that $\wv\not\equiv 0$.
Note that $\wu\not\equiv 0$ and $\oJ_\nu'(\wu, \wv)=0$,
then by testing (\ref{s7eq2-0}) with $\wu_-$ and $\wv_-$, we see from Hardy inequality (\ref{Hardy}) that $\wu\ge 0$ and $\wv\ge 0$.
By the maximum principle, one has that $\wu>0$ and $\wv>0$ in $\RN\setminus\{0\}$.
Hence, $(\wu, \wv)$ is a positive solution of (\ref{eq0-0}), which is radially symmetric. Moreover,
{\allowdisplaybreaks
\begin{align*}
d_\nu &\ge\lim_{n\to\iy}\left(\oJ_{\nu}(u_n, v_n)-\frac{1}{2^\ast}\oJ_{\nu}'(u_n, v_n)(u_n, v_n)\right)\nonumber\\
&=\frac{1}{N}\lim_{n\to\iy}\left(\|u_n\|^2_{\la_1}+\|v_n\|_{\la_2}^2\right)=\frac{1}{N}\lim_{n\to\iy}
\left(\|\wu_n\|^2_{\la_1}+\|\wv_n\|_{\la_2}^2\right)\nonumber\\
&\ge\frac{1}{N}\|\wu\|_{\la_1}^2+\frac{1}{N}\|\wv\|_{\la_2}^2=J_\nu(\wu, \wv),
\end{align*}
}%
that is, $J_\nu(\wu, \wv)\le d_\nu<d_0=\frac{1}{N}(S(\la_1)^{N/2}+S(\la_2)^{N/2})$. This completes the proof. \hfill$\square$

\vskip0.1in

\s{Proof of Theorem \ref{theorem2}: The case $N=3$}
\renewcommand{\theequation}{5.\arabic{equation}}

In this section, we give the proof of Theorem \ref{theorem2}. Assume that $N=3$, $\al+\bb=2^\ast=6$, $\al\ge 2$, $\bb\ge 2$, $\la_1, \la_2\in (0, \La_3)$ and
$\nu>0$.
Recall from (\ref{eq1-5}) that $S(\la_i)<S$, we take a $\e_0\in (0, 1/2)$ such that
\be\label{5eq1}\begin{cases}
\frac{\max\left\{S(\la_1)^{\frac{3}{2}},\,S(\la_2)^{\frac{3}{2}}\right\}}{1-\e_0}\le S(\la_1)^{\frac{3}{2}}+S(\la_2)^{\frac{3}{2}}\le(2-2\e_0)S^{\frac{3}{2}},\\
(1-2\e_0)S(\la_2)^{\frac{3}{2}}\ge \e_0 S(\la_1)^{\frac{3}{2}},\,\,\,(1-2\e_0)S(\la_1)^{\frac{3}{2}}\ge \e_0 S(\la_2)^{\frac{3}{2}}.
\end{cases}
\ee
Define
$$K_\nu:=\left\{(u, v)\in\BD \,\,: \,\,u\not\equiv0,\,\,v\not\equiv0,\,\,J_\nu'(u, v)=0\right\}$$
as the set of nontrivial critical point of $J_\nu$, and
\be\label{5eq2} b_\nu:=\inf_{(u, v)\in K_\nu}J_\nu (u, v).\ee
By Theorem \ref{th1-4} we see that for any $\nu\in (0, \nu_2]$, $K_\nu\neq\emptyset$, $b_\nu$ is well defined  and
$$b_\nu<\frac{1}{3}\left(S(\la_1)^{\frac{3}{2}}+S(\la_2)^{\frac{3}{2}}\right).$$
Note that $K_\nu\subset\N_\nu$, so $b_\nu\ge c_\nu>0$. Define
\be\label{5eq3}\mathcal{C}_0:=\max\left\{\left[S(\la_1)^{\frac{3}{2}}+S(\la_2)^{\frac{3}{2}}\right]^{\frac{\al}{6}},
\,\,\,\left[S(\la_1)^{\frac{3}{2}}+S(\la_2)^{\frac{3}{2}}\right]^{\frac{\bb}{6}}\right\}.\ee
Then for any $(u, v)\in \N_\nu$ with $J_\nu (u, v)<\frac{1}{3}(S(\la_1)^{\frac{3}{2}}+S(\la_2)^{\frac{3}{2}})$, since
$$\int_{\R^3}\left(u^6+v^6+6\nu |u|^\al|v|^\bb\right)=3J_\nu(u, v)<S(\la_1)^{\frac{3}{2}}+S(\la_2)^{\frac{3}{2}},$$
so
\be\label{5eq4}|u|_6^\al\le\mathcal{C}_0,\quad|v|_6^\bb\le\mathcal{C}_0.\ee
Then we see from H\"{o}lder inequality and (\ref{eq1-5-1}) that
$$S(\la_1)|u|_6^2 \le \|u\|^2_{\la_1}=|u|_6^6+\nu\al\int_{\R^3}|u|^\al|v|^\bb\le |u|_6^6+\nu\al\mathcal{C}_0 |u|_6^\al,$$
That is, we can obtain
\be\label{5eq5}
\begin{cases}
S(\la_1)|u|_6^2 \le  |u|_6^6+\nu\al\mathcal{C}_0 |u|_6^\al,\\
S(\la_2)|v|_6^2 \le  |v|_6^6+\nu\bb\mathcal{C}_0 |v|_6^\bb.
\end{cases}
\ee
The following result was introduced in \cite{AFP}.

\bl\label{lemma5-1} (see \cite[Lemma 3.3]{AFP}) Let $N\ge 3$, $A, B>0$, and $\theta\ge 2$ be fixed. For any $\nu>0$, let
$$S_\nu:=\left\{\sg>0 \,\,:\,\, A\sg^{\frac{2}{2^\ast}}\le \sg+ \nu B \sg^{\frac{\theta}{2^\ast}}\right\}.$$
Then for any $\e>0$, there exists $\nu_1>0$ depending only on $\e, A, B, \theta$ and $N$, such that
$$\inf S_\nu\ge (1-\e) A^{\frac{N}{2}} \quad\hbox{for all}\,\,\,0<\nu<\nu_1.$$

\el

Recall that $\al, \bb\ge 2$ and $2^\ast=6$. From Lemma \ref{lemma5-1}, we have the following result trivially.
\bl\label{lemma5-2}

Recall $\e_0$ in (\ref{5eq1}), then there exists $\widetilde{\nu}_1\in (0, \nu_2]$ such that for any $\nu\in (0, \widetilde{\nu}_1)$ there hold
{\allowdisplaybreaks
\begin{align}
\label{5eq6}&S(\la_1)\sg^{\frac{1}{3}}\le \sg+ \nu \al\mathcal{C}_0 \sg^{\frac{\al}{6}},\,\,\sg>0\,\,\Rightarrow \,\,\sg\ge(1-\e_0)S(\la_1)^{\frac{3}{2}},\\
\label{5eq7}&S(\la_2)\sg^{\frac{1}{3}}\le \sg+ \nu \bb\mathcal{C}_0 \sg^{\frac{\bb}{6}},\,\,\sg>0\,\,\Rightarrow \,\,\sg\ge(1-\e_0)S(\la_2)^{\frac{3}{2}},\\
\label{5eq8}&S\sg^{\frac{1}{3}}\le \sg+ \nu \al\mathcal{C}_0 \sg^{\frac{\al}{6}},\,\,\sg>0\,\,\Rightarrow \,\,\sg\ge(1-\e_0)S^{\frac{3}{2}},\\
\label{5eq9}&S\sg^{\frac{1}{3}}\le \sg+ \nu \bb\mathcal{C}_0 \sg^{\frac{\bb}{6}},\,\,\sg>0\,\,\Rightarrow \,\,\sg\ge(1-\e_0)S^{\frac{3}{2}}.
\end{align}
}%
\el

The following lemma is the counterpart of Lemma \ref{lemma3-3} for the case $N=3$,
and the first part of the proof is similar to that of Lemma \ref{lemma3-3}, so we do not give the details,
but the latter part of the proof is quite different.

\bl\label{lemma5-3} Assume that $\nu\in (0, \widetilde{\nu}_1)$. Let $(u_n, v_n)\in K_\nu$ be a minimizing
sequence of $b_\nu$, and $(u_n, v_n)\rightharpoonup (0, 0)$ weakly in $\BD$.
Then for any $r>0$ and for every $\e\in (-r, 0)\cup(0, r)$, there exists $\rho\in (\e, 0)\cup (0, \e)$ such that, up to a subsequence,
\be\label{5eq400}\hbox{either}\,\,\int_{B_{r+\rho}}(|\nabla u_n|^2+|\nabla v_n|^2)\to 0
\,\,\hbox{or}\,\,\int_{\RN\setminus B_{r+\rho}}(|\nabla u_n|^2+|\nabla v_n|^2)\to 0.\ee\el

\noindent{\bf Proof. } Fix any $\nu\in (0, \widetilde{\nu}_1)$. Without loss of generality, we only consider the case $\e\in (0, r)$
(the proof for the case $\e\in (-r, 0)$ is similar). Since $(u_n, v_n)\in K_\nu$ is a minimizing
sequence of $b_\nu$, we may assume that $J_\nu (u_n, v_n)<\frac{1}{3}(S(\la_1)^{\frac{3}{2}}+S(\la_2)^{\frac{3}{2}})$ and so $(u_n, v_n)$
satisfy (\ref{5eq4})-(\ref{5eq5}) for all $n$. Then by (\ref{5eq6})-(\ref{5eq7}) of Lemma \ref{lemma5-2} we have
\be\label{5eq150}|u_n|_6^6\ge (1-\e_0)S(\la_1)^{\frac{3}{2}},\quad |v_n|_6^6\ge (1-\e_0)S(\la_2)^{\frac{3}{2}},\quad\forall\,n\in\mathbb{N}.\ee

Note that $(u_n, v_n)$ are uniformly bounded in $\BD$ and $J_\nu'(u_n, v_n)=0$.
Then by repeating the argument of Lemma \ref{lemma3-3} with trivial modifications and
using the same notations $w_{i, n}, \sg_{i, n}, u_{i, n}, v_{i, n}$ with the same definitions as (\ref{5eq10})-(\ref{5eq13}),
there exists $\rho\in (0, \e)$ such that $u_{i, n}, v_{i, n}, i=1, 2$ satisfy (\ref{eq4-2})-(\ref{eq4-7}). Moreover, we can prove that
{\allowdisplaybreaks
\begin{align}\label{5eq4-2}|u_n|_{6}^6=|u_{1, n}|_{6}^6+|u_{2,n}|_{6}^6+o(1),\quad
 |v_n|_{6}^6=|v_{1, n}|_{6}^6+|v_{2,n}|_{6}^6+o(1).\end{align}
}%
Now we claim that
\be\label{5eq4-8}\hbox{either}\,\,\lim_{n\to\iy}(\|u_{1, n}\|^2+\|v_{1, n}\|^2)=0\,\,\hbox{or}\,\,\lim_{n\to\iy}(\|u_{2, n}\|^2+\|v_{2, n}\|^2)=0.\ee

  In fact, if (\ref{5eq4-8}) does not hold, then up to a subsequence,
\be\label{5eq44-8}\hbox{both}\,\,\lim_{n\to\iy}(\|u_{1, n}\|^2+\|v_{1, n}\|^2)>0\,\,\hbox{and}\,\,\lim_{n\to\iy}(\|u_{2, n}\|^2+\|v_{2, n}\|^2)>0.\ee
We have the following several cases.

{\bf Case 1.} Up to a subsequence, both $\lim\limits_{n\to\iy}\|u_{1, n}\|^2>0$ and $\lim\limits_{n\to\iy}\|v_{1, n}\|^2>0$.

Recall that norms $\|\cdot\|_{\la_i}, i=1, 2$ are equivalent to $\|\cdot\|$.  Note that (\ref{eq4-4})-(\ref{eq4-5}) yield
{\allowdisplaybreaks
\begin{align}
\label{5eq4-9}\|u_{1, n}\|_{\la_1}^2=|u_{1, n}|_6^6 + \nu \al\int_{\R^3}|u_{1, n}|^\al |v_{1, n}|^\bb+o(1),\\
\label{5eq4-10}\|v_{1, n}\|_{\la_2}^2=|v_{1, n}|_6^6 + \nu\bb \int_{\R^3}|u_{1, n}|^\al |v_{1, n}|^\bb+o(1).
\end{align}
}%
Hence, both $A_1:=\liminf\limits_{n\to\iy}|u_{1, n}|_6^6>0$ and $B_1:=\liminf\limits_{n\to\iy}|v_{1, n}|_6^6>0$. Since $(u_n, v_n)$
satisfy (\ref{5eq4}) for all $n$, by (\ref{5eq4-2}) and letting $n\to\iy$ in (\ref{5eq4-9})-(\ref{5eq4-10}), similarly as (\ref{5eq5}) we can prove that
\be\label{5eq14}
S(\la_1)A_1^{\frac{1}{3}} \le  A_1+\nu\al\mathcal{C}_0 A_1^{\frac{\al}{6}},\quad
S(\la_2)B_1^{\frac{1}{3}} \le  B_1+\nu\bb\mathcal{C}_0 B_1^{\frac{\bb}{6}}.
\ee
Then by  (\ref{5eq6})-(\ref{5eq7}) of Lemma \ref{lemma5-2} we have
\be\label{5eq15}A_1\ge (1-\e_0)S(\la_1)^{\frac{3}{2}},\quad B_1\ge (1-\e_0)S(\la_2)^{\frac{3}{2}}.\ee

{\bf Case 1.1.} Up to a subsequence, $\lim\limits_{n\to\iy}\|u_{2, n}\|^2>0$ and $\lim\limits_{n\to\iy}\|v_{2, n}\|^2>0$.

Then similarly as above, we can prove that
\be\label{5eq16}A_2:=\liminf\limits_{n\to\iy}|u_{2, n}|_6^6\ge (1-\e_0)S(\la_1)^{\frac{3}{2}},\,\,\, B_2:=\liminf\limits_{n\to\iy}|v_{2, n}|_6^6\ge (1-\e_0)S(\la_2)^{\frac{3}{2}}.\ee
Combining this with (\ref{5eq15}), (\ref{5eq1}) and (\ref{5eq4-2}), we deduce that
{\allowdisplaybreaks
\begin{align*}
b_\nu &=\lim_{n\to\iy}J_\nu(u_n, v_n)= \lim_{n\to\iy}\frac{1}{3}\left(|u_n|_{6}^6+|v_n|_{6}^6+6\nu\int_{\R^3}|u_n|^\al|v_n|^\bb\right)\\
&\ge\lim_{n\to\iy}\frac{1}{3}(|u_n|_{6}^6+|v_n|_{6}^6)\ge\frac{1}{3}(A_1+B_1+A_2+B_2)\\
&\ge\frac{2-2\e_0}{3}\left(S(\la_1)^{\frac{3}{2}}+S(\la_2)^{\frac{3}{2}}\right)\\
&\ge\frac{1}{3}\left(S(\la_1)^{\frac{3}{2}}+S(\la_2)^{\frac{3}{2}}\right)> b_\nu,
\end{align*}
}%
a contradiction. So Case 1.1 is impossible.

{\bf Case 1.2.} Up to a subsequence, $\lim\limits_{n\to\iy}\|u_{2, n}\|^2>0$ and $\lim\limits_{n\to\iy}\|v_{2, n}\|^2=0$.

Then $\int_{\R^3}|u_{2, n}|^\al|v_{2, n}|^\bb\to 0$ as $n\to\iy$, so (\ref{eq4-6}) yields
$$S(\la_1)|u_{2, n}|_6^2\le \|u_{2, n}\|_{\la_1}^2=|u_{2, n}|_6^6+o(1),$$
and so $A_2\ge S(\la_1)^{3/2}$. Then we conclude from (\ref{5eq15}), (\ref{5eq1}) and (\ref{5eq4-2}) that
{\allowdisplaybreaks
\begin{align*}
b_\nu
&\ge\lim_{n\to\iy}\frac{1}{3}(|u_n|_{6}^6+|v_n|_{6}^6)\ge\frac{1}{3}(A_1+B_1+A_2)\\
&\ge\frac{1-\e_0}{3}\left(S(\la_1)^{\frac{3}{2}}+S(\la_2)^{\frac{3}{2}}\right)+\frac{1}{3}S(\la_1)^{\frac{3}{2}}\\
&\ge\frac{1}{3}\left(S(\la_1)^{\frac{3}{2}}+S(\la_2)^{\frac{3}{2}}\right)> b_\nu,
\end{align*}
}%
a contradiction. So Case 1.2 is impossible.

{\bf Case 1.3.} Up to a subsequence, $\lim\limits_{n\to\iy}\|u_{2, n}\|^2=0$ and $\lim\limits_{n\to\iy}\|v_{2, n}\|^2>0$.

Then (\ref{eq4-7}) yields
$$S(\la_2)|v_{2, n}|_6^2\le \|v_{2, n}\|_{\la_2}^2=|v_{2, n}|_6^6+o(1),$$
and so $B_2\ge S(\la_2)^{3/2}$. Then we conclude from (\ref{5eq15}), (\ref{5eq1}) and (\ref{5eq4-2}) that
{\allowdisplaybreaks
\begin{align*}
b_\nu
&\ge\lim_{n\to\iy}\frac{1}{3}(|u_n|_{6}^6+|v_n|_{6}^6)\ge\frac{1}{3}(A_1+B_1+B_2)\\
&\ge\frac{1-\e_0}{3}\left(S(\la_1)^{\frac{3}{2}}+S(\la_2)^{\frac{3}{2}}\right)+\frac{1}{3}S(\la_2)^{\frac{3}{2}}\\
&\ge\frac{1}{3}\left(S(\la_1)^{\frac{3}{2}}+S(\la_2)^{\frac{3}{2}}\right)> b_\nu,
\end{align*}
}%
a contradiction. So Case 1.3 is impossible.

Since none of Cases 1.1-1.3 is true, so Case 1 is impossible.

{\bf Case 2.} Up to a subsequence, $\lim\limits_{n\to\iy}\|u_{1, n}\|^2=0$ and $\lim\limits_{n\to\iy}\|v_{1, n}\|^2>0$.

Then similarly as Case 1.3, we have $B_1\ge S(\la_2)^{3/2}$. Moreover, we see from (\ref{5eq150}) and (\ref{5eq4-2}) that
$$\liminf_{n\to\iy}|u_{2,n}|_6^6=\liminf_{n\to\iy}|u_n|_6^6-\lim_{n\to\iy}|u_{1,n}|_6^6>0.$$

{\bf Case 2.1.} Up to a subsequence, $\lim\limits_{n\to\iy}\|u_{2, n}\|^2>0$ and $\lim\limits_{n\to\iy}\|v_{2, n}\|^2>0$.

Then similarly as above, we see that (\ref{5eq16}) holds, and so
{\allowdisplaybreaks
\begin{align*}
b_\nu
&\ge\lim_{n\to\iy}\frac{1}{3}(|u_n|_{6}^6+|v_n|_{6}^6)\ge\frac{1}{3}(B_1+A_2+B_2)\\
&\ge\frac{1-\e_0}{3}\left(S(\la_1)^{\frac{3}{2}}+S(\la_2)^{\frac{3}{2}}\right)+\frac{1}{3}S(\la_2)^{\frac{3}{2}}\\
&\ge\frac{1}{3}\left(S(\la_1)^{\frac{3}{2}}+S(\la_2)^{\frac{3}{2}}\right)> b_\nu,
\end{align*}
}%
a contradiction. So Case 2.1 is impossible.

{\bf Case 2.2.} Up to a subsequence, $\lim\limits_{n\to\iy}\|u_{2, n}\|^2>0$ and $\lim\limits_{n\to\iy}\|v_{2, n}\|^2=0$.

Then similarly as Case 1.2, we have $A_2\ge S(\la_1)^{3/2}$, and so
{\allowdisplaybreaks
\begin{align*}
b_\nu
&\ge\lim_{n\to\iy}\frac{1}{3}(|u_n|_{6}^6+|v_n|_{6}^6)\ge\frac{1}{3}(B_1+A_2)\\
&\ge\frac{1}{3}\left(S(\la_1)^{\frac{3}{2}}+S(\la_2)^{\frac{3}{2}}\right)> b_\nu,
\end{align*}
}%
a contradiction. So Case 2.2 is impossible.

Since neither Case 2.1 or Case 2.2 is true, so Case 2 is impossible.

{\bf Case 3.} Up to a subsequence, $\lim\limits_{n\to\iy}\|u_{1, n}\|^2>0$ and $\lim\limits_{n\to\iy}\|v_{1, n}\|^2=0$.

By a similar argument as Case 2, we get a contradiction. So Case 3 is impossible.

Since none of Cases 1, 2 and 3 is true, we see that (\ref{5eq44-8}) is impossible, that is, (\ref{5eq4-8}) holds.
Recall the definitions (\ref{5eq12})-(\ref{5eq13}) of $(u_{i, n}, v_{i, n})$, (\ref{5eq400}) follows directly from (\ref{5eq4-8}). This completes the proof.\hfill$\square$\\

\noindent {\bf Proof of Theorem \ref{theorem2}. } Fix any $\nu\in (0, \widetilde{\nu}_1)$.
Take a sequence $(\ou_n, \ov_n)\in K_\nu$ such
that $J_\nu(\ou_n, \ov_n) \to b_\nu$ as $n\to\iy$. Recall that
$E(u, v)= |\nabla u|^2+|\nabla v|^2-\frac{\la_1}{|x|^2}|u|^2-\frac{\la_2}{|x|^2}|v|^2$, there exists $R_n>0$ such that
$$\int_{B_{R_n}} E(\ou_n, \ov_n)=\int_{\RN\setminus B_{R_n}} E(\ou_n, \ov_n)=\frac{1}{2}(\|\ou_n\|_{\la_1}^2+\|\ov_n\|_{\la_2}^2).$$
Define
$$(u_n(x), v_n(x)):=\left(R_n^{\frac{N-2}{2}}\ou_n(R_n x), \, R_n^{\frac{N-2}{2}}\ov_n(R_n x)\right),$$
Then by a direct computation, we see that $(u_n, v_n)\in K_\nu$ and $J_\nu(u_n, v_n) \to b_\nu$. Moreover,
\be\label{5eq4-11}\int_{B_{1}} E(u_n, v_n)=\int_{\RN\setminus B_{1}} E(u_n, v_n)=\frac{1}{2}(\|u_n\|_{\la_1}^2+\|v_n\|_{\la_2}^2)\to \frac{3}{2} b_\nu>0.\ee
 Besides,  we may assume that $J_\nu (u_n, v_n)<\frac{1}{3}(S(\la_1)^{\frac{3}{2}}+S(\la_2)^{\frac{3}{2}})$ and so $(u_n, v_n)$
satisfy (\ref{5eq4})-(\ref{5eq5}) for all $n$. Then by (\ref{5eq6})-(\ref{5eq7}) of Lemma \ref{lemma5-2} we have
\be\label{5eq1500}|u_n|_6^6\ge (1-\e_0)S(\la_1)^{\frac{3}{2}},\quad |v_n|_6^6\ge (1-\e_0)S(\la_2)^{\frac{3}{2}},\quad\forall\,n\in\mathbb{N}.\ee

Note that $(u_n, v_n)$ are uniformly bounded in $\BD$. Then up to a subsequence, we assume that $(u_n, v_n)\rightharpoonup (u, v)$ weakly in $\BD$.
Then $J_\nu'(u_n, v_n)=0$ implies $J_\nu'(u, v)=0$.\\

\noindent {\bf Step 1.} We show that both $u\not\equiv 0$ and $v\not\equiv 0$, that is $(u, v)\in K_\nu$. Moreover $J_\nu (u, v)=b_\nu$.

{\bf Case 1.} $(u, v)\equiv (0, 0)$.

Then we can apply twice Lemma \ref{lemma5-3} with $r=1$ and $\e=\pm 1/4$ respectively, and there exist $\rho^+\in (0, 1/4)$ and $\rho^-\in (-1/4, 0)$
such that the alternative (\ref{5eq400}) holds. Then by repeating the argument of Case 1 in the proof of Theorem \ref{th2} with trivial modifications, we get that $\int_{\R^3}\frac{u_n^2}{|x|^2}=o(1)$, $\int_{\R^3}\frac{v_n^2}{|x|^2}=o(1)$ and so
{\allowdisplaybreaks
\begin{align*}
&S|u_n|_6^2\le\int_{\R^3}|\nabla u_n|^2=\int_{\R^3}u_n^6+\nu\al\int_{\R^3}|u_n|^\al |v_n|^\bb+ o(1),\\
&S|v_n|_6^2\le\int_{\R^3}|\nabla v_n|^2=\int_{\R^3}v_n^6+\nu\bb\int_{\R^3}|u_n|^\al |v_n|^\bb+ o(1).
\end{align*}
}%
Denote $A=\liminf_{n\to\iy}|u_n|_6^6$ and $B=\liminf_{n\to\iy}|v_n|_6^6$, then (\ref{5eq1500}) yields $A>0$ and $B>0$. Then by H\"{o}lder inequality it is easy to prove that
$$
SA^{\frac{1}{3}} \le  A+\nu\al\mathcal{C}_0 A^{\frac{\al}{6}},\quad
SB^{\frac{1}{3}} \le  B+\nu\bb\mathcal{C}_0 B^{\frac{\bb}{6}}.
$$
Then by  (\ref{5eq8})-(\ref{5eq9}) of Lemma \ref{lemma5-2} we have
$$ A\ge (1-\e_0)S^{\frac{3}{2}},\quad B\ge (1-\e_0)S^{\frac{3}{2}}.$$
So we conclude from (\ref{5eq1}) that
{\allowdisplaybreaks
\begin{align*}
b_\nu
&\ge\lim_{n\to\iy}\frac{1}{3}(|u_n|_{6}^6+|v_n|_{6}^6)\ge\frac{1}{3}(A+B)
\ge\frac{2-2\e_0}{3}S^{\frac{3}{2}}\\
&\ge\frac{1}{3}\left(S(\la_1)^{\frac{3}{2}}+S(\la_2)^{\frac{3}{2}}\right)> b_\nu,
\end{align*}
}%
a contradiction. So Case 1 is impossible.

{\bf  Case 2.} Either $u\equiv 0,\,v\not\equiv 0$ or $u\not\equiv 0, v\equiv 0$.

Without loss of generality, we assume that $u\not\equiv 0, v\equiv 0$. We see from $J_\nu'(u, v)(u, 0)=0$ that
$$S(\la_1)|u|_6^2\le\|u\|_{\la_1}^2=|u|_6^6,$$
which implies $|u|_{6}^6\ge S(\la_1)^{3/2}$.

{\bf Case 2.1.} Up to a subsequence, $\lim_{n\to\iy}\|u_n-u\|>0$.

Denote $w_n=u_n-u$. Note that $J_\nu'(u_n, v_n)=0$. Then by Brezis-Lieb Lemma (\cite{BL}) and Lemma \ref{lemma3} we conclude that
{\allowdisplaybreaks
\begin{align*}
&\|w_n\|_{\la_1}^2=\int_{\R^3}w_n^6+\nu\al\int_{\R^3}|w_n|^\al |v_n|^\bb+ o(1),\\
&\|v_n\|_{\la_2}^2=\int_{\R^3}v_n^6+\nu\bb\int_{\R^3}|w_n|^\bb |v_n|^\bb+ o(1).
\end{align*}
}%
Denote $C=\liminf_{n\to\iy}|w_n|_6^6$, then $C>0$. Then by H\"{o}lder inequality it is easy to prove that
$$
S(\la_1)C^{\frac{1}{3}} \le  C+\nu\al\mathcal{C}_0 C^{\frac{\al}{6}},\quad
S(\la_2) B^{\frac{1}{3}} \le  B+\nu\bb\mathcal{C}_0 B^{\frac{\bb}{6}}.
$$
Then by  (\ref{5eq6})-(\ref{5eq7}) of Lemma \ref{lemma5-2} we have
$$ C\ge (1-\e_0)S(\la_1)^{\frac{3}{2}},\quad B\ge (1-\e_0)S(\la_2)^{\frac{3}{2}}.$$
So we conclude from (\ref{5eq1}) that
{\allowdisplaybreaks
\begin{align*}
b_\nu
&\ge\lim_{n\to\iy}\frac{1}{3}(|u_n|_{6}^6+|v_n|_{6}^6)=\frac{1}{3}|u|_6^6+\lim_{n\to\iy}\frac{1}{3}(|w_n|_{6}^6+|v_n|_{6}^6)\\
&\ge\frac{1}{3}S(\la_1)^{\frac{3}{2}}+\frac{1}{3}(B+C)
\ge\frac{1}{3}S(\la_1)^{\frac{3}{2}}+\frac{1-\e_0}{3}\left(S(\la_1)^{\frac{3}{2}}+S(\la_2)^{\frac{3}{2}}\right)\\
&\ge\frac{1}{3}\left(S(\la_1)^{\frac{3}{2}}+S(\la_2)^{\frac{3}{2}}\right)> b_\nu,
\end{align*}
}%
a contradiction. So Case 2.1 is impossible.

{\bf Case 2.2.} $u_n\to u$ strongly in $D^{1, 2}(\R^3)$.

Then $u_n\to u$ strongly in $L^{6}(\R^3)$. Recall that $v_n\rightharpoonup 0$ in $D^{1, 2}(\R^3)$, up to a subsequence, $u_n\to u$ and $v_n\to 0$ almost everywhere in $\R^3$. So Lemma \ref{lemma3} yields
$$\int_{\R^3}|u_n|^\al |v_n|^\bb=\int_{\R^3}|u_n-u|^\al|v_n|^\bb+ o(1)=o(1).$$
Then we have
$$S(\la_2)|v_n|_6^2\le\|v_n\|_{\la_2}^2=|v_n|_6^6+o(1),$$
so $B\ge S(\la_2)^{3/2}$, and we conclude from (\ref{5eq1}) that
{\allowdisplaybreaks
\begin{align*}
b_\nu
&\ge\lim_{n\to\iy}\frac{1}{3}(|u_n|_{6}^6+|v_n|_{6}^6)\ge\frac{1}{3}(|u|_6^6+B)
\ge\frac{1}{3}\left(S(\la_1)^{\frac{3}{2}}+S(\la_2)^{\frac{3}{2}}\right)> b_\nu,
\end{align*}
}%
a contradiction. So Case 2.2 is impossible, and so Case 2 is impossible.

Since neither Case 1 nor Case 2 is true, we obtain that $u\not\equiv 0$ and $v\not\equiv 0$. Since $J'_\nu(u, v)=0$, so $(u, v)\in K_\nu$.
Then
$$b_\nu\le J_\nu(u, v)=\frac{1}{3}\|(u, v)\|_{\BD}^2\le \liminf_{n\to\iy} \frac{1}{3}\|(u_n, v_n)\|_{\BD}^2=\liminf_{n\to\iy}J_\nu(u_n, v_n)=b_\nu,$$
so $J_\nu(u, v)=b_\nu$, and $(u_n, v_n)\to (u, v)$ strongly in $\BD$.
Then (\ref{5eq1500}) implies that
\be\label{5eq17}|u|_6^6\ge (1-\e_0)S(\la_1)^{\frac{3}{2}},\quad |v|_6^6\ge (1-\e_0)S(\la_2)^{\frac{3}{2}}.\ee

\noindent {\bf Step 2.} We show that neither $u$ or $v$ is sign-changing, so $(|u|, |v|)$ is a positive ground state solution of (\ref{eq0-0}).

Assume by contradiction that $u_+\not\equiv 0$ and $u_-\not\equiv 0$. By $J_\nu'(u, v)( u_{\pm}, 0)=0$ we obtain
$$S(\la_1)|u_{\pm}|_6^2\le\|u_{\pm}\|_{\la_1}^2=\int_{\R^3}u_{\pm}^6+\nu\al\int_{\R^3}|u_{\pm}|^\al |v|^\bb.$$
Then
$$S(\la_1)|u_{\pm}|_6^2\le |u_{\pm}|_6^6+\nu\al\mathcal{C}_0 |u_{\pm}|_6^\al,$$
By (\ref{5eq6}) of Lemma \ref{lemma5-2} we have
$$ |u_{\pm}|_6^6\ge (1-\e_0)S(\la_1)^{\frac{3}{2}}.$$
So we conclude from (\ref{5eq1}) and (\ref{5eq17}) that
{\allowdisplaybreaks
\begin{align*}
b_\nu
&\ge\frac{1}{3}(|u|_{6}^6+|v|_{6}^6)\ge\frac{2-2\e_0}{3}S(\la_1)^{\frac{3}{2}}+\frac{1-\e_0}{3}S(\la_2)^{\frac{3}{2}}\\
&\ge\frac{1}{3}\left(S(\la_1)^{\frac{3}{2}}+S(\la_2)^{\frac{3}{2}}\right)> b_\nu,
\end{align*}
}%
a contradiction. So $u$ is not sign-changing. Similarly, $v$ is not sign-changing. That is, $(|u|, |v|)$ is a solution of $J_\nu$. By the maximum principle, we see that $|u|>0$ and $|v|>0$ in $\R^3\setminus\{0\}$. Since $J_\nu(|u|, |v|)=b_\nu$,
so $(|u|, |v|)$ is a positive ground state solution of (\ref{eq0-0}).\\

\noindent {\bf Step 3.} We show that $b_\nu\to
\frac{1}{3}\left(S(\la_1)^{3/2}+S(\la_2)^{3/2}\right)$ as $\nu\to 0$.

From the above argument, we may assume that $(u_\nu, v_\nu)$ is a positive ground state solution of (\ref{eq0-0}) with $b_\nu=J_\nu(u_\nu, v_\nu)$ for any $\nu\in (0, \widetilde{\nu}_1)$. The rest argument is similar to that in the proof of Theorem \ref{th2}, and we omit the details.
This completes the proof.\hfill$\square$

\vskip0.1in

\s{Proof of Theorem \ref{th4}: The moving planes method}
\renewcommand{\theequation}{6.\arabic{equation}}

In this section, we will use the moving planes method to prove Theorem \ref{th4}. In the sequel, we assume
that $N=3 $ or $N= 4$, $\al+\bb=2^\ast$, $\al\ge 2$, $\bb\ge 2$ and $\la_1, \la_2\in (0, \La_N)$. Fix any $\nu>0$.
Let $(u, v)$ be any a positive solution of (\ref{eq0-0}). For $\la<0$ we consider the reflection
$$x=(x_1, x_2, \cdots, x_N)\mapsto x^\la=(2\la-x_1, x_2,\cdots, x_N),$$
where $x\in \Sg^\la:=\{x\in\RN : x_1<\la\}$. Define $u^\la(x):=u(x^\la)$ and $v^\la(x):=v(x^\la)$, then
$$u(x)=u^\la(x),\,\,\,v(x)=v^\la(x),\quad\hbox{for}\,\,x\in\partial\Sg^\la=\{x\in\RN : x_1=\la\}.$$
Define $w^\la(x):=u^\la(x)-u(x)$ and $\sg^\la(x):=v^\la(x)-v(x)$ for $x\in\Sg^\la$. Then
\be\label{eq5-0}w^\la(x)=\sg^\la(x)=0,\quad\forall\, x\in\partial\Sg^\la.\ee
Recall that $(u, v)$ satisfies (\ref{eq0-0}), we have
{\allowdisplaybreaks
\begin{align}\label{eq5-1}
-\Delta w^\la(x)=&\frac{\la_1}{|x|^2}w^\la(x)+a_1^\la(x)w^\la(x)+a_2^\la(x)\sg^\la(x)
+\la_1\left(\frac{1}{|x^\la|^2}-\frac{1}{|x|^2}\right)u^{\la}(x)\nonumber\\
\ge &\frac{\la_1}{|x|^2}w^\la(x)+a_1^\la(x)w^\la(x)+a_2^\la(x)\sg^\la(x)
\end{align}
}%
holds in $\Sg^\la\setminus\{0^\la\}$, where
{\allowdisplaybreaks
\begin{align}\label{eq5-2}&a_1^\la:=\frac{(u^\la)^{2^\ast-1}-u^{2^\ast-1}}{u^\la-u}+\nu\al v^\bb \frac{(u^\la)^{\al-1}-u^{\al-1}}{u^\la-u}\ge0,\\
\label{eq5-2-1}&a_2^\la:=\nu\al (u^\la)^{\al-1} \frac{(v^\la)^{\bb}-v^{\bb}}{v^\la-v}\ge 0.
\end{align}
}%
Similarly,
\begin{align}\label{eq5-3}
-\Delta \sg^\la(x)\ge\frac{\la_2}{|x|^2}\sg^\la(x)+b_1^\la(x)\sg^\la(x)+b_2^\la(x)w^\la(x)
\end{align}
holds in $\Sg^\la\setminus\{0^\la\}$, where
{\allowdisplaybreaks
\begin{align}\label{eq5-4}&b_1^\la:=\frac{(v^\la)^{2^\ast-1}-v^{2^\ast-1}}{v^\la-v}+\nu\bb u^\al \frac{(v^\la)^{\bb-1}-v^{\bb-1}}{v^\la-v}\ge0,\\
\label{eq5-4-1}&b_2^\la:=\nu\bb (v^\la)^{\bb-1} \frac{(u^\la)^{\al}-u^{\al}}{u^\la-u}\ge 0.
\end{align}
}%
Define
$$\Om_1^\la:=\left\{x\in\Sg^\la : w^\la(x)<0\right\},\quad \Om_2^\la:=\left\{x\in\Sg^\la : \sg^\la(x)<0\right\}.$$
Since $u, v\in L^{2^\ast}(\RN)$ and $\Om_i^\la\subset \Sg^\la$, there exists $\la_0<0$ such that for any $\la\le \la_0$, we have
{\allowdisplaybreaks
\begin{align}
&\label{eq5-5}\left\|(2^\ast-1)u^{2^\ast-2}+\nu\al(\al-1)u^{\al-2}v^\bb\right\|_{L^{N/2}(\Om_1^\la)}\le\frac{1}{4}\left(1-\frac{\la_1}{\La_N}\right)S,\\
&\label{eq5-6}\left\|(2^\ast-1)v^{2^\ast-2}+\nu\bb(\bb-1)u^{\al}v^{\bb-2}\right\|_{L^{N/2}(\Om_2^\la)}\le\frac{1}{4}\left(1-\frac{\la_2}{\La_N}\right)S,\\
&\label{eq5-7}(\nu\al\bb)^2\left\|u^{\al-1}v^{\bb-1}\right\|^2_{L^{N/2}(\Om_1^\la\cap\Om_2^\la)}
\le\frac{1}{16}\left(1-\frac{\la_1}{\La_N}\right)\left(1-\frac{\la_2}{\La_N}\right)S^2.
\end{align}
}%

\noindent{\bf Step 1.} We claim that for any $\la\le\la_0$, both $w^\la>0$ and $\sg^\la>0$ in $\Sg^\la\setminus\{0^\la\}$.

Fix any $\la\le \la_0$. Define $w^\la_-:=\max\{-w^\la, 0\}$ and $\sg^\la_-:=\max\{-\sg^\la, 0\}$, then $w^\la_-, \sg_-^\la\in D^{1, 2}(\RN)$.
Testing (\ref{eq5-1}) with $w^\la_-$ and using H\"{o}lder inequality and Hardy inequality (\ref{Hardy}), we obtain
\begin{align*}
\int_{\Om^\la_1}|\nabla w^\la_-|^2 &\le\int_{\Om^\la_1}\frac{\la_1}{|x|^2}|w^\la_-|^2+\int_{\Om^\la_1}a_1^\la|w^\la_-|^2
+\int_{\Om^\la_1\cap \Om^\la_2}a_2^\la w^\la_-\sg^\la_-\\
&\le\frac{\la_1}{\La_N}\int_{\Om^\la_1}|\nabla w^\la_-|^2+\|a_1^\la\|_{L^{\frac{N}{2}}(\Om_1^\la)}\|w^\la_-\|^2_{L^{2^\ast}(\Om_1^\la)}\\
&\quad+\|a_2^\la\|_{L^{\frac{N}{2}}(\Om_1^\la\cap\Om_2^\la)}\|w^\la_-\|_{L^{2^\ast}(\Om^\la_1)}\|\sg^\la_-\|_{L^{2^\ast}( \Om^\la_2)}.
\end{align*}
When $\theta\ge 1$, we see from the mean value theorem that
$$\frac{s^\theta-t^\theta}{s-t}\le \theta t^{\theta-1},\quad\forall\,0<s<t.$$
Recall that $u^\la<u$ in $\Om_1^\la$ and $v^\la<v$ in $\Om_2^\la$. Since $\al\ge 2$ and $\bb\ge 2$,
we see from (\ref{eq5-2})-(\ref{eq5-2-1}) and (\ref{eq5-4})-(\ref{eq5-4-1}) that
{\allowdisplaybreaks
\begin{align}
&a_1^{\la}\le(2^\ast-1)u^{2^\ast-2}+\nu\al(\al-1)u^{\al-2}v^\bb,\quad\hbox{in}\,\,\Om_1^\la,\nonumber\\
\label{eq5-0-0}&a_2^\la,\,\, b_2^\la\le\nu\al\bb u^{\al-1}v^{\bb-1},\quad\hbox{in}\,\,\Om_1^\la\cap\Om_2^\la,\\
& b_1^\la\le (2^\ast-1)v^{2^\ast-2}+\nu\bb(\bb-1)u^{\al}v^{\bb-2}\quad\hbox{in}\,\,\Om_2^\la.\nonumber
\end{align}
}%
Then we see from (\ref{eq5-5}) and (\ref{sobolev}) that
$$\|a_1^\la\|_{L^{\frac{N}{2}}(\Om_1^\la)}\|w^\la_-\|^2_{L^{2^\ast}(\Om_1^\la)}\le \frac{1}{4}\left(1-\frac{\la_1}{\La_N}\right)\int_{\Om^\la_1}|\nabla w^\la_-|^2.$$
From above we obtain
{\allowdisplaybreaks
\begin{align}\label{eq5-10}
\frac{3}{4}\left(1-\frac{\la_1}{\La_N}\right)\int_{\Om^\la_1}|\nabla w^\la_-|^2 \le
\|a_2^\la\|_{L^{\frac{N}{2}}(\Om_1^\la\cap\Om_2^\la)} \|w^\la_-\|_{L^{2^\ast}(\Om^\la_1)}\|\sg^\la_-\|_{L^{2^\ast}(\Om^\la_2)}.
\end{align}
}%
Similarly, testing (\ref{eq5-3}) with $\sg^\la_-$ we can prove that
{\allowdisplaybreaks
\begin{align}\label{eq5-11}
\frac{3}{4}\left(1-\frac{\la_2}{\La_N}\right)\int_{\Om^\la_2}|\nabla\sg^\la_-|^2 \le
\|b_2^\la\|_{L^{\frac{N}{2}}(\Om_1^\la\cap\Om_2^\la)} \|w^\la_-\|_{L^{2^\ast}(\Om^\la_1)}\|\sg^\la_-\|_{L^{2^\ast}(\Om^\la_2)}.
\end{align}
}%
Let $|\Om|$ denotes the Lebesgue measure of $\Om$ in $\RN$. If $\left|\Om^\la_1\cap \Om^\la_2\right|>0$, then
{\allowdisplaybreaks
\begin{align}
&\label{eq5-14}\int_{\Om^\la_1}|\nabla w_-^\la|^2>0,\,\,\,
 \int_{\Om^\la_2}|\nabla \sg_-^\la|^2>0.
\end{align}
}%
Combining this with (\ref{eq5-0-0})-(\ref{eq5-11}), (\ref{sobolev}) and (\ref{eq5-7}), we get
{\allowdisplaybreaks
\begin{align*}
0&<\frac{9}{16}\left(1-\frac{\la_1}{\La_N}\right)\left(1-\frac{\la_2}{\La_N}\right)\int_{\Om^\la_1}|\nabla w^\la_-|^2\int_{\Om^\la_2}|\nabla\sg^\la_-|^2\\
&\le \|a_2^\la\|_{L^{\frac{N}{2}}(\Om_1^\la\cap\Om_2^\la)}\|b_2^\la\|_{L^{\frac{N}{2}}(\Om_1^\la\cap\Om_2^\la)}\|w^\la_-\|^2_{L^{2^\ast}(\Om^\la_1 )}\|\sg^\la_-\|^2_{L^{2^\ast}(\Om^\la_2)}\\
&\le \frac{1}{16}\left(1-\frac{\la_1}{\La_N}\right)\left(1-\frac{\la_2}{\La_N}\right)\int_{\Om^\la_1}|\nabla w^\la_-|^2\int_{\Om^\la_2}|\nabla\sg^\la_-|^2,
\end{align*}
}%
a contradiction. Hence $\left|\Om^\la_1\cap \Om^\la_2\right|=0$, and so (\ref{eq5-10})-(\ref{eq5-11}) yield
{\allowdisplaybreaks
\begin{align}
\label{eq5-16}\int_{\Om^\la_1}|\nabla w^\la_-|^2 \le 0,\quad
\int_{\Om^\la_2}|\nabla\sg^\la_-|^2 \le 0.
\end{align}
}%
This implies that $\left|\Om^\la_1\right|=0$ and $\left|\Om^\la_2\right|=0$. That is, both $w^\la\ge 0$ and $\sg^\la\ge 0$ in $\Sg^\la\setminus\{0^\la\}$. If $w^\la\equiv 0$ in $\Sg^\la\setminus\{0^\la\}$, then we see from (\ref{eq5-1}) that
$$-\Delta w^\la(x)\ge \la_1\left(\frac{1}{|x^\la|^2}-\frac{1}{|x|^2}\right)u^{\la}(x)>0\,\,\,\hbox{in}\,\,\Sg^\la\setminus\{0^\la\},$$
a contradiction. So $w^\la\not\equiv 0$ in $\Sg^\la\setminus\{0^\la\}$. Then by the maximum principle,
we conclude that $w^\la> 0$ in $\Sg^\la\setminus\{0^\la\}$. Similarly, $\sg^\la> 0$ in $\Sg^\la\setminus\{0^\la\}$. \\

\noindent{\bf Step 2.} Define $\la^*=\sup\big\{\overline{\la}<0 : w^\la>0, \sg^\la>0 \,\,\hbox{in}\,\,\Sg^\la\setminus\{0^\la\}, \,\,\forall\,\la<\overline{\la}\big\}$. Then we claim that $\la^*=0$.

Assume by contradiction that $\la^*<0$. Clearly we have both $w^{\la^*}\ge 0$ and $\sg^{\la^*}\ge 0$ in $\Sg^{\la^*}\setminus\{0^{\la^*}\}$.
By a similar argument as in Step 1, in fact we have both $w^{\la^*}> 0$ and $\sg^{\la^*}> 0$ in $\Sg^{\la^*}\setminus\{0^{\la^*}\}$.
Take $\e>0$ such that
$$\e<\frac{1}{2}\min\left\{\frac{1}{4}\left(1-\frac{\la_1}{\La_N}\right)S, \,\,\frac{1}{4}\left(1-\frac{\la_2}{\La_N}\right)S,\,\, \frac{1}{16}\left(1-\frac{\la_1}{\La_N}\right)\left(1-\frac{\la_2}{\La_N}\right)S^2 \right\}.$$
Then there exists small $\dd_1\in (0, |\la^*|)$ such that for any $\la\in [\la^*, \la^*+\dd_1)$, there hold
{\allowdisplaybreaks
\begin{align}\label{eq5-18}
&\left\|(2^\ast-1)u^{2^\ast-2}+\nu\al(\al-1)u^{\al-2}v^\bb\right\|_{L^{N/2}(\Sg^\la\setminus\Sg^{\la^*})}\le\e,\\
&\left\|(2^\ast-1)v^{2^\ast-2}+\nu\bb(\bb-1)u^{\al}v^{\bb-2}\right\|_{L^{N/2}(\Sg^\la\setminus\Sg^{\la^*})}\le\e,\\
&\label{eq5-19}(\nu\al\bb)^2\left\|u^{\al-1}v^{\bb-1}\right\|^2_{L^{N/2}(\Sg^\la\setminus\Sg^{\la^*})}\le\e.
\end{align}
}%
Meanwhile, since $w^{\la^*}> 0$ and $\sg^{\la^*}> 0$ in $\Sg^{\la^*}\setminus\{0^{\la^*}\}$, by convergence almost everywhere and thereby in the measure sense
of $(w^\la, v^\la)\to (w^{\la^*}, v^{\la^*})$ in $\Sg^{\la^*}$,
we have that
{\allowdisplaybreaks
\begin{align*}
&\lim_{\la\to\la^*}\left\|(2^\ast-1)u^{2^\ast-2}+\nu\al(\al-1)u^{\al-2}v^\bb\right\|_{L^{N/2}(\Om_1^\la\cap\Sg^{\la^*})}=0,\\
&\lim_{\la\to\la^*}\left\|(2^\ast-1)v^{2^\ast-2}+\nu\bb(\bb-1)u^{\al}v^{\bb-2}\right\|_{L^{N/2}(\Om_2^\la\cap\Sg^{\la^*})}=0,\\
&\lim_{\la\to\la^*}(\nu\al\bb)^2\left\|u^{\al-1}v^{\bb-1}\right\|^2_{L^{N/2}(\Om_1^\la\cap\Om_2^\la\cap\Sg^{\la^*})}=0.
\end{align*}
}%
Then there exists $\dd_2\in (0, \dd_1)$ such that for any $\la\in [\la^*, \la^*+\dd_2)$, there hold
{\allowdisplaybreaks
\begin{align*}
&\left\|(2^\ast-1)u^{2^\ast-2}+\nu\al(\al-1)u^{\al-2}v^\bb\right\|_{L^{N/2}(\Om_1^\la\cap\Sg^{\la^*})}\le\e,\\
&\left\|(2^\ast-1)v^{2^\ast-2}+\nu\bb(\bb-1)u^{\al}v^{\bb-2}\right\|_{L^{N/2}(\Om_2^\la\cap\Sg^{\la^*})}\le\e,\\
&(\nu\al\bb)^2\left\|u^{\al-1}v^{\bb-1}\right\|^2_{L^{N/2}(\Om_1^\la\cap\Om_2^\la\cap\Sg^{\la^*})}\le\e.
\end{align*}
}%
Recall that $\Om_i^\la\subset \Sg^\la$. Combining these with (\ref{eq5-18})-(\ref{eq5-19}), we see that (\ref{eq5-5})-(\ref{eq5-7}) hold for any $\la\in [\la^*, \la^*+\dd_2)$.
Then repeating the proof of Step 1, we conclude that for any $\la\in [\la^*, \la^*+\dd_2)$, $w^{\la}> 0$ and $\sg^{\la}> 0$
in $\Sg^{\la}\setminus\{0^{\la}\}$, which contradicts with the definition of $\la^*$. Therefore $\la^*=0$.\\

\noindent{\bf Step 3.} We claim that both $u$ and $v$ are radially symmetric with respect to the origin.

With the help of Steps 1 and 2, this argument is standard. Since $\la^*=0$, then we can carry out the above procedure in the opposite direction,
namely moving the parallel planes in the negative $x_1$ direction from positive infinity. Then they must stop
at the origin again, and so we get the symmetry of both $u$ and $v$ with respect to $0$ in the $x_1$ direction by combining the
two inequalities obtained in the two opposite directions. Since the direction can be
chosen arbitrarily, we conclude that both $u$ and $v$ are radially symmetric with respect to the origin.
This completes the proof of Theorem \ref{th4}.\hfill$\square$

\vskip0.1in

\s{Uniqueness results for the special case $\la_1=\la_2$ and $\al=\bb=\frac{2^\ast}{2}$}
\renewcommand{\theequation}{7.\arabic{equation}}

When $\la_1=\la_2$,
some uniqueness results about ground state solutions of (\ref{CZ}) were obtained by the authors in \cite{CZ1, CZ2}.
We remark that, by using the same ideas as in \cite{CZ1, CZ2}, these results also hold for problem (\ref{eq0-0}) if we assume  $N\ge 4$, $\la_1=\la_2$ and $\al=\bb=\frac{2^\ast}{2}$.
First we consider the case $N=4$, then we have the following result, which improves Theorem \ref{th2} in case $\la_1=\la_2$.

\bt\label{th5} Assume that $N= 4$, $\la_1=\la_2\in (0, 1)$, $\al=\bb=2$ and $\nu>0$.

\begin{itemize}

\item[(1)] If $\nu\neq1/2$, then for any $\mu>0$, $((1+2\nu)^{-1/2} z_\mu^1,\,  (1+2\nu)^{-1/2} z_\mu^1)$ is a ground state solution  of (\ref{eq0-0}), with
\be\label{s5eq2-28}c_\nu=J_\nu \left((1+2\nu)^{-1/2}z_\mu^1,  (1+2\nu)^{-1/2} z_\mu^1\right)=\frac{1}{2(1+2\nu)}S(\la_1)^{2}.\ee
Moreover, the set $\{((1+2\nu)^{-1/2} z_\mu^1,\,  (1+2\nu)^{-1/2} z_\mu^1): \mu>0\}$ contains all positive ground state solutions of (\ref{eq0-0}).

\item[(2)] If $\nu=1/2$, then for any $\mu>0$ and $\theta\in (0, \pi/2)$, $(\sin \theta\, z_\mu^1, \cos\theta\, z_\mu^1)$ is a ground state solution  of (\ref{eq0-0}) and $c_{1/2}=\frac{1}{4}S(\la_1)^2$. Moreover, the set $\{(\sin\theta\, z_\mu^1,\, \cos\theta \,z_\mu^1): \mu>0, \theta\in (0, \pi/2)\}$ contains all positive ground state solutions of (\ref{eq0-0}).

\end{itemize}
\et

\noindent{\bf Proof.} (1) This result can be obtained by repeating the proofs of \cite[Theorem 1.1 and Theorem 1.2]{CZ1} with trivial modifications. We omit the details.

(2) This result can be obtained by repeating the proofs of Theorem \ref{th3}-(2) with trivial modifications. We omit the details.\hfill$\square$

\br As pointed out in the Introduction, by \cite{Terracini} we know that $Z_i$ contains all positive solutions of (\ref{eq1-2}). Here,
for the case where $N= 4$, $\la_1=\la_2\in (0, 1)$, $\al=\bb=2$, $\nu>0$ and $\nu\neq 1/2$,
we conjecture that the set $\{((1+2\nu)^{-1/2} z_\mu^1,\,  (1+2\nu)^{-1/2} z_\mu^1): \mu>0\}$ contains all positive solutions of (\ref{eq0-0}).\er

Now we consider the case $N\ge 5$. Denote $p=\frac{2^\ast}{2}$ for simplicity. Consider
\be\label{s5eq8}
\begin{cases} k^{p-1}+p\nu k^{\frac{p}{2}-1}l^{\frac{p}{2}}  =
1,\\
p\nu k^{\frac{p}{2}} l^{\frac{p}{2}-1}+ l^{p-1}  =1,\\
k>0,\,\,\, l>0.\end{cases}\ee
Let $\nu>0$. By a direct computation, it was proved in \cite[Lemma 2.1]{CZ2} that there exists $(k_0, l_0)$, such that
\be\label{eqqqq} \hbox{$(k_0, l_0)$ satisfies (\ref{s5eq8}) and}\,\,\,   k_0=\min\{k\,: \,\hbox{$(k, l)$ is a solution of (\ref{s5eq8})}\}.\ee
Then we have the following uniqueness result.

\bt\label{th6}Assume that $N\ge 5$, $\la_1=\la_2\in (0, \La_N)$ and $\al=\bb=p=\frac{2^\ast}{2}$. If
 $\nu\ge\frac{2}{N}$, then for any $\mu>0$, $(\sqrt{k_0}z_\mu^1, \sqrt{l_0}z_\mu^1)$ is a positive ground state solution of (\ref{eq0-0}). Moreover,
the set $\{(\sqrt{k_0}z_\mu^1, \sqrt{l_0}z_\mu^1) : \mu>0\}$ contains all positive ground state solutions of (\ref{eq0-0}).
\et

\noindent{\bf Proof.} This result can be obtained by repeating the proofs of \cite[Theorem 1.1 and Theorem 1.2]{CZ2} with trivial modifications. We omit the details.\hfill$\square$

\vskip0.1in

\noindent{\it Acknowledgements.} The authors wish to thank the anonymous referee very much for  his/her careful reading and valuable comments.

\end{document}